\documentclass[leqno, a4paper,11pt]{amsart}
\pdfoutput=1

\usepackage{enumitem}
\usepackage{amsmath,mathtools,mathrsfs,amsfonts,amssymb}
\usepackage[margin=2.4cm]{geometry}
\usepackage{newpxtext}
\usepackage{tikz-cd}
\usepackage{comment}
\usepackage[notextcomp]{stix}

 \mathchardef\ordinarycolon\mathcode`\:
\mathcode`\:=\string"8000
\begingroup \catcode`\:=\active
  \gdef:{\mathrel{\mathop\ordinarycolon}}
\endgroup

\usepackage{graphicx}
\usepackage{multirow}

\usepackage{hyperref}
\hypersetup{colorlinks,
 linkcolor={blue!90!black},
 citecolor={red!80!black},
 urlcolor={blue!50!black},
 breaklinks=true}
 \usepackage{color}

\usepackage{setspace}
\usepackage{cleveref}

\newtheorem{thm}{Theorem}[section]
\newtheorem{corollary}[thm]{Corollary}
\newtheorem{prop}[thm]{Proposition}
\newtheorem{lemma}[thm]{Lemma}
\newtheorem{fact}[thm]{Fact}

\newtheorem{conj}{Conjecture}

\theoremstyle{definition}
\newtheorem{defn}[thm]{Definition}
\newtheorem{example}[thm]{Example}

\theoremstyle{remark}
\newtheorem{remark}[thm]{Remark}

\newcommand{\bt}{\begin{thm}}
\newcommand{\et}{\end{thm}}
\newcommand{\bp}{\begin{prop}}
\newcommand{\ep}{\end{prop}}
\newcommand{\bd}{\begin{defn}}
\newcommand{\ed}{\end{defn}}
\newcommand{\bl}{\begin{lemma}}
\newcommand{\el}{\end{lemma}}
\newcommand{\bfa}{\begin{fact}}
\newcommand{\efa}{\end{fact}}
\newcommand{\bc}{\begin{corollary}}
\newcommand{\ec}{\end{corollary}}
\newcommand{\bex}{\begin{example}}
\newcommand{\eex}{\end{example}}
\newcommand{\br}{\begin{remark}}
\newcommand{\er}{\end{remark}}
\newcommand{\ben}{\begin{enumerate}}
\newcommand{\een}{\end{enumerate}}

\newcommand{\ds}{\displaystyle}
\newcommand{\mf}{\mathbf}
\newcommand{\floor}[1]{\lfloor #1 \rfloor}
\newcommand{\ceil}[1]{\lceil #1 \rceil}

\newcommand{\scrM}{\mathscr{M}}
\newcommand{\scrR}{\mathscr{R}}
\newcommand{\scrD}{\mathscr{D}}
\newcommand{\scrH}{\mathscr{H}}
\newcommand{\scrK}{\mathscr{K}}
\newcommand{\scrhn}{\mathscr{H}_n}
\newcommand{\scrkn}{\mathscr{K}_n}

\newcommand{\ra}{\rightarrow}

\newcommand{\Z}{\mathbb{Z}}

\newcommand{\N}{\mathbb{N}}
\newcommand{\mscr}{\mathscr{M}}
\newcommand{\rscr}{\mathscr{R}}
\newcommand{\pscr}{\mathscr{P}}
\newcommand{\ptscr}{\widetilde{\mathscr{P}}}
\newcommand{\dscr}{\mathscr{D}}
\newcommand{\cac}{{\mathscr C}}

\newcommand{\cab}{{\mathscr B}}

\newcommand{\hscr}{{\mathscr H}}
\newcommand{\kscr}{{\mathscr K}}

\newcommand{\w}{\mathrm{wt}}
\newcommand{\wtr}{\w_{\rscr}}
\newcommand{\wtp}{\w_{\pscr}}
\newcommand{\wtm}{\w_{\mscr}}
\newcommand{\wtd}{\w_{\dscr}}
\newcommand{\ttr}{t_{\rscr}}
\newcommand{\ttm}{t_{\mscr}}
\newcommand{\ttd}{t_{\dscr}}

\newcommand{\atl}{{\tilde{a}}}
\newcommand{\btl}{{\tilde{b}}}
\newcommand{\ctl}{{\tilde{c}}}
\newcommand{\rw}{\rscr_w}
\newcommand{\mw}{\mscr_w}
\newcommand{\rtri}{\scrR^\triangle}
\newcommand{\mtri}{\scrM^\triangle}

\newcommand{\rect}{\scrR^\hrectangle}
\newcommand{\mrect}{\scrM^\hrectangle}
\newcommand{\prect}{\pscr^\hrectangle}
\newcommand{\rtt}{\scrR^T}
\newcommand{\mtt}{\scrM^T}
\newcommand{\dtt}{\scrD^T}
\newcommand{\rtriw}{\scrR^\triangle_w}
\newcommand{\mtriw}{\scrM^\triangle_w}
\newcommand{\rectw}{\scrR^\hrectangle_w}

\newcommand{\mrectw}{\scrM^\hrectangle_w}

\newcommand{\bfU}{\mathbf{U}}
\newcommand{\bfe}{\mathbf{e}}
\newcommand{\bfu}{\mathbf{u}}
\newcommand{\bfd}{\mathbf{d}}
\newcommand{\bfv}{\mathbf{v}}

\newcommand{\bfq}{\mathbf{q}}

\newcommand{\beone}{\mathbf{e}_1}
\newcommand{\betwo}{\mathbf{e}_2}

\newcommand{\hphi}{\hat{\varphi}}
\newcommand{\vphi}{\varphi}
\newcommand{\hbfq}{\hat{\bfq}}

\newcommand{\npr}{\nu^{+}_{r}}
\newcommand{\nmr}{\nu^{-}_{r}}

\makeatletter
\newcommand\myitem[1][]{\item[#1]\refstepcounter{dummy}\def\@currentlabel{#1}}
\makeatother

\makeatletter
\@namedef{subjclassname@2020}{
  \textup{2020} Mathematics Subject Classification}
\makeatother

\usepackage{imakeidx}
\makeindex[options=-c -s dotted.ist, title={Notation Index},]

\begin{document}
\author[A. \,Sammartano, E. \,Schlesinger]{Alessio~Sammartano and Enrico~Schlesinger}
\address{Dipartimento di Matematica \\ Politecnico di Milano \\ Milan \\ Italy}
\email{alessio.sammartano@polimi.it}
\email{enrico.schlesinger@polimi.it}

\title{Initial ideals of weighted forms and the  genus of locally Cohen-Macaulay curves}

%\subjclass[2020]{Primary: ; Secondary: }
%\keywords{}

\begin{abstract}
Let $\mathcal{C}$ be a locally Cohen-Macaulay curve in complex projective 3-space.
The maximum genus problem  predicts the largest possible arithmetic genus $g(d,s)$ that $\mathcal{C}$ can achieve assuming that it has degree $d$ and does not lie on surfaces of degree less than $s$.
In this paper, we prove that this prediction is correct when $d=s$ or $d\geq 2s-1$.
We obtain  this result by proving another conjecture, by Beorchia, Lella, and the second author, about initial ideals  associated to certain homogeneous forms in a non-standard graded polynomial ring.
\end{abstract}

\maketitle

\section{Introduction}

The classification of  curves is a central theme in algebraic  geometry, 
which has prompted a tremendous amount of  research in the last two centuries.
An important step in this classification effort is  the determination of  the genera and degrees for which there exist curves in a projective space with(out) prescribed singularities. 
For smooth curves, this  problem  dates back to the nineteenth century
\cite{Halphen}
and is yet not completely solved 
\cite{BBEM,GP,H,HH,MaSc}.
In this work, we consider  curves in $ \mathbb{P}^3_\mathbb{C}$ that are
locally Cohen-Macaulay,
that is, schemes of dimension one without embedded or isolated points;
this is the natural class of curves from the point of view of liaison theory 
\cite{Migliore}.
If $\mathcal{C}\subseteq \mathbb{P}_\mathbb{C}^3$ is a locally Cohen-Macaulay curve of degree $d$
that is not contained in any surface of degree less than $s$,
then $d\geq s$, and  the arithmetic genus of $\mathcal{C}$ is bounded above by the following  function
\begin{equation}\label{Bound}
g(d,s) = \begin{cases}
(s-1)d+1 - {s+2 \choose 3} & \text{ if } s \leq d \leq 2s,\\
{d-s\choose 2} - {s-1\choose 3} & \text{ if } d \geq 2s+1,
\end{cases}
\end{equation}
see \cite{B,S}.
The \emph{maximum genus problem for locally Cohen-Macaulay curves} asks whether this upper bound is always attained.

The case of this problem when $d=s$  is of particular interest, since it  implies all the cases when $d \geq 2s-1$
\cite[Proposition 6.1]{BLS}.
In order to settle the case $d=s$,
a  construction of primitive multiple lines  is proposed in  \cite{BLS}, 
where it is  conjectured that they achieve the  bound \eqref{Bound}, yielding  thus an affirmative answer to the maximum genus problem for $d=s$ and  $d \geq 2s-1$.
In turn, this geometric conjecture is implied by the following algebraic conjecture, which is the main subject of this work.

\begin{conj}[\protect{\cite[Conjecture B]{BLS}}]\label{Conj1}
Let $P = \mathbb{C}[x,y,z]$ be a  polynomial ring  with weights $\w(x)=1, \w(y)=2, \w(z)=3$.
Let $f\in P$ be a general $\w$-homogeneous polynomial with $\w(f)=3m$, for some  $m \in \mathbb{N}_{>0}$.
Then, the ideal $I = (x,y)^{3m-2}+ (g)$ contains no polynomial of standard degree less than $3m-2$.
\end{conj}

An equivalent formulation in terms of Gr\"obner bases is the following:
if $<$ is a term order in $P$ refining the standard grading on $P$ such that $x>y>z$,
then the initial ideal of $I$ is
\begin{equation}
\mathrm{in}_<(I) = (x,y,z)^{3m-2}.
\end{equation}
This formulation bears a strong resemblance to a long-standing conjecture of  Moreno-Socias 
\cite{MS} concerning the structure of initial ideals of general forms in standard graded polynomial rings.
While the latter is known to have an affirmative answer  in dimension three, with proof essentially due to Anick \cite{A},
the non-standard weights represent a considerable complication,
 and Anick's argument does not extend to the weighted context.

Conjecture \ref{Conj1} can be also  translated  in terms of linear algebra. 
Let $R=\mathbb{C}[x,y,z]/(x,y)^{3m-2}$ and $M=(x,y,z)^{3m-2}/(x,y)^{3m-2}$,
then the graded $P$-modules 
$R[-3m]$ and $M$ have the same  Hilbert function.
A $\w$-homogeneous  $f\in P$ with $\w(f)=3m$ 
induces a  vector space map $f \cdot: R[-3m] \ra M$  by multiplication by $f$ followed by the projection that kills all monomials of degree less than $3m-2$. 
We show that Conjecture \ref{Conj1} is equivalent to the following.

\begin{conj}\label{Conj2}
Let $f\in P$ be a general $\w$-homogeneous polynomial of weight $3m$ for some $m \in \mathbb{N}_{>0}$.
The vector space map $ f \cdot :R_{w} \rightarrow M_{w+3m}$
is an isomorphism for all  $w\in \mathbb{N}$.
\end{conj}

In this paper, we prove the following result.

\begin{thm}\label{ThmConjAreTrue}
Conjectures \ref{Conj1} and \ref{Conj2} are true.
As a consequence, 
if either $d=s \geq 1$ or $d \geq 2s+1 \geq 3$, 
the maximum genus of a locally Cohen-Macaulay curve in $\mathbb{P}^3_{\mathbb{C}}$ of degree $d$ that does not lie on a surface of degree $s-1$ is equal to $g(d,s)$.
\end{thm}

We now discuss the general outline of the proof.
Let $\mathbf{A}_w=(a_{\bfv,\bfu})$ be the square matrix representing the  map $f \cdot:R_w \to M_{w+3m}$
with respect to the monomial bases $\rscr_w$ of $R_w$ and $\mscr_w$ of $M_{w+3m}$.
Denote by $a_\bfq$ the coefficient of the monomial $\bfq$ in $f$.
The entries of $\mathbf{A}_w $ are
$a_{\bfv,\bfu}
= a_\frac{\bfv}{\bfu}$ if $\bfu \text{ divides }\bfv$,
$a_{\bfv,\bfu}
= 0 $ otherwise.
Conjecture \ref{Conj2} amounts to the non-vanishing of the determinant of $\mathbf{A}_w$  for all $w \in \mathbb{N}$.
 A natural combinatorial approach to proving this non-vanishing is the search for a  bijection $\varphi : \rscr_w \to \mscr_w$
satisfying  two properties:
\begin{description} \label{divis-uniq}
\item[Divisibility]
 $\bfu$ divides $\varphi(\bfu)$ for every monomial $\bfu \in \rscr_w$.
\item[Uniqueness]
the multiset of monomials $\mathscr{Q}_\varphi = \big\{	\bfq(\bfu)= \frac{\varphi(\bfu)}{\bfu}\mid \bfu \in \rscr_w\big\}$
is different from $\mathscr{Q}_{\hat{\varphi}}$ for every other bijection $\hat{\varphi} : \rscr_w \to \mscr_w$ satisfying {\bf Divisibility}.
\end{description}
The property
 {\bf Divisibility} guarantees that the product $\Gamma_\varphi = \prod_{\bfu \in \rscr_w} a_{\vphi(\bfu),\bfu}$
is non-zero, while
the property {\bf Uniqueness} guarantees that the term $\Gamma_\varphi$ appears with coefficient $\pm 1$ in  $\det(\mathbf{A}_w)$,
implying  the desired non-vanishing  for  general $f$.
We  call the monomials $\bfq(\bfu)= \frac{\varphi(\bfu)}{\bfu}$ in $\mathscr{Q}_\varphi$ the {\em multipliers} of $\vphi$.
Unfortunately, in general there exists no such bijection, see Section \ref{examples}.
Therefore, we need to consider a subtler property  that implies 
the non-vanishing
of $\det(\mathbf{A}_w)$:
\begin{description} \label{noncancellation}
\item[Non cancellation]
for every bijection $\hat{\varphi} : \rscr_w \to \mscr_w$  with the same
 multiset of multipliers as $\varphi$, the permutation $\hphi^{-1} \circ \vphi$ is even.
\end{description}
If $\vphi$ satisfies {\bf Divisibility} and {\bf Non cancellation}, then the term
 $\Gamma_\varphi$ appears with a
coefficient $\pm n$ in $\det(\mathbf{A}_w)$, 
where $n$ is the  number of bijections $\hat{\varphi}$.
Thus, it follows that   $\det(\mathbf{A}_w) \ne 0$
for  general $f$.
Our main technical theorem is the following:

\bt \label{mainthm}
For every  $m,w\in \mathbb{N}$ with $m>0$ there is a bijection $\varphi : \rscr_w \to \mscr_w$ that satisfies both {\bf Divisibility} and {\bf Non cancellation}. 
\et

Theorem \ref{mainthm} implies, more generally,  the statement of Theorem \ref{ThmConjAreTrue} 
over arbitrary fields of characteristic zero.
We remark that the problems considered above are meaningful and open over any infinite field.
However, in positive characteristic it seems unlikely that Conjecture \ref{Conj2} could be attacked with a combinatorial approach similar to that outlined above, 
due to the presence of the numerical coefficient in the term that certifies the non-vanishing of 
$\det(\mathbf{A}_w)$.
As explained above,
the existence of bijections $\varphi$ satisfying both {\bf Divisibility} and {\bf Uniqueness}
would imply the validity of Theorem \ref{ThmConjAreTrue}  in all characteristics,
but such bijections do not exist in general.

The proof of  Theorem \ref{mainthm} is intricate and long.
In some sense, this is to be expected.
Whenever one wants to  verify  that a certain property holds generically, 
one must  prove that the property holds in a Zariski open subset, 
and  that this subset is non-empty. 
The difficulty is typically concentrated in the latter step, 
which requires an explicit construction, and often one that relies on a number of non-natural choices and ad-hoc degenerations.
There are abundant and illustrious examples of this phenomenon, including
Fr\"oberg's conjecture \cite{A,MS},
the Maximal rank conjecture \cite{BE,Larson},
the Lefschetz properties \cite{Adiprasito,Lefschetz},
and the Minimal resolution conjecture \cite{FL,HS},
just to name a few.
We have done our best effort to clarify the nature of the arguments 
and to improve the readability of the proofs, by writing all the details explicitly  and recalling the technical definitions each time they are used.
We have included an appendix to illustrate all the aspects of the construction of the bijection $\varphi$ of Theorem \ref{mainthm}.
We have also collected all the notation in an index at the end of the paper.

\subsection*{Structure of the paper}
In Section \ref{weights} we collect some preliminary results and  we exhibit an explicit isomorphism of graded vector spaces $\psi:R[-3m] \rightarrow M$,
cf. Proposition \ref{Mstructure}. 
In Section \ref{reduction}
we reduce the size of the problem. 
In Section \ref{varphi} we subdivide the monomial bases $\rscr=\rect \cup \rtri$ and $\mscr=\mrect \cup \mtri$
as the disjoint union of a {\em rectangular} and a {\em triangular region}. 
We then construct a bijection $\varphi^\hrectangle: \rect \ra \mrect$ that satisfies {\bf Divisibility} and {\bf Uniqueness}. 
In Section \ref{triangular} we construct a bijection $\varphi^\triangle  : \rtri \ra \mtri$ that satisfies {\bf Divisibility} and {\bf Uniqueness}.
Gluing $\varphi^\hrectangle$ and $\varphi^\triangle$ we obtain, for all weights $w$, a bijection
$\vphi_w: \rscr_w \ra \mscr_w$ satisfying {\bf Divisibility}. 
In Section \ref{wcongruozero} we show that $\vphi_w$ satisfies
{\bf Uniqueness} when $w\equiv 0 \pmod{3}$. 
In Section \ref{examples} we show that {\bf Uniqueness} fails for $\phi_w$
already when $m=3$ and $w=8$; we also show that, when $m=4$ and $w=14$, there is no bijection $\rscr_w \ra \mscr_w$ that satisfies both
{\bf Divisibility} and {\bf Uniqueness}. 
In Section \ref{spblocks} we analyze the failure of uniqueness on certain subdomains which we call {\em special blocks},
and we show that {\bf Non cancellation} still holds for these.
In Section \ref{wcongruounodue} we show that, when $w\equiv 1,2 \pmod{3}$,  uniqueness fails for $\vphi_w$ only along special blocks,  allowing us to prove that $\vphi_w$ satisfies {\bf Non cancellation} and completing the proof of Theorem \ref{mainthm}.
Appendix \ref{AppendixTable} contains six tables which illustrate 
the details of the construction of $\varphi$
for $m=7$
and   $18 \leq w \leq 23$, 
so as to include one example for each congruence class of the weight $w$ modulo 6.

\subsection*{Notation} \label{notation} The reader can find an index of notation at the end of the paper.
Given a real number $x$, we denote by $\floor{x}$ and $\ceil{x}$ respectively the largest integer
$\leq x$ and the smallest integer $\geq x$.
For an integer $n$, we 
define $\epsilon(n)$ \index{aaep @ $\epsilon(n)=n-2\floor{\frac{n}{2}}$}
 as the unique integer $\epsilon \in \{0,1\}$ congruent to $n$ modulo $2$,
 so that
\begin{equation}\label{residuemod2def}
n=2\floor{\frac{n}{2}}+\epsilon(n).
\end{equation}
Similarly, $\eta(n)$\index{aaet @ $\eta(n)=n-3\floor{\frac{n}{3}}$} will denote the unique integer $\eta \in \{0,1,2\}$ congruent to $n$ modulo $3$, so that
for every integer $n$
\begin{equation}\label{residuemod3def}
n=3\floor{\frac{n}{3}}+\eta(n).
\end{equation}
The residues of $n$ modulo $2$ and modulo $3$ determine the residue modulo $6$, precisely,
\begin{equation}\label{residuemod6def}
n=6\floor{\frac{n}{6}}+\eta(n)+3 \epsilon(n+\eta(n)).
\end{equation}

\section{Preliminaries} \label{weights}
Let $\Bbbk$ be a field of characteristic zero, and
let $P=\Bbbk[x,y,z]$ \index{P @ Polynomial ring $P=\Bbbk[x,y,z]$}be the polynomial ring in $3$ variables. 
We write $\bfu=x^ay^bz^c$ for monomials in $P$ and, 
when convenient, 
we identify the monomial $\bfu$ with the integral vector $(a,b,c) \in \N^3$. 
The weight of such a monomial
is $\w(\bfu)=a+2b+3c$. 
Let $P_w$ denote the linear span of monomials of weight $w$, that is, the set of $\w$-homogeneous polynomials of weight
$w$.

Fix an integer $m \geq 2$. 
Let $R$ denote the graded $\Bbbk$-algebra $\Bbbk[x,y,z]/(x,y)^{3m - 2}$\index{R @ $R=k[x,y,z]/(x,y)^{3m - 2}$}. 
We will still denote by $\bfu$ the class of a monomial $\bfu$ in $R$.  
The monomial $\Bbbk$-basis $\rscr$\index{ms @ $\mscr$ the monomial basis of $M$} of $R$ consists of  monomials $\mf{u}=x^a y^b z^c$ for which $a+b \leq 3(m - 1)$.

Let $M$ denote the graded ideal of $R$ defined as
$
M = (x,y,z)^{3m - 2}/(x,y)^{3m - 2} \index{M @ $M = (x,y,z)^{3m - 2}/(x,y)^{3m - 2}$}.
$
The monomial $\Bbbk$-basis $\mscr$ of $M$ \index{ms @ $\mscr$ the monomial basis of $M$}consists of  monomials $\bfv=x^i y^j z^k$ that satisfy $i+j\leq 3(m - 1)$  and $i+j+k \geq 3m - 2$. Note that such monomials have weight at least $3m$, 
and only one of them, namely
$\mf{v}_0=x^{3(m - 1)}z$, has weight exactly $3m$. In fact, as we show in the next Proposition, $M$ is isomorphic as a graded vector space to $R[-3m]$, a statement that follows also from the remark that both ideals $(x,y,z)^{3m - 2}$ and
$(x,y)^{3m - 2}+(z^m)$ have the same weighted Hilbert function as an ideal $(g_1,g_2,g_3)$ where $g_i \in P$ are $\w$-homogeneous polynomials of weight $\w (g_i)=3m+1-i$ and form a regular sequence.
\bp \label{Mstructure}
There is an  isomorphism of graded vector spaces
$\psi:R[-3m] \rightarrow M$ given on the monomial bases of  $R$ and $M$ by the formula:
\begin{equation}\label{psi}
\psi\left(x^a y^b z^c\right)=
x^{3(m-1)-a-b}y^{a}z^{1+b+c}   \qquad (a+b \leq 3(m - 1)).
\end{equation}
In particular, $
\dim_\Bbbk \, R_w \, = \, \dim_\Bbbk \, M_{w+3m}$  for all $w \in \N$.

\ep

\begin{proof}
We define a different graded $P$-module structure on $M = (x,y,z)^{3m - 2}/(x,y)^{3m - 2}$, 
with respect to which $M$ is isomorphic  to $R[-3m]$ with generator $\mf{v}_0=[x^{3(m-1)}z]$. 
To avoid confusion, we use capital letters
$X$, $Y$ and $Z$ to denote the action of $x$, $y$ and $z$ on $M$ with respect to this new module structure.
Let $X$ act on $M$
as multiplication by $x^{-1}y$, 
$Y$ as multiplication by $x^{-1}z$ and $Z$ as
multiplication by $z$. 
This gives the formula for $\psi$ in the statement: note
that $\psi$ is well defined because the class in $R$ of the monomial $X^a Y^b
Z^c$ is non zero if and only if $a+b \leq 3(m - 1)$; 
it is $\w$-homogeneous of weight
$3m$, that is,  $\w (\psi (x^a y^b z^c))= \w( x^a y^b z^c) + 3m$; 
 it is bijective, with inverse given on the monomial bases by
\begin{equation}\label{psiinverse}
\psi^{-1}([x^i y^j z^k])= x^{j}y^{3(m-1)-i-j}z^{i+j+k -3m - 2}   \qquad (i+j\leq 3(m - 1) \text{~and~}i+j+k\geq 3m - 2).
\end{equation}
\end{proof}

\bc \label{same-hilb}
Fix an integer $m \geq 2$.
Let
\[
g(x,y,z)= a_0 z^m+a_1(x,y)z^{m\!-\!1}+\cdots + a_{m\!-\!1} (x,y)z + a_m(x,y)
\]
be a $\w$-homogeneous polynomial in $\Bbbk[x,y,z]$ of weight $3m$, with $a_0\ne 0$. 
Consider in $\Bbbk[x,y,z]$ the $\w$-homogeneous ideals
$I_g=(x,y)^{3m - 2}+(g)$\index{Ig @ $I_g=(x,y)^{3m - 2}+(g)$}
 and $J=(x,y,z)^{3m - 2}$\index{J @ $J=(x,y,z)^{3m - 2}$}.
 Then, $\Bbbk[x,y,z]/I_g$ and $\Bbbk[x,y,z]/J$
are isomorphic as graded $\Bbbk$-vector spaces. \ec
\begin{proof}
Since $g$ is $\w$-homogeneous of weight $3m$ and is not a zero divisor (because
$a_0\neq 0$) of $R$, the ideal $I_g/(x,y)^{3m - 2}$ generated by $g$ is isomorphic
to $R[-3m]$. 
Thus, $\Bbbk[x,y,z]/I \cong R/R[-3m]$ as $\Z$-graded $\Bbbk$-vector spaces.
By Proposition \ref{Mstructure},  the same is true for
$\Bbbk[x,y,z]/J \cong R/M$.
\end{proof}

\noindent
It follows that Conjectures \ref{Conj1} can be restated in terms of linear algebra.
\bp \label{good g} Fix an integer $m \geq 2$. Let $g$ be
a $\w$-homogeneous polynomial in $P=\Bbbk[x,y,z]$ of weight $3m$. The following
conditions are equivalent:
\begin{enumerate}
\item\label{good g:1}
The ideal $I_g=(x,y)^{3m-2}+(g)$ does not contain a nonzero polynomial of standard degree $3(m-1)$.
\item\label{good g:2}
If  $\succcurlyeq$ is a term order on the monomials of $\Bbbk[x,y,z]$ that refines the standard degree,
the initial ideal of $I_g=(x,y)^{3m-2}+(g)$ with respect to $\succcurlyeq$ is $J=(x,y,z)^{3m-2}$.
\item\label{good g:3}
 The $\w$-homogeneous $\Bbbk$-linear map $ g \cdot :R[-3m] \rightarrow M$ defined by multiplication by $g$ followed by the $\Bbbk$-linear projection from $R$
to $M$
is an isomorphism.
\end{enumerate}
\ep

\section{Reduction of  size for the matrix $\mathbf{A}_w$} \label{reduction}
In this section, we reduce the problem to a finite number of values of $w$ and to a smaller dimension of
the spaces involved. 
To this end, we define $I_0 \subset \Bbbk[x,y,z]$ to be the ideal generated as a $\Bbbk$-vector subspace by those monomials
$\mf{u}=x^ay^bz^c$ that satisfy
\begin{equation}\label{I0}
2a+2b+3c > 6m-9
\end{equation}
and $J_0 \subset \Bbbk[x,y,z]$ to be the ideal generated by those monomials
$\mf{v}=x^iy^jz^k$ that satisfy
\begin{equation}\label{J0}
2i+2j+3k > 9(m\!-\!1).
\end{equation}
Note that
\begin{enumerate}
\item
$(x,y)^{3m-2} \subset I_0$ because $a+b \geq 3m-2$ implies $2a+2b+3c \geq 6m-4$;
\item
$J_0 \subset J=(x,y,z)^{3m-2} $ because
the condition $2i+2j+3k > 9(m\!-\!1)$ implies $i+j+k \geq 3m-2$;
\item \label{defineprime}
the set $\rscr'$ of monomials $\mf{u}=x^ay^bz^c$ such that $2a+2b+3c \leq 6m-9$ is a basis of
$R'=\Bbbk[x,y,z]/I_0$\index{rp @ $R'=\Bbbk[x,y,z]/I_0$}
because the  condition $2a+2b+3c > 6m-9$ is closed under multiplication by a monomial;
by the same token,  the set $\mscr'$ of monomials $\mf{v}=x^iy^jz^k$ such that $2i+2j+3k \leq 9(m\!-\!1)$, $i+j \leq 3(m\!-\!1)$ and $i+j+k \geq 3m-2$ is a basis of
$M'=(x,y,z)^{3m-2} /(J_0+(x,y)^{3m-2} )$\index{mp @ $M'=(x,y,z)^{3m-2} /(J_0+(x,y)^{3m-2} )$}.
\end{enumerate}

\bp \label{Upper triangular}
Let $g$ be a $\w$-homogeneous polynomial in $P=\Bbbk[x,y,z]$ of weight $3m$.
Then
\begin{enumerate}
  \item the map $g\cdot$ maps $I_0$ into $J_0$, therefore,
  it induces a $\w$-homogeneous vector space map  $R'[-3m] \ra M'$ which we still denote
  by $g\cdot$;
  \item if the coefficient of $z^m$ in $g$ is nonzero, then $g \cdot:R[-3m] \ra M$ is an isomorphism
  if and only if the induced map $g \cdot: R'[-3m] \ra M'$ is an isomorphism.
  \end{enumerate}
  \ep
\begin{proof}
The first statement follows from the fact that, if $\w(\bfq=x^{a'}y^{b'}z^{c'})=3m$,
then
$$
  2(a+a')+ 2(b+b')+3(c+c') \geq 3m+ 2a+2b+3c.
  $$

The second statement will follow from the Snake's Lemma if we can show that
$$g \cdot: I_0/(x,y)^{3m-2} [-3m] \ra (J_0+(x,y)^{3m-2} )/((x,y)^{3m-2}$$
is an isomorphism.

We first prove this for $g=z^m$. By construction, a class $\bf{f}$ is in the kernel
of $z^m \cdot$ if and only if it can be represented by a polynomial of degree
$<2m-2$. 
The condition $2a+2b+3c \geq 6m-8$ implies $a+b+c \geq 2m-2$, that is,
every monomial in $I_0$ has degree at least $2m-2$. We conclude that $z^m \cdot$ is
injective on $I_0/(x,y)^{3m-2}$. To show that it is surjective,
observe that monomials $\bfv=x^{i}y^{j}z^{k}$ satisfying $i+j \leq 3(m\!-\!1)$ and
$2i+2j+3k \geq 9m-8$ form a basis of $(J_0+(x,y)^{3m - 2})/((x,y)^{3m-2}$.
For such a monomial $\bfv$ we have
$$3k \geq 9m-8-2(i+j)\geq 3m-2,$$
hence $k \geq m$, that is, $\bfq$ is in the image by $z^m\cdot$ of $I_0$.
This proves that $z^m \cdot: I_0/(x,y)^{3m-2}[-3m] \ra (J_0+(x,y)^{3m - 2})/((x,y)^{{3m-2}}$ is an isomorphism.

Now, suppose that the coefficient $a_0$ of $z^m$ in $g$ is nonzero. We fix a monomial order $\prec$ in $\Bbbk[x,y,z]$ such
that  $z^m \prec \mf{u}$ for every other monomial $\mf{u}$ of weight $3m$. Given a weight $w$,
consider a monomial basis $\{\bfu_1,\ldots,\bfu_r\}$ of $ (I_0/(x,y)^{3m-2})_w$ such that $i<j$ implies
$\bfu_i \prec \bfu_j$. A nonzero element of $ (I_0/(x,y)^{3m-2})_w$ has the form
$
\mathbf{f}=b_j \bfu_j + b_{j+1} \bfu_{j+1}+ \cdots b_r \bfu_r
$
with $b_j \neq 0$. Then
$$
g \cdot \mathbf{f}=a_0 b_j z^m\bfu_j + \mbox{higher order terms}
$$
is nonzero because $a_0 \neq 0$ and $\deg (z^m\bfu_j) \geq 3m-2$.
Hence, $(g \cdot)_w$ is injective, and therefore an isomorphism, for every $w$.
\end{proof}

\bc\label{finitecases} Let $g$ be
a $\w$-homogeneous polynomial in $P=\Bbbk[x,y,z]$ of weight $3m$.
Assume $z^m$ appears in $g$ with nonzero coefficient. Then,
in order to show $g \cdot:R[-3m] \rightarrow M$ is an isomorphism, it is
enough to show that
$
(g \cdot)_{w}: R_{w} \rightarrow M_{w+3m}
$
is an isomorphism for $0 \leq w \leq 6m-9$.
\ec

\begin{proof}
If $ w >6m-9$, then $2a+2b+3c=a+w>6m-9$, so $R'_w=0$ and the statement follows from Proposition \ref{Upper triangular}.
\end{proof}

\bc \label{Trimming}
Fix an integer $m \geq 2$ and define $\ell=3m-2$. Let $
g(x,y,z) $ be a $\w$-homogeneous polynomial in $\Bbbk[x,y,z]$ of weight $3m$, in which
$z^m$ appears with a nonzero coefficient. Let $I_g$ denote the $\w$-homogeneous ideal
$I_g=(x,y)^\ell+(g)$. Any monomial $ \bfv=x^{i}y^jz^k $ satisfying
$$
2i+2j+3k > 9(m\!-\!1)
$$
belongs to $I_g$, that is, $J_0 \subset I_g$.
\ec
\begin{proof}
 If $i+j \geq \ell$, then $\bfv$ is in $I$. Suppose $i+j < \ell$. Then, $k \geq m$.
Let $w=\w(\bfv)-3m$. 
By  Proposition \ref{Upper triangular} we have
$
\bfv=g \bfu+\hat{\bfv}
$
where $\bfu \in I_0$ and $\hat{\bfv}$ is a $\w$-homogeneous polynomial of standard degree at most $3(m\!-\!1)$.
But every monomial in $\bfu$ has degree at least $2m-2$, while we may assume that every monomial
in $g$ has degree at least $m$. Hence $\hat{\bfv}$ must be zero, that is, $\bfv$ belongs to $I_g$.
\end{proof}

\section{Construction of the bijection $\varphi$ satisfying divisibility - part one} \label{varphi}
We recall the notation we introduced for the monomial
bases of the graded vector spaces we are interested in. Given an integer $m \geq 2$, we denote by $\rscr$ and $\mscr$ the monomial bases of $R=\Bbbk[x,y,z]/(x,y)^{3m-2}$ and $M=(x,y,z)^{3m-2}/(x,y)^{3m-2}$ respectively. 
For the remainder of the paper, we will identify a monomial
$\bfu=x^ay^bz^c$ with the vector $(a,b,c)$ in $\N^3$.
Thus,
\begin{equation} \label{rscrdef}\index{rs @ $\rscr$ the monomial basis of $R$}
  \rscr=
\left\{\bfu=
(a,b,c) \in \N^3: a+b\leq 3m-3
\right\},
\end{equation}
and
\begin{equation} \label{mscrdef}\index{ms @ $\mscr$ the monomial basis of $M$}
  \mscr=
\left\{\bfv=
(i,j,k) \in \N^3:   i\!+\!j \leq 3(m\!-\!1),\; i\!+\!j\!+\!k \geq 3m\!-\!2
\right\}.
\end{equation}
The monomial basis $\rscr'$ of $R'$ is identified with the subset of $\rscr$ of monomials satisfying the additional condition $2a+2b+3c \leq 6m-9$. 
Since this  implies $a+b \leq 3(m\!-\!1)$, we have
\begin{equation} \label{rprimedef}\index{rsp @ $\rscr'$ the monomial basis of $R'$}
  \rscr'=
\left\{\bfu=
(a,b,c) \in \N^3:  2a+2b+3c \leq 6m-9
\right\}.
\end{equation}
The monomial basis $\mscr'$ of $M'$ is identified with the subset of $\mscr$ of monomials satisfying the additional condition $2i+2j+3k \leq 9(m-1)$. 
Hence,
\begin{equation} \label{mprimedef}\index{msp @ $\mscr'$ the monomial basis of $M'$}
  \mscr'=
\left\{\bfv=
(i,j,k) \in \N^3:  2i+2j+3k \leq 9(m\!-\!1), i\!+\!j \leq 3(m\!-\!1),\; i\!+\!j\!+\!k \geq 3m\!-\!2
\right\}
\end{equation}
Note that  $1 \leq k \leq 3(m\!-\!1)$ and $ i+j \geq 1$ for every $(i,j,k) \in
\mscr'$.

Recall the weight of a monomial $\bfu=(a,b,c)$ in  $\rscr$ is $\wtr (\bfu)=a+2b+3c$, whereas for
a monomial $\bfv=(i,j,k)$ in  $\mscr$ it is convenient to define its $\mscr$-weight as $\wtm (\bfv)=i+2j+3k-3m$. 
Then, by
Propositions \ref{Mstructure} and \ref{Upper triangular},
for any integer $w$ the set $\rscr'_w$ of monomials in $\rscr'$ of weight $w$ is in one-to-one
correspondence with the set $\mscr'_w$ of monomials in $\mscr'$ of $\mscr$-weight $w$. Note that
$\rscr'_w$ and $\mscr'_w$ are nonempty if and only if $0 \leq w \leq 6m-9$.

In this section we construct an explicit  bijection $\varphi: \rscr' \ra \mscr'$ that preserves weights and satisfies divisibility, that is,
such that for every $\bfu=(a,b,c) \in \rscr'$ its image $\bfv=(i,j,k) \in \mscr'$ satisfies
$
i \geq a,\, j \geq b, \, k \geq c
$.
The bijection we construct also preserves a subtler invariant which we now introduce, motivated
by the remarkable fact that for each integer $0 \leq h \leq m-2$, the set of monomials $\bfu=(a,b,c)$ in
$\rscr'$ satisfying $c \leq 2h+1$ is in one-to-one
correspondence with the set of elements of $\mscr'$ satisfying $i+j \geq 3(m-1)-3h -2$.
This suggests to define
for $\bfu=(a,b,c) \in \rscr$ and $\bfv =(i,j,k)\in \mscr$
\begin{equation}\label{tinvariants}\index{tr @ $\ttr$ the $t$ invariant of $R$}\index{trm @ $\ttm$ the $t$ invariant of $M$}
  \ttr(\bfu)= \floor{\ds \frac{c}{2}}; \quad \ttm(\bfv)  = \floor{\ds \frac{3(m\!-\!1)-(i+j)}{3}}.
\end{equation}
A more conceptual introduction of these invariants can be found in the proof of Theorem \ref{phitriregion}.
The fact that, for every integer $t$, the set $\rscr'_{w,t}$ of monomials  $\bfu \in \rscr'_w$ satisfying $\ttr(\bfu)=t$ is in one-to-one correspondence with the  set $\mscr'_{w,t}$ of monomials $\bfv \in \mscr'_w$ satisfying $\ttm(\bfv)=t$ will follow from the construction of the bijection $\varphi$ that  preserves weights and $t$-invariants (and satisfies divisibility).

It will be useful to parametrize monomials $\bfu$ in $\rscr'$ by means of their weight $\wtr(\bfu)$, their $t$-invariant $\ttr(\bfu)$, and their first coordinate $a(\bfu)$.
Recall that, for  an integer $n$, we set
$\epsilon(n)=0$ if $n$ is even and $\epsilon(n)=1$ if $n$ is odd.

\bp \label{rprimeparamwta} \index{ur @ $\bfu_{\rscr}(w,t,a)$ parametrization of $\rscr'$ }
Fix an  $m \geq 2$. Recall that
$\rscr'=
\left\{\bfu=
(a,b,c) \in \N^3:  2a+2b+3c \leq 6m-9 \right\}$, and let
$$\mathscr{P}_{\rscr}=
\left\{(w,t,a) \in \N^3: \;
w+a \leq 6m-9, \; w-a-6t-3 \epsilon(w+a) \geq 0 \right\}.
$$
The function $\bfu_{\rscr}:\mathscr{P}_{\rscr} \ra \rscr' $ defined by
\begin{equation}\label{pmapcomplete}
\bfu_{\rscr}(w,t,a)=\left(a\, , \; -3t + \, \frac{w-a -3 \epsilon(w+a)}{2}\, , \; 2t + \epsilon(w+a)\right)
\end{equation}
is a one-to-one correspondence with inverse $\bfu \mapsto (\wtr(\bfu),\ttr(\bfu), a(\bfu))$.
\ep
\begin{proof}
For any integer $n$ we have  $n=2 \floor{\frac{n}{2}}+\epsilon(n)$. 
If $w=a+2b+3c$,
then $\epsilon(w+a)=\epsilon(c)$, hence $c= 2 \floor{\frac{c}{2}}+\epsilon(w+a)$. 
Thus,
the map that sends $\bfu=(a,b,c)$ to $\left(\wtr(\bfu)=a+2b+3c,\ttr(\bfu)=\floor{\frac{c}{2}}, a(\bfu)\right)$ is
invertible with inverse $\bfu_{\rscr}$. 
It remains to show that $\bfu_{\rscr}$ maps $\mathscr{P}_{\rscr}$
onto $\rscr'$. 
To see this, observe that $\rscr'$ is defined by the conditions $a \geq 0$, $b \geq 0$,
  $c \geq 0$ and $2a+2b+3c\leq 6m-9$. 
  These conditions translate to $a \geq 0$, $ w-a-6t-3 \epsilon(w+a) \geq 0$,
$t \geq 0$, and $a+w \leq 6m-9$; they imply that $w \geq 0$. 
It follows that
  $(\bfu_{\rscr})^{-1}$ maps $\rscr'$ into $\mathscr{P}_{\rscr}$, and it is clear that
  $\bfu_{\rscr}$ maps  $\mathscr{P}_{\rscr}$  into $\rscr'$.
\end{proof}

Similarly, we can use the weight $\wtm$, the $t$-invariant $\ttm$ and the first coordinate $i$ to parametrize $\mscr'$.
Recall that, for an integer $n$, we set
$\eta(n)=n-3 \floor{\frac{n}{3}} \in \{0,1,2\}$.
\bp \label{mprimeparamwta}
Fix an integer $m \geq 2$. Recall that
\begin{equation*}
  \mscr'=
\left\{\bfv=
(i,j,k) \in \N^3:  2i+2j+3k \leq 9(m\!-\!1), i\!+\!j \leq 3(m\!-\!1),\; i\!+\!j\!+\!k \geq 3m\!-\!2
\right\},
\end{equation*}
and let
$$\mathscr{P}_{\mscr}=
\left\{(w,t,i) \in \N^3: \,
w+i \leq 6m\!-\!9, \, 3(m\!-\!1\!+\!t) - w +\eta(w\!+\!i) \leq i \leq 3(m\!-\!1\!-\!t)\!-\!\eta(w\!+\!i) \right\}.
$$
The function $\bfv_{\mscr}:\mathscr{P}_{\mscr} \ra \mscr' $ defined by
\begin{equation*}\label{pmmapcomplete}
\bfv_{\mscr}(w,t,i)=
\left(i \, ,
 \; -3t + \, 3(m\!-\!1) -i - \eta(w+i)\, ,
 \; 2t -(m-2)+ \frac{w+i+ 2\eta(w+i)}{3}
\right)
\end{equation*}
is a one-to-one correspondence with inverse $\bfv \mapsto (\wtm(\bfv),\ttm(\bfv), i(\bfv))$.
\ep
\begin{proof}
Given $\bfv=(i,j,k)$, let $w=\wtm(\bfv)=i+2j+3k-3m$. 
Then, $w+i \equiv 2i+2j \equiv -(i+j)$ modulo $3$, hence
$$
t=\ttm(\bfv)  = \floor{\ds \frac{3(m\!-\!1)-(i+j)}{3}} =(m-1)+\frac{-i-j-\eta(w+i)}{3}.
$$
Solving this equation for $j$ we obtain
$$
j= 3(m-1-\ttm(\bfv))-i-\eta(w+i).
$$
We can  solve the equation $w=i+2j+3k-3m$ for $k$ to obtain the formula for $k$ as a function of $(w,t,i)$. 
This shows that
$\bfv_{\mscr}$ is bijective from $\Z^3$ to $\Z^3$, and it remains to check that it maps $\mathscr{P}_{\mscr}$ onto
$\mscr'$.
To see this, observe that $\mscr'$ is defined by the conditions $i \geq 0$, $j \geq 0$, $i+j \leq 3(m-1)$,
$2i+2j+3k\leq 9(m-1)$ and $i\!+\!j\!+\!k \geq 3m\!-\!2$, which imply $k \geq 1$.
These conditions translate in $(w,t,i)$ parameters to $i \geq 0$,  $i \leq 3(m-1-t)-\eta(w+i)$, $t \geq 0$,
$w+i \leq 6m-9$ and $3(m-1+t) - w +\eta(w+i) \leq i$ respectively, and they imply $w \geq 0$.
Thus, $\mathscr{P}_{\mscr}$ corresponds to $\mscr'$ under $\bfv_{\mscr}$.
\end{proof}

\noindent
The geometry of $\rscr'_w$ and $\mscr'_w$ depends on the  {\em threshold number}
\begin{equation}\label{threshold-def}\index{tz @ $\tau_w$ the threshold number}
  \tau_w=\floor{w/3}-(m\!-\!1).
\end{equation}
In fact, as we will show presently, the two parameter spaces $\mathscr{P}_\rscr$ and $\mathscr{P}_\mscr$
coincide in the region $0 \leq t \leq \tau_w$; in this region, the section of $\mathscr{P}_\rscr = \mathscr{P}_\mscr$  with the plane $w=w_0$ is just the rectangle $0 \leq a \leq 6m-9-w_0$, $0 \leq t \leq \tau_{w_0}$. We are thus led to the following definition.

\bd \label{ronedef}\index{rspr @ $\rect$ the rectangular region of $\rscr'$}\index{rtri @ $\rtri$ the triangular region of $\rscr'$}
\index{mspr @ $\mrect$ the rectangular region of $\mscr'$}\index{mtri @ $\mtri$ the triangular region of $\mscr'$}
We let $\rect$  denote the set of monomials $\bfu \in \rscr'$ that satisfy $\ttr(\bfu) \leq \tau_{\wtr(\bfu)}$, and
$\mrect$  denote the set of monomials $\bfv \in \mscr'$ that satisfy $\ttm(\bfv) \leq \tau_{\wtm(\bfv)}$.
 We call them the {\em rectangular} regions.
We let $\rtri$ (resp. $\mtri$) denote the complement of $\rect$ in $\rscr'$ (resp. of $\mrect$ in $\mscr'$),
and call them the {\em triangular} regions.
We define $\rectw$ and $\mrectw$ as the set of monomials of weight $w$
in $\rect$ and $\mrect$ respectively, and likewise for
 $\rtriw$ and $\mtriw$.
\ed

\br
Recall that $\rscr'_w$ and $\mw'$ are nonempty if and only if $0 \leq w \leq 6m-9$.
For $0 \leq w \leq 3m-4$, the threshold number $\tau_w$ is negative, so $\rw'=\rtriw$ and $\mw'=\mtriw$ consist only of the triangular regions. In the range $3(m\!-\!1) \leq w \leq 6m\!-\!10$ we have both the triangular and the rectangular regions, while $\rw'=\rectw$  and $\mw'=\mrectw$ for $w=6m-9$.
\er

\subsection{The bijection $\varphi^\hrectangle$ between the rectangular regions}\label{rectangular}
In Theorem \ref{phirectregion} we will show that composing the parametrization $\bfv_\mscr$ of Proposition \ref{mprimeparamwta}
with the inverse of  the parametrization $\bfu_\rscr$ of Proposition \ref{rprimeparamwta} yields a bijection
$\varphi^\hrectangle: \rect \ra \mrect$ that satisfies divisibility. For a given $\bfu =(a,b,c)\in \rect$, it turns out that the {\em multiplier}
$\bfq^\hrectangle (\bfu)=\varphi^\hrectangle(\bfu)-\bfu$ has the form
\begin{equation}\label{rectmultiplier}\index{ql @ $\bfq_\lambda=(0,3\lambda, m -2\lambda)$ multiplier having $a'=0$}
\bfq_\lambda=
(0,0,m)+ \lambda(0,3,-2)=(0,3\lambda, m -2\lambda),
\end{equation}
where $\lambda$ is an integer that depends only on $a+w$.

\br \label{rmkqlambda}
Suppose that $\bfq =(a',b',c')$ is a monomial of weight $a'+2b'+3c'=3m$. 
If $a'=0$,
then $b' \equiv 0 \pmod{3} $, thus,  there is a unique integer $\lambda$ such that
$b'=3 \lambda$, and  $\bfq=\bfq_\lambda$.
\er

We first explain how to compute the integer
$\lambda$ that determines the multiplier $\bfq^\hrectangle (\bfu)$.

\bd \label{lambdarhodefinition} \index{aalambda @ $\lambda(n)$}\index{aarho @ $\rho(n)$}
For an integer $n$, define $\rho(n)$ as the unique integer in
$\{-6,-4,-3,-2,-1,1\}$ that is congruent  to $n$ modulo $6$, and we define
$$
\lambda(n)=\frac{n-\rho(n)}{6}.
$$
We introduce the notation
\begin{equation}\label{deltadefinition}\index{aadeltaw @ $\delta_w$ the largest value of $a(\bfu)$ for $\bfu \in \rscr'_w$}
\delta_w=6m-9-w.
\end{equation}
Note that $\delta_w$ is the largest value of $a(\bfu)$ for $\bfu$ in $\rscr'_w$, since  $a+w \leq 6m-9$ holds in $\rscr'_w$.

Given $\bfu=(a,b,c) \in \rect$ of weight $w$, define  $\lambda(\bfu)=\lambda(\delta_w-a(\bfu))$ and $\rho(\bfu)=\rho(\delta_w-a(\bfu))$.\index{aalambdau @ $\lambda(\bfu)=\lambda(\delta_w-a(\bfu))$}
\index{aarhou @ $\rho(\bfu)=\rho(\delta_w-a(\bfu))$}
Then,
\begin{equation}\label{lambda-u-definition} 
\delta_w-a(\bfu)=6 \lambda(\bfu)+\rho(\bfu)
\end{equation}
Note that $\lambda(\bfu)$ and $\rho(\bfu)$ only depend on the sum $w+a$ of the weight $w$ and the first entry $a$ of $\bfu$.
\ed

\br \label{remark2lambda}
One could also define  $ \lambda(n)$ as the unique function of $n \in \Z$ that satisfies the properties:
\begin{itemize}
  \item $\lambda(n)=0$  if $n=1$,
  \item $\lambda(n)=1$ if $n=0,2,3,4,5$,
  \item $\lambda(n+6)=\lambda(n)+1$.
\end{itemize}
Then, one would recover $\rho(n)$ as $\rho(n)=n -6\lambda(n)$.
Explicit formulas for $\rho(n)$ and $\lambda(n)$ are
\begin{equation} \label{lambdarhoformulas}
\rho(n)=3 \epsilon(n)+2 \eta(-n)-6, \quad  \lambda(n) = \ceil{\frac{n+1}{2}}- \ceil{\frac{n}{2}}.
\end{equation}
\er

\bt \label{phirectregion}\index{aavphir @ $\varphi^\hrectangle$ the bijection between the rectangular regions}
The composition of the parametrization $\bfv_\mscr$ of Proposition \ref{mprimeparamwta}
with the inverse of  the parametrization $\bfu_\rscr$ of Proposition \ref{rprimeparamwta} is  a bijection
$\varphi^\hrectangle: \rect \ra \mrect$ given by the formula
\begin{equation}\label{phirect}
\varphi^\hrectangle(\bfu)=\bfu+ \bfq_{\lambda(\bfu)}.
\end{equation}
This bijection satisfies divisibility and preserves both the weight and the $t$-invariant.
\et

\begin{proof}
Let
$
\prect=\{(w,t,a) \in \Z^3: 3(m\!-\!1) \leq w \leq 6m-9, \, 0 \leq t \leq \tau_w, \,  0 \leq a \leq \delta_w\}.
$
We claim that $\prect$ is the intersection of the region $t \leq \tau_w$ with the parameter space $\pscr_{\rscr}$
of Proposition \ref{rprimeparamwta}, as well as the intersection of the region $t \leq \tau_w$ with the parameter space $\pscr_{\mscr}$ of Proposition \ref{mprimeparamwta}. 

We first treat the case of $\pscr_{\rscr}=
\left\{(w,t,a) \in \N^3: \;
w+a \leq 6m-9, \; w-a-6t-3 \epsilon(w+a) \geq 0 \right\}.
$
If $(w,t,a) \in \pscr_{\rscr} \cap \{t \leq \tau_w\}$, then we have $\tau_w \geq t \geq 0$, hence $w \geq 3(m-1)$.
Next we show that, if $t \leq \tau_w$, then the condition $w+a \leq 6m-9$ implies $ w-a-6t-3 \epsilon(w+a) \geq 0$.
For this, note that $t \leq \tau_w$ implies $6t \leq 2w-6(m\!-\!1)$, hence, if $w+a \leq 6m-9$, then
\begin{equation}\label{bpositive}
w-a-6t-3 \epsilon(w+a) \geq -a-3 \epsilon(w+a)+6(m\!-\!1)-w \geq 3-3 \epsilon(w+a) \geq 0.
\end{equation}
This shows $\pscr_{\rscr} \cap \{t \leq \tau_w\}=\prect$.

We now treat the case of $\pscr_{\mscr}$, recall that
$$
\mathscr{P}_{\mscr}=
\left\{(w,t,i) \in \N^3: \,
w+i \leq 6m-9, \, 3(m-1+t) - w +\eta(w+i) \leq i \leq 3(m-1-t)-\eta(w+i) \right\}
$$
The equality $\pscr_{\mscr} \cap \{t \leq \tau_w\}=\prect$ will follow once we show
that the inequalities $w+i \leq 6m-9$ and $t \leq \tau_w$ imply
$$
3(m-1+t) - w +\eta(w+i) \leq i \leq 3(m-1-t)-\eta(w+i).
$$
We start with the inequality on the left: from $t \leq \tau_w=\frac{w-\eta(w)}{3}-m+1$ it follows
$$
3(m-1+t) - w +\eta(w+i) \leq \eta(w+i)-\eta(w) \leq i.
$$
To show the inequality on the right, we use both constraints $t \leq \tau_w$ and $w+i \leq 6m-9$:
$$
3(m-1-t)-\eta(w+i) \geq
6(m-1)-w +\eta(w)-\eta(w+i)> 6m-9-w \geq i.
$$

To construct the bijection $\varphi^\hrectangle$, denote by $\bfu^{\hrectangle}$\index{uparam@$\bfu^{\hrectangle}$} and
$\bfv^{\hrectangle}$ the restrictions to $\prect$ of the bijections $\bfu_{\rscr}$ and $\bfv_{\mscr}$
of Propositions \ref{rprimeparamwta}  and \ref{mprimeparamwta}. 
Then, define
 $\varphi^\hrectangle= \bfv^\hrectangle \circ (\bfu^{\hrectangle})^{-1}$. 
 By construction,
 $\varphi^\hrectangle: \rect \ra \mrect$ is a bijection that preserves the weight and the $t$-invariant.

To finish, we need to show that $\varphi^\hrectangle$
satisfies divisibility, i.e.,  that for each $(w,t,a) \in \prect$ the {\em multiplier vector}
$
\bfq^{\hrectangle}(w,t,a)=\bfv^{\hrectangle}(w,t,a)-\bfu^{\hrectangle}(w,t,a)
$
has non-negative entries. 
To compute $\bfq^{\hrectangle}(w,t,a)$,  let $\bfu=\bfu^{\hrectangle}(w,t,a)=(a,b,c)$ and write
$\delta_w-a=6\lambda (\bfu) +\rho(\bfu) $. 
Then, $a=\delta_w-6\lambda - \rho$ and
\begin{equation}\label{crholambda}
c=2t+\epsilon(w+a)=2t+1-\epsilon(\rho).
\end{equation}
Although  we do not need an explicit formula for $b$ in terms of $(w,t,\lambda,\rho)$ at this point, we compute it now for later use. As $b=\frac{w-a-3c}{2}$, we obtain
\begin{equation}\label{brholambda}
b=w-3(m-1)+3(\lambda-t)+\frac{\rho+3 \epsilon(\rho)}{2} =
3(\lambda +\tau_w-t)+\eta(w)+\frac{\rho+3 \epsilon(\rho)}{2} .
\end{equation}
Similarly, if we let $(i,j,k)=\bfv^{\hrectangle}(w,t,a)$, then $i=a=\delta_w-6\lambda - \rho$ and, substituting in the formula of Proposition \ref{mprimeparamwta}
for $\bfv^{\hrectangle}$ we obtain
\begin{equation}\label{krholambda}
k=2t+1+m-2\lambda + \frac{2\eta(-\rho)-\rho-6}{3}.
\end{equation}
We check  that $\frac{2\eta(-\rho)-\rho-6}{3}=-\epsilon(\rho)$ for all $6$ possible values of $\rho$. It follows that
$$
\bfq^{\hrectangle}(w,t,a)=(0, \frac{3m-3m+6\lambda}{2},m -2 \lambda)=
\bfq_\lambda
$$
as claimed.
Note that the multiplier $\bfq^{\hrectangle}(w,t,a)$ does not depend on $t$.

The vector $\bfq^{\hrectangle}(w,t,a)$
has non-negative entries if and only if  $ 0 \leq \lambda(\delta_w\!-\!a) \leq m/2$. 
We observe that for
$(w,t,a) \in \prect $ we have $ 0 \leq \delta_w-a \leq \delta_w$. 
The equality $\delta_w-a=6 \lambda +\rho$ together with $\delta_w-a  \geq 0$ implies that $\lambda \geq 0$.
 On the other hand,
in $\prect$ the inequality  $3(m\!-\!1) \leq w$ holds; together with $\delta_w-a \leq \delta_w$ it implies
that $
6 \lambda +\rho \leq 6m-9-w \leq 3m-6.
$
As $\rho \leq 1$, we conclude that  $2 \lambda \leq m$, and the proof is complete.
\end{proof}

\noindent
\br
The maximum value of $\lambda(\delta_w-a)$ for $0 \leq a \leq \delta_w$ is $m\!-\!1 \!-\! \ds
\floor{\frac{w}{6}}$ when $w$ is congruent modulo $6$  to a number
between $0$ and $3$, and is $m-2 - \ds \floor{\frac{w}{6}}$
when $w$ is congruent modulo $6$ to either $4$ or $5$ (the maximum is
attained at either $a=0$ or $a=1$).
\er

\subsection{An order relation on the rectangular regions} \label{orderrectangular}
Our analysis of the uniqueness properties of  $\varphi^\hrectangle$ requires an order relation
on $\prect_w$ that is suggested by the decomposition $\delta_w-a=6 \lambda+\rho$. To introduce this order relation,
fix a $w$ in the interval $3(m-1) \leq w \leq 6m-9$, so that the rectangular region $\rectw$
is nonempty and parametrized by the set $\prect_w$ of integral points of a rectangle:
$$
\prect_w=\{(t,a) \in \Z^2: 0 \leq t \leq \tau_w, \, 0 \leq  a \leq \delta_w= 6m-9-w\}.
$$
As in the proof of Theorem \ref{phirectregion}, we can parametrize the interval $0 \leq  a \leq \delta_w$ using
the decomposition $\delta_w-a=6\lambda+\rho$: here $\lambda$ and $\rho$ should be considered as independent parameters,
subject to the restrictions that $\rho \in \{-6,-4,-3,-2,-1,1\}$ and $0 \leq 6 \lambda +\rho \leq \delta_w$.
Then, we order pairs $(t,a(\lambda,\rho))$ lexicographically with respect to $(\lambda, \rho,t)$:

\bd \label{order-rect}\index{or @ $\prec$ the order relation on the rectangular regions}
We define a total order relation $\prec$ on $\prect_w$ as follows: $(t,a) \prec (t',a')$ if either $ \lambda(\delta_w\!-\!a)<\lambda(\delta_w\!-\!a')$,
or $ \lambda(\delta_w\!-\!a)=\lambda(\delta_w\!-\!a')$ and $a > a'$, or $ \lambda(\delta_w\!-\!a)=\lambda(\delta_w\!-\!a')$, $a = a'$
and $t<t'$.

Given $\bfu$, $\bfu'$ in $\rectw$, we write
$\bfu \prec \bfu'$ if $(\ttr(\bfu),a(\bfu)) \prec (\ttr(\bfu'),a(\bfu'))$.
As $\lambda(\bfu)=\lambda(\delta_w-a(\bfu))$ and $\rho(\bfu)=\rho(\delta_w-a(\bfu))$, it follows that
$\bfu \prec \bfu'$ if and only if $(\lambda(\bfu),\rho(\bfu),\ttr(\bfu))$ precedes $(\lambda(\bfu'),\rho(\bfu'),\ttr(\bfu'))$
lexicographically.
Similarly, given $\bfv$, $\bfv'$ in $\mrectw$, we write
$\bfv \prec \bfv'$ if $(\ttm(\bfv),i(\bfv)) \prec (\ttm(\bfv'),i(\bfv'))$, so that the bijection $\vphi^{\hrectangle}:\rectw \ra \mrectw$ is order preserving.
\ed

In the tables in Appendix \ref{AppendixTable}, 
monomials in $\rectw$ are listed according to the order relation $\prec$.

\br \label{orderdeltawn}
Define a new order $\succ_6$\index{os @ $\succ_6$ an order relation on the integers} on the set of integers
modifying the usual order by declaring
that $\delta_w-6n-1$ is larger than $\delta_w-6n$ for every $n$. Thus, if we set
$\delta_{w,n}=\delta_w-6n$\index{aadeltawn @ $\delta_{w,n}=\delta_w-6n$}, the new order is
\begin{equation}\label{neworder}
\cdots \delta_{w,n}\!-\!1 \succ_6 \delta_{w,n} \succ_6 \delta_{w,n} \!-\!2 \succ_6\delta_{w,n} \!-\!3 \succ_6 \delta_{w,n} \!-\!4 \succ_6 \delta_{w,n} \!-\!5 \succ_6 \delta_{w,n} \!-\!7 = \delta_{w,n+1}\!-\!1\succ_6 \delta_{w,n+1}
\cdots
\end{equation}
Then, $(t,a) \prec (t',a')$ if and only if either $a \succ_6 a'$ or $a=a'$ and $t <t'$.
\er

\bex
As $\rho \in \{-6,-4,-3,-2,-1,1\}$ and $0 \leq 6 \lambda +\rho \leq \delta_w$, 
we have $\lambda(\bfu)=0$ if and only if $a(\bfu)=\delta_w-1$
(and in this case $\rho(\bfu)=1$):  monomials with $a=\delta_w-1$ are smaller than those with $a=\delta_w$, which have $\lambda=1$ and
$\rho=-6$, then come those with $a=\delta_w-2$, which have $\lambda=1$ and $\rho=-4$, and so on as in (\ref{neworder}).
\eex

For later use,  for each integer $r \geq 0$
we identify in Lemma \ref{urrectangular} below the largest monomial $\bfu^{(r)}$ in $\rectw$ whose
second entry $b$ is at most $r +\eta(w)$, and then in Lemma \ref{aneworder} we describe the segments
$\bfu^{(r-1)} \prec \bfu \preceq \bfu^{(r)}$  in terms of the first entry $a$ of $\bfu$.

\bl \label{urrectangular}
Fix an integer $w$ in the interval  $3(m-1) \leq w \leq 6m-9$. Given an integer $0 \leq r \leq \ds \floor{\frac{\delta_w}{2}}$, let
$\bfu^{(r)}=\bfu^{\hrectangle}(w,\tau_w,\delta_w-2r)$.
Then, $\bfu^{(r)}=(\delta_w-2r, r+\eta(w),  2\tau_w +1)$, $\lambda (\bfu^{(r)})=\floor{\frac{r+3}{3}}$, and
$\rho (\bfu^{(r)})$ is even and congruent to $2r$ modulo $6$.
Furthermore,
$$
b (\bfu) \geq r+\eta(w)+1 \quad \mbox{if $\bfu \in \rectw$ and $\bfu \succ \bfu^{(r)}$}.
$$
\el

\begin{proof}
By definition of $\bfu^{\hrectangle}$, the first entry $a$ of $\bfu^{(r)}$ is $\delta_w-2r$. In particular,
$a+w=6m-9-2r$ is odd, thus formula (\ref{crholambda})
gives $c(\bfu^{(r)})= 2\tau_w +1$.

To compute $b(\bfu^{(r)})$, we use  formula (\ref{brholambda}). 
For this, we recall
$$
 \lambda(\bfu^{(r)})=  \lambda(\delta_w\!-\!a(\bfu^{(r)})=\lambda(2r), \quad  \rho(\bfu^{(r)})= \rho(2r)
$$
As $2r= 6\lambda(2r)+\rho(2r)$, we see $\rho(2r)$ is even, that is, one of the integers $\{-6,-4,-2\}$,
and $r=3 \lambda(2r)+ \rho(2r)/2$. Formula (\ref{brholambda}) now gives $b(\bfu^{(r)})= r+\eta(w)$.
The fact that  $\lambda(2r)=\floor{\frac{r+3}{3}}$ also follows from  $r=3 \lambda(2r)+ \rho(2r)/2$.

It remains to show that
$b (\bfu) > b(\bfu^{(r)})$ if $\bfu \in \rectw$ and $\bfu \succ \bfu^{(r)}$. So suppose $\bfu \succ \bfu^{(r)}$. As
$\ttr(\bfu^{(r)})=\tau_w$, the largest value of $\ttr$ on $\rectw$, we must have
either $\lambda(\bfu) >\lambda (\bfu^{(r)})$ or $\lambda(\bfu)=\lambda (\bfu^{(r)})$ and $\rho(\bfu)>\rho (\bfu^{(r)})$.
If $\lambda(\bfu) >\lambda (\bfu^{(r)})$, then by   (\ref{brholambda}) we obtain
$$
b(\bfu) -b(\bfu^{(r)}) \geq 3+\frac{\rho(\bfu)+3 \epsilon(\rho(\bfu))}{2}   -\frac{\rho (\bfu^{(r)})}{2}  \geq 1
$$
because $\rho(\bfu)+3 \epsilon(\rho(\bfu) )\geq  -6$ and $\rho (\bfu^{(r)}) \leq -2$.
Finally, if $\lambda(\bfu)=\lambda (\bfu^{(r)})$ and $\rho(\bfu)>\rho (\bfu^{(r)})$, then by   (\ref{brholambda}) we obtain
$$
b(\bfu) -b(\bfu^{(r)}) \geq \frac{\rho(\bfu)-\rho (\bfu^{(r)})+3 \epsilon(\rho(\bfu))}{2} >0.
$$
\end{proof}

\bl \label{aneworder}
Suppose $\bfu=(a,b,c) \in \rectw$, $r\in \mathbb{N}$ satisfies $r \equiv 1 \pmod{3}$, 
and $r_{max}=\floor{\frac{\delta_w}{2}}$. Then,

\begin{enumerate}
  \item   $\bfu \preceq \bfu^{(0)}$ if and only if $a \in \{\delta_w,\delta_w-1\}$;
  \item  $\bfu^{(r-1)} \prec \bfu \preceq \bfu^{(r)}$ if and only if $a=\delta_w-2r$;
  \item  $\bfu^{(r)} \prec \bfu \preceq \bfu^{(r+1)}$ if and only if $a \in \{\delta_w-2r-1,\delta_w-2r-2 \}$;
  \item  $\bfu^{(r+1)} \prec \bfu \preceq \bfu^{(r+2)}$ if and only if $a \in \{\delta_w-2r-3,\delta_w-2r-4,\delta_w-2r-5 \}$;
  \item $\bfu^{(r_{max})} \prec \bfu$ if and only if $w$ is congruent to either $0$ or $4$ modulo $6$, and $a=0$.
\end{enumerate}
\el

\begin{proof}
It follows from (\ref{neworder}), since $\bfu^{(r)}$ is the largest monomial in $\rectw$ having  $a= \delta_w -2r$ because $\ttr (\bfu^{(r)})=\tau_w$, and $\tau_w$ is the largest value of $\ttr$ for monomials in the rectangular region.
\end{proof}

\subsection{Uniqueness in the rectangular regions}
In this subsection we prove (Corollary \ref{unirectangle2}) that
the bijection $\varphi^\hrectangle:\rectw \ra \mrectw$ of Theorem \ref{phirectregion} satisfies uniqueness: any bijection $\rectw \ra \mrectw$ with the same multipliers multiset as $\varphi^\hrectangle$ must in fact coincide with $\varphi^\hrectangle$.
In order to avoid  confusion, we fix the following terminology: given a bijection $\vphi:\rectw \ra \mrectw$, we say that the vector $\bfq(\bfu)= \vphi (\bfu)-\bfu$ is the $\vphi$-multiplier of $\bfu$; the multiplicity of a $\vphi$-multiplier
$\bfq$ is the number of monomials $\bfu \in \rectw$ such that $\bfq=\bfq(\bfu)$. The multiset
$\mathscr{Q}_\varphi$ of multipliers of $\varphi$
is the set of $\varphi$-multipliers $\bfq$ counted with their multiplicities.

\bl \label{lemmaunirectangle}
Fix a monomial $\bfu \in \rectw$ and let $n_1 < n_2$ be integers. If $\bfu+\bfq_{n_1}$ and $\bfu+\bfq_{n_2}$
belong to $\mrectw$, then  $\bfu +\bfq_{n_2} \prec \bfu +\bfq_{n_1}$.
\el
\begin{proof}
Note that the first entry of both $\bfu+\bfq_{n_1}$ and $\bfu+\bfq_{n_2}$ is equal to the first entry $a$ of $\bfu$,
and the $\scrM$-weight of both  $\bfu+\bfq_{n_1}$ and $\bfu+\bfq_{n_2}$ is equal to $w=\wtr(\bfu)$.
Hence,  by Proposition \ref{mprimeparamwta}, there exist integers $t_1$ and $t_2$ such that
$\bfu+\bfq_{n_1}=\bfv_{\scrM}(w,t_1,a) $  and $\bfu+\bfq_{n_2}=\bfv_{\scrM}(w,t_2,a) $.

Using the formula in Proposition \ref{mprimeparamwta}  for $\bfv_{\scrM}$ we compute
 \begin{equation}\label{lemnew}
 \bfq_{n_2}-\bfq_{n_1}=\bfv_{\scrM}(w,t_2,a) - \bfv_{\scrM}(w,t_1,a)= (t_2-t_1) (0,-3,2)
  \end{equation}
On the other hand, as $\bfq_n=(0,3n, m-2n)$, we have
 $ \bfq_{n_2}-\bfq_{n_1}=(n_2-n_1)(0,3,-2)$.
We conclude  that $t_2-t_1=-(n_2-n_1)<0$, hence $(a,t_2) \prec (a,t_1)$, which is equivalent to
$\bfv_{\scrM}(w,t_2,a) \prec \bfv_{\scrM}(w,t_1,a)$.
\end{proof}

\bc \label{unirectangle}
Let $3(m-1) \leq w \leq 6m\!-\!9, 
\bfu_0 \in \rect_w$.
If $\hphi:\rscr'_w \ra \mscr'_w$ is a bijection such that
\begin{enumerate}
\item  $ \hphi(\bfu_0) \in \mrectw$,
\item $ \hphi(\bfu_0) =\bfu_0+\bfq_n$ for some $n \geq \lambda(\bfu_0)$,
\item  $\varphi^\hrectangle (\bfu_0) \preceq \hphi(\bfu_0) $ in $\mrectw$,
\end{enumerate}
then $n=\lambda(\bfu_0)$ and $\hphi(\bfu_0)=\varphi^\hrectangle (\bfu_0)$.

\ec
\begin{proof}
Let $\lambda_0=\lambda(\bfu_0)$, so that $\varphi^\hrectangle (\bfu_0)=\bfu_0+\bfq_{\lambda_0}$ by Theorem \ref{phirectregion}.
By assumption $ \hphi(\bfu_0) =\bfu_0+\bfq_n$ for some $n \geq \lambda_0$.  If we had
$n > \lambda_0$, then by Lemma \ref{lemmaunirectangle} we would have
$ \hphi(\bfu_0) = \bfu_0+\bfq_n \prec \bfu_0+\bfq_{\lambda_0} = \varphi^\hrectangle (\bfu_0),$
a contradiction. 
Thus, $ n= \lambda_0$ and  $\hphi(\bfu_0)=\varphi^\hrectangle (\bfu_0)$.
\end{proof}

\bc[Uniqueness of $\varphi^{\hrectangle}$] \label{unirectangle2}
Let $3(m-1) \leq w \leq 6m\!-\!9$ and  $\varphi^\hrectangle:\rectw \ra \mrectw$ be as in  Theorem \ref{phirectregion}. 
If
$\hphi:\rectw \ra \mrectw$ is a bijection with the same multiset of multipliers as $\varphi^\hrectangle$,
then $\hphi=\varphi^\hrectangle$.
\ec
\begin{proof}
We prove the result by induction on the order $\prec$.  
Fix a monomial  $\bfu_0$ in $\rectw$ and assume
$\hphi(\bfu)=\varphi^\hrectangle(\bfu)$ for all monomials $\bfu$ smaller than $\bfu_0$. We need to show $\hphi(\bfu_0)=\varphi^\hrectangle(\bfu_0)$.

Let $\lambda_0=\lambda(\bfu_0)$, so that $\varphi^\hrectangle(\bfu_0)=\bfu_0+\bfq_{\lambda_0}$.
 The $\hphi$-multiplier $\hphi(\bfu_0)-\bfu_0$  of $\bfu_0$ is, by assumption, also a $\varphi^\hrectangle$-multiplier,  thus it has the form
$\bfq_n$ for some integer $n$. 
The two bijections $\hphi$ and $\varphi^\hrectangle$ agree on monomials strictly smaller than $\bfu_0$, hence they have the same multipliers on these monomials. 
A monomial  $\bfq_\lambda$ with $\lambda < \lambda_0$ arises as a
$\varphi^\hrectangle$-multiplier $\varphi^\hrectangle(\bfu)-\bfu$ precisely for those monomials $\bfu$ that satisfy $\lambda(\bfu)=\lambda$,
and these monomials are all smaller than $\bfu_0$. 
Since the multipliers $\bfq_\lambda$ have the same multiplicity as multipliers of
$\hphi$ and of $\varphi^\hrectangle$, 
and since the two bijections coincide on monomials smaller than $\bfu_0$, we conclude that the $\hphi$-multiplier $\bfq_n$ of $\bfu_0$ has $n \geq \lambda_0$.

Finally, let   $\bfv=\hphi(\bfu_0)$.
Suppose $ \bfv \prec \varphi^\hrectangle (\bfu_0)$. Then, as $\varphi^\hrectangle$ is order preserving, the inverse image $\bfu=(\varphi^\hrectangle)^{-1}(\bfv)$ in $\rectw$ is smaller than $\bfu_0$. 
By induction, $\hphi$ and $\varphi^\hrectangle$ coincide on $\bfu$, hence
$  \hphi(\bfu)=\varphi^\hrectangle(\bfu)=\hphi(\bfu_0).$
But this is impossible because $\hphi$ is a bijection. We conclude that $ \bfv \succeq \varphi^\hrectangle (\bfu_0)$, and then Corollary \ref{unirectangle}
implies  $\hphi(\bfu_0)=\varphi^\hrectangle(\bfu_0)$.
\end{proof}

\section{Construction of the bijection $\varphi$ satisfying divisibility - part two} \label{triangular}
In this section, we construct a bijection $\varphi^{\triangle}:\rtri \ra \mtri$ with the divisibility property.
One could adapt the argument by which we obtained $\varphi^{\hrectangle}$:
in the region $t \geq \tau_w+1$, the two parameter spaces $\pscr_{\rscr}$ and $\pscr_{\mscr}$ do not coincide,
but, fixed a pair $(w,t)$ satisfying  $t \geq \tau_w+1$, the set $A_{w,t}$ of $a$'s such that $(w,t,a) \in \pscr_{\rscr}$
and the set $I_{w,t}$ of $i$'s such that $(w,t,i) \in \pscr_{\mscr}$ are in one-to-one correspondence, they essentially
differ by a translation. 
Thus, there is a unique order preserving bijection $o_{w,t}:A_{w,t} \ra I_{w,t}$,
and one obtains a bijection $\rtri \ra \mtri$ sending $\bfu_{\rscr}(w,t,a)$ to $\bfv_{\mscr}(w,t,o_{w,t}(a))$.
This bijection  satisfies divisibility, except for a few cases where it can be modified in such a way that divisibility holds without exceptions.
However, it is unclear whether this bijection has the uniqueness properties we will need in the sequel. 
We will therefore construct a different bijection $\varphi^{\triangle}: \rtri \ra \mtri$, which is more canonical, in the sense that preserves more invariants,  and has better properties for our purposes.

We construct the map  $\varphi^{\triangle}$ on a domain slightly larger than $\rtri$, and then show that it restricts to a bijection
$\rtri \rightarrow \mtri$. More precisely, we define
\begin{equation} \label{rtdef}
\index{rtt @ $\rtt$ the extended triangular region of $\rscr$}
\index{mtt @ $\mtt$ the extended triangular region of $\mscr$}
  \rtt=
\left\{\bfu=
(a,b,c) \in \rscr: \ttr (\bfu) \geq \tau_{\wtr(\bfu)}+1
\right\}
\end{equation}
and
\begin{equation} \label{mtdef}
  \mtt=
\left\{\bfv=
(i,j,k) \in \mscr:  \ttm (\bfv) \geq \tau_{\wtm(\bfv)}+1
\right\}.
\end{equation}
\br \label{remarkrtt}
Recall that $\rscr$ is the set of $(a,b,c)\in\N^3$ that satisfy $a+b \leq 3(m-1)$.
The condition $\ttr (\bfu) \geq \tau_{\wtr(\bfu)}+1$
implies that $a+b \leq 3(m-1)$. Indeed, $\ttr (\bfu) \geq \tau_{\wtr(\bfu)}+1$ is equivalent to
\begin{equation}\label{ttrconditionexplicit}
  2a+4b+3c \leq 6m-12-3\epsilon(c)+2 \eta(a+2b)
\end{equation}
and implies that $a+b \leq 3m-4$.
Thus,
$$
  \rtt=
\left\{\bfu=
(a,b,c) \in \N^3: \ttr (\bfu) \geq \tau_{\wtr(\bfu)}+1
\right\}.
$$
\er

By Definition \ref{ronedef}, the triangular regions $\rtri$ and $\mtri$ are subsets of $\rtt$ and $\mtt$, respectively.

\bt \label{phitriregion} \index{aavphit @ $\varphi^\triangle$ the bijection between the triangular regions}
Fix an integer $m \geq 2$.
Let
$\dscr=
\left\{\bfd=
(w,r,s) \in \N^3: 2r+3s \leq w
\right\}$. Given $\bfd=(w,r,s) \in \dscr$, let $\wtd(\bfd)=w$ and $\ttd(\bfd)= \ds\floor{\frac{w-2r-3s}{6}}.$
Define
\begin{equation}\label{Dtdefinition}
  \dtt=
\left\{\bfd=
(w,r,s) \in \dscr: \ttd (\bfd) \geq \tau_{w}+1
\right\}
\end{equation}
and maps
 $\bfu^{\triangle}, \bfv^{\triangle}: \dtt \ra \N^3$ by the formulas
\begin{equation}\label{pmap}
\bfu^{\triangle}(w,r,s)=
\left(
3s+\eta(w+r), r, \frac{w-2r-3s-\eta(w+r)}{3}
\right)
\end{equation}
and
\begin{equation}\label{p'map-explicit}
\bfv^{\triangle}(w,r,s)=
\left(
\frac{6(m\!-\!1)\!-\!w \!-\!2r  \!+\!3s\!-\!\epsilon(w\!+\!s)}{2},
2r\!+\!\epsilon(w\!+\!s),
\frac{w\!-\!2r\!-\!s\!+2-\!\epsilon(w\!+\!s)}{2}
\right)
\end{equation}
Then, $\bfu^{\triangle}: \dtt \ra \rtt$ and $\bfv^{\triangle}: \dtt \ra \mtt$
are bijective. Furthermore, the composition
$\varphi^\triangle = \bfv^{\triangle} \circ (\bfu^{\triangle})^{-1} : \rtt \ra \mtt$
 satisfies divisibility and preserves weights and the $t$-invariants.
\et

\begin{proof}
We first explain how the parametrization (\ref{pmap}) of $\N^3$ by $(w,r,s)$ arises.
Consider, as in the beginning of the paper, the polynomial ring  $P=\Bbbk[x,y,z]$ with the grading obtained
assigning weight $1$ to the variable $x$, weight $2$ to $y$ and weight $3$ to $z$.
View $P=\Bbbk[x,y,z]$ as a graded module over its subring $P_0=\Bbbk[y,x^3,z^2]$. Then,
\begin{equation}\label{freebasisP0}
  \{1,x,x^2,z,xz,x^2z\}=\{x^\eta z^\epsilon: 0 \leq \eta \leq 2, 0 \leq \epsilon \leq 1\}
\end{equation}
 is a free basis of $P$ over $P_0$, whose elements have weight $0,1,2,3,4,5$ respectively.
 Given a monomial $y^{r}x^{3s}z^{2t}$ in $P_0$, we denote by $r_\pscr$, $s_\pscr$ and $t_\pscr$ its degree with respect to
 $y$, $x^3$ and $z^2$ respectively. We extend this multigrading to all monomials of $\Bbbk[x,y,z]$
 setting $r=s=t=0$ on the above chosen basis. A monomial in $P$ is then determined by its weight $w$ and its multidegree $(r,s,t)$,
 and the last parameter $t$ is not needed as we can recover it from $w$, $r$ and $s$.  We make this explicit.
Identifying as usual a monomial $\bfu=x^ay^bz^c$ of $P$ with $(a,b,c)$,
 we have
   \begin{equation}\label{Pmultigrad}
     r_\pscr(\bfu)=b, \;
     s_\pscr(\bfu)= \floor{\frac{a}{3}}, \;
     t_\pscr(\bfu)= \floor{\frac{c}{2}}, \; \wtp(\bfu)=a+2b+3c.
   \end{equation}
Note that $t_\pscr=t_\rscr$ is the $t$-invariant we introduced for monomials in $\rscr$.  We can recover $a$, $b$, $c$ and
$t$ from $(w=\wtp(\bfu),r=r_\pscr(\bfu),s=s_\pscr(\bfu))$: for this, observe that $\eta(a)=\eta(w-2r)=\eta(w+r)$,
hence $a=3s+\eta(w+r)$, $b=r$, and $c=\ds \frac{w-3s-\eta(w+r)-2r}{3}$. 
In fact,  the function  $\bfu^{\triangle}:\dscr \ra \N^3$  defined by equation (\ref{pmap}) is the inverse of the function $\bfu \mapsto (\wtp(\bfu),    r_\pscr(\bfu),  s_\pscr(\bfu))$. It is clear that $\bfu^{\triangle}$
preserves weights, and it also preserves the $t$-invariants:
$$
t_{\pscr}(\bfu^\triangle (\bfd))= \floor{\frac{w-2r-3s-\eta(w+r)}{6}}= \floor{\frac{w-2r-3s}{6}}=\ttd(\bfd).
$$
It follows that the image of $\bfu^\triangle(\dtt)$ is $\rtt$, the set of monomials $\bfu=(a,b,c)  \in \N^3$ that satisfy
$\ttr ( \bfu) \geq \tau_{\wtp(\bfu)}+1$ - cf. Remark \ref{remarkrtt}.

We obtain a second parametrization of $\N^3$ by $(w,r,s)$ replacing in the above argument the weight two monomial
$y$ with $X^2$, the weight three monomial $x^3$ with $Z$, and the the weight six monomial
$z^2$ with $Y^3$. Namely, consider the polynomial ring  $\tilde{P}=\Bbbk[X,Y,X]$  with the grading obtained
assigning weight $1$ to the variable $X$, weight $2$ to $Y$ and weight $3$ to $Z$.
View $\tilde{P}$ as a weighted module over its subring $\tilde{P}_0=\Bbbk[X^2,Z,Y^3]$. Then
\begin{equation}\label{freebasisP0tilde}
\{1,X,Y,XY,Y^2,XY^2\}=\{X^{\epsilon_1} {Y^{\eta_1}}: 0 \leq \eta_1 \leq 2, 0 \leq \epsilon_1 \leq 1\}
\end{equation}
is a free basis of $\tilde{P}$ over $\tilde{P}_0$, whose elements have weight $0,1,2,3,4,5$, respectively.
 Given a monomial $X^{2r}Z^{s}Y^{3t}$ in $\tilde{P}_0$, we denote by $r_{\ptscr}$, $s_{\ptscr}$ and $t_{\ptscr}$ its degree with respect to
 $X^2$, $Z$ and $Y^3$, respectively. 
We extend this multigrading to all monomials of $\tilde{P}$
 setting $r=s=t=0$ on the above chosen basis: then for $\bfU=X^\atl Y^\btl Z^\ctl$ in $\tilde{P}$
    \begin{equation}\label{Ptildamultigrad}
     r_{\ptscr}(\bfU)=\floor{\frac{\atl}{2}}, \;
     s_{\ptscr}(\bfU)=\ctl, \;
     t_{\ptscr}(\bfU)= \floor{\frac{\btl}{3}},\; \w_{\ptscr}(\bfU)=\atl+2\btl+3\ctl.
   \end{equation}
A monomial $\bfU=(\atl,\btl, \ctl)$ in $\tilde{P}$ is then determined by its weight $w$ and its multidegree $(r,s)$: the function
$\bfU \mapsto \left(\w_{\ptscr}(\bfU), r_{\ptscr}(\bfU),  t_{\ptscr}(\bfU)\right)$ is a bijection $\N^3 \ra \dscr$ with inverse
\begin{equation}\label{utildeinverse}
\bfU(w,r,s)=\left(2r+\epsilon(w+s),  \frac{w-3s-\epsilon(w+s)}{2}-r,s\right)
\end{equation}
(note $\epsilon(\atl)=\epsilon(w+s)$).
This bijection preserves weights and $t$-invariants:
$$
t_{\ptscr}(\bfU(\bfd))= \floor{\frac{w-3s-2r-\epsilon(w+s)}{6}}= \floor{\frac{w-2r-3s}{6}}=\ttd(\bfd).
$$
It follows that the image of $\bfU(\dtt)$ is  the set $\widetilde{\pscr}^T$ of monomials $\bfU=(\atl,\btl,\ctl)  \in \N^3$ that satisfy
$t_{\ptscr} ( \bfU) \geq \tau_{\widetilde{\pscr}(\bfU)}+1$.

Now observe that, if $\bfU=(\atl,\btl,\ctl) $ belongs to $\widetilde{\pscr}^T$, then $\atl+\btl \leq 3(m-1)$. 
Indeed, the condition for
$\bfU$ to belong to $\widetilde{\pscr}^T$ is $\floor{\frac{\btl}{3}} \geq \floor{\frac{w}{3}}-m+2$, that is,
\begin{equation}\label{ttildeconditionexplicit}
\atl+\btl\leq 3m-6+\eta(\atl+2\btl)-\eta(\btl)-3\ctl,
\end{equation}
which implies $\atl+\btl\leq 3m-4$.
We remark that this is crucial for the rest of the argument, as it shows that  $\widetilde{\pscr}^T$
is contained in the set $\rscr$ of monomials $(\atl,\btl,\ctl) \in \N^3$ satisfying $\atl+\btl \leq 3(m-1)$.
We can then compose
$\bfU: \dtt \ra \widetilde{\pscr}^T $ with the bijective map $\Psi: \rscr \ra \mscr$ of Proposition \ref{Mstructure} .
Recall that  $ \mscr$ is the set of monomials $\bfv=
(i,j,k) \in \N^3$ satisfying  $ i\!+\!j \leq 3(m\!-\!1)$ and $ i\!+\!j\!+\!k \geq 3m\!-\!2$, and the formula for $\Psi$ is
\begin{equation}\label{Psi}
\Psi ((\atl, \btl, \ctl))=
\left(
3(m\!-\!1)\!-\!\atl\!-\!\btl, \atl, 1\!+\!\btl\!+\!\ctl
\right).
\end{equation}
The map $\Psi$ preserves weights because $\wtm(i,j,k)=i+2j+3k-3m$, and it also preserves $t$-invariants:
if $(i,j,k)=\Psi ((\atl, \btl, \ctl))$, then
$$
\ttm ((i,j,k))=
\floor{\frac{3(m\!-\!1)\!-\!i\!-\!j}{3}}= \floor{\frac{\btl}{3}}=t_{\ptscr} ((\atl, \btl, \ctl)).
$$
It follows that the map $\bfv^{\triangle} (\bfd)=\Psi  (\bfU(\bfd))$ of equation (\ref{p'map-explicit}) is a bijection
of $\dtt$ onto $\mtt$ that preserves weights and $t$-invariants, hence
$\varphi^\triangle = \bfv^{\triangle} \circ (\bfu^{\triangle})^{-1} : \rtt \ra \mtt$ too
is a bijection  preserving  weights and $t$-invariants.

It remains to show that $\varphi^\triangle$ satisfies divisibility, i.e., that for every
$(w,r,s) \in \dtt$ the multiplier
$\bfq^{\triangle}(w,r,s)=\bfv^{\triangle}(w,r,s)- \bfu^{\triangle}(w,r,s)$
has non-negative entries. 
Let  $\bfq^{\triangle}(w,r,s)=(a',b',c')$ where
\begin{equation}\label{phiprime}
\begin{cases}
 a'(w,r,s)&=\ds i-a= \displaystyle
 \frac{6(m\!-\!1)\!-\!w \!-\!2r  \!-\!3s\!-\!\epsilon(w\!+\!s)-2 \eta(w+r)}{2}
 \\
 b'(w,r,s)&=\ds j-b =\ds r + \epsilon(w+s), \\
 c'(w,r,s)&=\ds k-c= 
  \displaystyle
\frac{w-2r+3s+2\eta(w+r)-3\epsilon(w+s)+6}{6}.
\end{cases}
\end{equation}

Obviously, we have $b' \geq 0$. 
For $c'$, note that $w -2r \geq 3s$ on $\dtt$, hence
$ c' \geq s +1> 0$. In order to prove that $a' \geq 0$, we need the condition that $(w,r,s)$ satisfies
$\ttd ((w,r,s)) \geq \tau_w+1$. Using formula (\ref{residuemod6def}) to compute the residue modulo $6$ of $w-2r-3s$,
we compute
\begin{equation}\label{ttdexplicit}
\ttd ((w,r,s))=\floor{\frac{w-2r-3s}{6}}=\frac{w-3s-2r-\eta(w+r)-3\epsilon(w+s+\eta(w+r))}{6},
\end{equation}
hence the condition $\ttd ((w,r,s)) \geq \tau_w+1$ is equivalent to
\begin{equation}\label{tconditionwrs}
6m-12-w-2r-3s-\eta(w+r)-3\epsilon(w+s+\eta(w+r))+2 \eta(w) \geq 0.
\end{equation}
Since $2a'= 6(m\!-\!1)  -w-2r-3s-\epsilon(w+s)-2\eta(w+r),$ 
we conclude that for $(w,r,s) \in \dtt$ we have
\begin{equation}\label{2aprime2}
2a' \geq (4 -2\eta(w))+(2-\eta(w+r)+3\epsilon(w+s+\eta(w+r))-\epsilon(w+s)) \geq 4-2\eta(w)).
\end{equation}
Thus $a' \geq 2-\eta(w) \geq 0$, as desired. 
We will refine this inequality in Lemma \ref{aprime}.
\end{proof}

\bc \label{phitriangles}
The bijection $\varphi^{\triangle}:\rtt \ra \mtt$ of Theorem \ref{phitriregion} maps $\rtri$ onto $\mtri$. More precisely,
for a positive integer $w$ congruent to $2$ modulo $3$, let $\bfu_w= \left(6m-8-w, 0,  \ds  \frac{2w-6m+8}{3}\right)$ and
$\bfv_w=\left(6m-8-w, 0,  \ds  \frac{2w-3m+8}{3}\right)$. Then,
\begin{eqnarray*}
\rtt \setminus \rtri &=& \{\bfu_w: \, 3m-4 \leq w \leq 6m-10, \, w\equiv 2 \; (\mbox{mod. }3)\},
\\
\mtt \setminus \mtri &=& \{\bfv_w: \, 3m-4 \leq w \leq 6m-10, \, w\equiv 2 \; (\mbox{mod. }3)\},
\end{eqnarray*}
and $\bfv_w =\bfu_w+(0,0,m)=\varphi^{\triangle} (\bfu_w)$ for every $3m-4 \leq w \leq 6m-10$,  $w\equiv 2 \; (\mbox{mod. }3)$.
\ec

\begin{proof}
Recall that $\mtt$ is the set of monomials $\bfv=(i,j,k)$ such that
$i+j \leq 3(m-1)$, $i+j+k \geq 3m-2$ and $\ttm(\bfv) \geq \tau_w+1$.

Suppose that $\bfv=(i,j,k) $ belongs to $\mtt$ and has $\mscr$-weight equal to  $w$, that is, $i+2j+3k-3m=w$.
Then, $i \equiv j+w \pmod{3}$, hence $\eta(-i-j)=\eta(j-w)$ and the condition $\ttm(\bfv) \geq \tau_w+1$ -- see (\ref{tinvariants}) -- is equivalent to
$  i+j+\eta(j-w) \leq 6m-9-w+\eta(w)=\delta_w+\eta(w).$
Therefore, we have
\begin{equation}\label{iwdelta}
i\leq \delta_w+\eta(w)-j-\eta(j-w).
\end{equation}
Now observe that the right hand side of (\ref{iwdelta}) is $\leq \delta_w$ unless $j=0$ and $w \equiv 2 \pmod{3}$.
On the other hand, if $\bfv  \in \mtt \setminus \mtri$, then $i \geq \delta_w+1$. 
Comparing with (\ref{iwdelta}),
we conclude that $w \equiv 2 \pmod{3}$, $j=0$ and $i=\delta_w+1=6m-8-w$, therefore
$\bfv=\bfv_w $. 
As $0 \leq i \leq 3(m-1)$ and  $w \equiv 2 \pmod{3}$, we see that $ 3m-4 \leq w \leq 6m-10$. For these values of $w$, one checks immediately
that $\bfv_w \in \mtt$, hence
$$
\mtt \setminus \mtri = \{\bfv_w: \, 3m-4 \leq w \leq 6m-10, \, w\equiv 2 \; (\mbox{mod. }3)\}.
$$

Now suppose that $\bfu=(a,b,c) \in \rtt \setminus \rtri$, and let $\bfv=(i,j,k)= \varphi^{\triangle} (\bfu)$. 
Then,
$\bfv \in \mtt$ and $\wtm(\bfv)=\wtr(\bfu)$, hence $i \geq a$ by divisibility, and $i+w \geq a+ w \geq \delta_w+1$
because $\bfu \notin \rtri$.
It follows that    $\bfv \in \mtt \setminus \mtri$, and by the first part of the proof $w \equiv 2 \pmod{3}$, $ 3m-4 \leq w \leq 6m-10$ and $\bfv=\bfv_w$. As $i(\bfv_w)=6m-8-w=\delta_w+1$, it follows from the inequalities above that $a=i$.
As $j=0$, by divisibility we must also have $b=0$.  As $\wtr(\bfu)=w$, we conclude $\bfu=\bfu_w$.
\end{proof}

\subsection{Technical results on the triangular regions}
We will  need to parametrize monomials $\bfv=(i,j,k)$ in $\mtt_w$ using the first two entries $i$ and $j$
of $\bfv$. The following lemma identifies the { domain} of this parametrization.
\bl \label{ijtriangle}
Assume $3(m-1) \leq w \leq 6m-9$. Let $\delta_w=6m-9-w$. Recall
that 
$\mtt_w$ is the set of monomials $\bfv=(i,j,k)$ such that $i+2j+3k-3m=w$,
$i+j \leq 3(m-1)$, $i+j+k \geq 3m-2$ and $\ttm(\bfv) \geq \tau_w+1$.

There exists a $\bfv \in \mtt_w$ with $j(\bfv)=j$ if and only if $j \in \N$ satisfies
\begin{equation}\label{monewj}
j+2\eta(j-w) \leq \delta_w+2 \eta(w)-3.
\end{equation}
Given such an integer $j$, there exists a $\bfv \in \mtt_w$ with $j(\bfv)=j$ and $i(\bfv)=i$ if and only if
$i \in \N$ is congruent to $j+w \pmod{3}$ and satisfies

\begin{equation}\label{monewij}
\frac{\delta_w+3-j+3 \epsilon(w+j)}{2}
\leq i
\leq
\delta_w+\eta(w)-j-\eta(j-w).
\end{equation}
\el

\begin{proof}
Let $\bfv=(i,j,k) \in \mtt_w$. 
The upper bound in $(\ref{monewij})$ for $i$ is proven in Corollary \ref{phitriangles}, formula (\ref{iwdelta}).
To prove the lower bound for $i$, we use the condition $i+j+k \geq 3m-2$: solving  the relation
$i+2j+3k=m+3w$ for $k$,  we obtain $2i \geq \delta_w+3-j$. 
If $w+j$ is even, this gives the lower bound 
in $(\ref{monewij})$. 
If $w+j$ is odd, then $\delta_w+3-j=6(m-1)-w-j$ is odd, therefore $2i \geq \delta_w+3-j+1$.
We use the congruence $i \equiv j+w$ (modulo $3$) to conclude that $2i \neq  \delta_w+3-j+1$. 
Hence, $2i$
is at leas the next even number $\delta_w+3-j+3$, and this gives the lower bound in (\ref{monewij}) when $w+j$ is odd.

The interval for $i$ in (\ref{monewij}) is nonempty if and only if
\begin{equation}\label{inequmod6}
\delta_w+2 \eta(w)-3  -j-2\eta(j-w)-3 \epsilon(w+j) \geq 0
\end{equation}
Note that the the left hand side of (\ref{inequmod6}) is congruent to $0$ modulo $2$ and modulo $3$, hence
is congruent to $0$ also modulo $6$. We conclude that there exists an integer $i$ satisfying (\ref{monewij}) if and only
if $j+2\eta(j-w) \leq \delta_w+2 \eta(w)-3$.

Conversely, suppose $j$ satisfies (\ref {monewj}) and $i$ is an integer congruent to $j+w$ modulo $3$ satisfying  (\ref{monewij}) -
one such integer is $i=\delta_w+\eta(w)-j-\eta(j-w)$. 
Set $k=\ds \frac{w+3m-i-2j}{3}$
and $\bfv=(i,j,k)$. 
By reversing the above argument, we see that $i+j+k \geq 3m-2$ and
$\ttm(\bfv) \geq \tau_w+1$. 
To conclude that $\bfv \in \mtt_w$ we  need to check that $i+j \leq 3(m-1)$, but this follows from $\ttm(\bfv) \geq \tau_w+1$ by the assumption $w \geq 3(m-1)$:
we have
$
i+j \leq \delta_w=6m-9-w \leq 3m-6.
$\end{proof}

The proof or our main theorem  will require more details on the first coordinate $a'$ of the multiplier
$\bfq^\triangle(w,r,s)$. 
By Theorem \ref{phitriregion},  we know  that $a' \geq0$. The following Lemma refines this
 result.

\bl \label{aprime}
In the notation of Theorem \ref{phitriregion}, given $(w,r,s) \in \dtt$, write
$(a',b',c')=\bfq^{\triangle}(w,r,s)=\bfv^{\triangle}(w,r,s)- \bfu^{\triangle}(w,r,s)$.
Let $(a,b,c)= \bfu^{\triangle}(w,r,s)$ and $t=t_\dscr (w,r,s)=\floor{\frac{c}{2}}$. 
Then,
\begin{equation}\label{aprimetformula}
 a'=3 (t-\tau_w-1)+ h(a,c)+2-\eta(w)
\end{equation}
where
\begin{equation}\label{hformula}
h(a,c)=\frac{2+3\epsilon(c)-\eta(a)-\epsilon(c+\eta(a))}{2} \in\{0,1,2\}
\end{equation}
depends only on the residue classes $\eta(a)$ and $\epsilon(c)=\epsilon(w+a)$.
Thus $h(a,c)$ is determined by
$$h(1,0)=h(2,0)=0, \, h(0,0)=h(2,1)=1, \, h(0,1)=h(1,1)=2.$$
\el

\begin{proof}
Formula (\ref{phiprime}) for $a'$ gives
$a'(w,r,s)=\displaystyle
 \frac{6(m\!-\!1)\!-\!w \!-\!2r  \!-\!3s\!-\!\epsilon(w\!+\!s)-2 \eta(w+r)}{2}$.
On the other hand, we have $t=\ds \floor{\frac{w-2r-3s}{6}}$, hence by (\ref{residuemod6def}) we have
$
6t=w-2r-3s-\eta(w+r)-3\epsilon(w+s+\eta(w+r)).
$
By definition of  threshold number $\tau_w$, we have
$
3 \tau_w =w-\eta(w)-3(m-1),
$
hence
$$
2 a'-6(t-\tau_w-1)=6-\epsilon(w+s)-\eta(w+r)-2\eta(w)+3\epsilon(w+s+\eta(w+r)).
$$
Finally, observe that $\epsilon(w+a)=\epsilon(2a+2b+3c)=\epsilon(c)$ and
$\eta(w+r)=\eta(a+3b+3c)=\eta(a)$, while
$$
\epsilon(w+s)=\epsilon(3w+3s)=\epsilon(w+a+\eta(a))
\quad \text{and} \quad
\epsilon(w+s + \eta(w+r))=\epsilon(w+a+2\eta(a))=\epsilon(w+a).
$$
\end{proof}

\subsection{Uniqueness in the triangular regions} \label{utriangle}
We can now prove uniqueness for the bijection $\varphi^\triangle$  on the triangular regions:
\bt \label{unitriangle}
Let $0 \leq w \leq 6m\!-\!9$ and
let $\varphi^\triangle:\rtriw \ra \mtriw$ be the restriction to $\rtriw$ of the bijection $\varphi^\triangle$ of Theorem \ref{phitriregion}. 
If
$\hat{\varphi}:\rtriw \ra \mtriw$ is a bijection that has the same multiset of multipliers as $\varphi^\triangle$, 
then $\hat{\varphi}=\varphi^\triangle$.
\et

\begin{proof}
We fix $w$ throughout the proof.
 We write
$\bfu=(a,b,c)$ for a monomial in $\rscr_w$, $\bfv=(i,j,k)$ for a monomial in $\mscr_w$, $\bfq=(a',b',c')$ for a multiplier in $\mathscr{R}_{3m}$.

Let $D=\left(\bfu^\triangle\right)^{-1}(\rtriw) = \left(\bfv^\triangle\right)^{-1}  (\mtriw)$.
Since $w$ is fixed,
$D$ is the set of those $(r,s)$ that parametrize monomials in $\rtriw$ and $\mtriw$ as
$\bfu^\triangle(w,r,s)$ and $\bfv^\triangle(w,r,s)$ respectively. 
To further simplify notation,
we write $\bfu(r,s)$ for $\bfu^\triangle(w,r,s)$, $\bfv(r,s)$ for $\bfv^\triangle(w,r,s)$ and $\varphi$ for $\varphi^\triangle$: explicit formulas for these parametrizations
are given in Theorem \ref{phitriregion}. Note that $\bfv(r,s)=\varphi(\bfu(r,s))$ by definition of $\varphi$.
The $\varphi$-multiplier of $\bfu(r,s)$ is $\bfq(r,s)= \bfv(r,s)-\bfu(r,s)$. 
Similarly, given a bijection
$\hat{\varphi}: \rtriw \ra \mtriw$, we write $\hat{\bfv}(r,s)=\hat{\varphi}(\bfu(r,s))$ and $\hat{\bfq}(r,s)= \hat{\bfv}(r,s)-\bfu(r,s)$.
The assumption that $\hat{\varphi}$ has the same multiset of multipliers as $\varphi$ means that
for every $\bfq \in \mathscr{R}_{3m}$, the number of pairs $(r,s) \in D$ such that $\bfq(r,s)=\bfq$ equals the number
of pairs  $(r_0,s_0) \in D$ such that  $\hat{\bfq}(r_0,s_0)=\bfq$. If this is the case, we have to show that $\hat{\varphi}=\varphi$, that is,
$\hat{\bfv}(r,s)=\bfv(r,s)$ for every $(r,s) \in D$.

We will prove this by descending induction on the first coordinate $a(r,s)$ of $\bfu(r,s)$, using Claim 1 and Claim 2 below. 
Recall that $a(r,s)=3s+\eta(w+r)$.
\\

\noindent {\bf Claim 1:}  The bijection $\hat{\varphi}$ preserves the invariant $r$, that is, if we define for an integer $r$
$$
\scrR^{\triangle}_{w,r}= \big\{\bfu \in \rtriw: b(\mf{u})=r\big\}, \quad
\scrM^{\triangle}_{w,r}= \big\{\bfv \in \mtriw:  \floor{j(\mf{v})/2}=r\big\},
$$
then $\hat{\varphi} (\scrR^{\triangle}_{w,r})=\scrM^{\triangle}_{w,r}$.
\\

\noindent {\bf Proof of Claim 1.}
We proceed by induction on $r$. Given an integer $r$, define
$$
\npr=\#\big\{\bfv \in \scrM^{\triangle}_{w,r}: j(\mf{v})=2r\big\}, \quad
\nmr=\#\big\{\bfv \in \scrM^{\triangle}_{w,r}: j(\mf{v})=2r+1\big\}.
$$
Note that $\# \scrR^{\triangle}_{w,r}= \# \scrM^{\triangle}_{w,r}= \npr+\nmr$
because $\varphi$ preserves the $r$-invariant.
By the explicit formula (\ref{p'map-explicit}) for $\bfv(r,s)$, the number
$\npr$ (resp. $\nmr$) is the number of integers $s$ such that $(r,s)$ belongs to $D$ and
$\epsilon (w+s) =0$ ( resp. $\epsilon (w+s) =1$).

By  formula (\ref{phiprime}) for the $\varphi$-multipliers, we see that, given an integer $h$,
the number of $\varphi$-multipliers  $\bfq=(a',b',c')$ such that $b'=h$ is equal to
$\nu^{+}_{h}+\nu^{-}_{h\!-\!1}$ (counting multiplicities).

We will prove slightly more than Claim 1: by induction on $r$ we will prove
\\

\noindent {\bf Statement $S_r$}:
Given  $\bfv$ in $\scrM^{\triangle}_{w,r}$, let $\hat{\bfu}(\bfv)=\hat{\varphi}^{-1}(\bfv)=(a,b,c)$.
Then, $b=r$, hence  $\hat{\varphi} \big(\scrR^{\triangle}_{w,r}\big)=\scrM^{\triangle}_{w,r}$.
Furthermore, if $\hat{\bfq}(\bfv)=\bfv-\hat{\bfu}(\bfv)=(a',b',c')$ denotes the corresponding multiplier,
then  $b'$ is equal to $r$ for precisely $\npr$ monomials $\bfv$ in $\scrM^{\triangle}_{w,r}$ and to $r+1$ for the remaining $\nmr$ monomials in $\scrM^{\triangle}_{w,r}$.
\\

\noindent
The statement is trivial for $r<0$, so we can assume that $S_{r'}$ holds for all $r'\leq r$ and prove that $S_{r+1}$ holds. We use the hypothesis that $\varphi$ and $\hat{\varphi}$ have the same multiset of multipliers:
by the induction hypothesis, all multipliers $\bfq=(a',b',c')$ with $b' \leq r$
occur as  $\hat{\bfq}(\bfv)$ for some $\bfv \in \scrM^{\triangle}_{w,r'}$ with $r' \leq r$, therefore
$b' \geq r+1$ if $(a',b',c')=\hat{\bfq}(\bfv)$ with $\bfv \in \scrM^{\triangle}_{w,r+1}$. Since the
multipliers $\bfq=(a',b',c')$ satisfying
$b'=r+1$ are exactly $\nmr+\nu^{+}_{r+1}$, the induction hypothesis implies also that
$b'\big(\hat{\bfq}(\bfv)\big) = r+1$ for at most $\nu^{+}_{r+1}$ monomials $\bfv$ in $\scrM^{\triangle}_{w,r+1}$.

There are $\nu^{+}_{r+1}$ monomials $\bfv=(i,j,k)$ in $\scrM^{\triangle}_{w,r+1}$ with $j=2(r+1)$. Let
$\bfv$ be one of these and let $\hat{\bfu}=\hat{\varphi}^{-1}(\bfv)=(a,b,c)$, $\hat{\bfq}=\bfv-\hat{\bfu}=(a',b',c')$.
By  induction, we have $b\geq r+1$ and $b'  \geq r+1$. 
On the other hand,
 $b+b'=j=2(r+1)$ implies
$b= r+1$ and $b' = r+1$. 
This accounts for all the $\nu^{+}_{r+1}$
remaining multipliers $\bfq$ with $b'=r+1$. So for the $\nu^{-}_{r+1}$ monomials $\bfv$ in $\scrM^{\triangle}_{w,r+1}$ with $j(\bfv)=2(r+1)+1$,
we must have $b'=b'\big(\hat{\bfq}(\bfv)\big)  \geq r+2$. As $j=b+b'$, we conclude that
$b= r+1$ and $b'  = r+2$ for the monomials $\bfv$ in $\scrM^{\triangle}_{w,r+1}$ with $j=2(r+1)+1$. This proves $S_{r+1}$ and Claim 1.
\\

\noindent {\bf Claim 2:}
If $s>\bar{s}$ and $\bfv(r,\bar{s})-\bfu(r,s)=\bfq(r_1,s_1)$ for some $(r_1,s_1)\in D$,
then $a(r_1,s_1)>a(r,s)$.
\\

\noindent {\bf Proof of Claim 2.}
Suppose that  $\bfv(r,\bar{s})-\bfu(r,s)=\bfq(r_1,s_1)$: then
\begin{equation}\label{deltaq}
\bfq(r_1,s_1)=\bfv(r,\bar{s})-\bfu(r,s)=\bfv(r,\bar{s})+\bfq(r,s)-\bfv(r,s).
\end{equation}
We can compute
the first entry $\partial a'=a'(r_1,s_1)-a(r,s)$ of $\bfq(r_1,s_1)-\bfq(r,s)$ using formula (\ref{phiprime}):
\begin{equation}\label{deltaaprime}
2 \partial a'= -2(r_1-r)-3(s_1-s)- (\epsilon(w+s_1)- \epsilon(w+s))- 2 ( \eta(w+r_1)-\eta(w+r)).
\end{equation}

\noindent{\bf Case 1.} Assume first that $\bar{s}=s-2h$ with $h \geq 1$. By (\ref{p'map-explicit}) we have
\begin{equation}\label{deltavcase1}
\bfq(r_1,s_1)-\bfq(r,s) =  \bfv(r,s-2h)-\bfv(r,s)=(-3h,0,h),
\end{equation}
hence $2 \partial a'= -6h$ and $b'(\bfq(r_1,s_1))=b'(\bfq(r,s))$, that is, $r_1+\epsilon(w+s_1)=r+\epsilon(w+s)$.

\noindent{\bf Case 1.1.} Suppose further that $\epsilon(w+s_1)=\epsilon(w+s)$. Then, $r_1=r$,
and $2\partial a'= 3(s-s_1)$ by  (\ref{deltaaprime}).
On the other hand, $\partial a'=-3h$  by (\ref{deltavcase1}), and we conclude $s_1-s=2h \geq 2$.
Hence $s_1 >s$ and a fortiori  $a(r_1,s_1)>a(r,s)$.

\noindent{\bf Case 1.2.} Suppose now that $\epsilon(w+s_1)=0$ and $\epsilon(w+s)=1$. Then, $r_1=r+1$, so that $\eta(w+r_1)\equiv \eta(w+r)+1 $
$\pmod{3}$ and $\eta(w+r)-\eta(w+r_1)\geq -1$. Using formulas (\ref{deltaaprime}) and (\ref{deltavcase1}) we compute
$$
3(s_1-s)=6h+2(\eta(w+r)-\eta(w+r_1))-1 \geq 3,
$$
hence $s_1 >s$ and a fortiori $a(r_1,s_1)>a(r,s)$.

\noindent{\bf Case 1.3.} Suppose finally that $\epsilon(w+s_1)=1$ and $\epsilon(w+s)=0$. Then, $r_1=r-1$, and $\eta(w+r)-\eta(w+r_1)\geq -2$.  Using formulas (\ref{deltaaprime}) and (\ref{deltavcase1}) we compute
$$
3(s_1-s)=6h+2(\eta(w+r)-\eta(w+r_1))+1 \geq 3
$$
hence $s_1 >s$ and a fortiori $a(r_1,s_1)>a(r,s)$.

\noindent{\bf Case 2.} Assume now that $\bar{s}=s-2h-1$ with $h \geq 0$, and $\epsilon(w+s)=0$. As in Case 1, we can use formula (\ref{deltaq}) to conclude
\begin{equation}\label{deltavcase2}
\bfq(r_1,s_1)-\bfq(r,s) =  \bfv(r,s-2h-1)-\bfv(r,s)=(-2-3h,1,h),
\end{equation}
hence  $2 \partial a'= -6h-4$ and $b'(\bfq(r_1,s_1))=b'(\bfq(r,s))+1$, that is, $r_1+\epsilon(w+s_1)=r+1$.

\noindent{\bf Case 2.1.} Suppose further that $\epsilon(w+s_1)=0$. Then, $r_1=r+1$, $\eta(w+r_1) \equiv \eta(w+r)+1 \pmod{3}$ and
$$
3(s_1-s)=6h+2(\eta(w+r)-\eta(w+r_1))+2 \geq 0.
$$
Thus $s_1 \geq s$. If $s_1>s$, a fortiori $a(r_1,s_1)>a(r,s)$. If $s_1=s$, then  $\eta(w+r)-\eta(w+r_1)=-1$ and so
$$
a(r_1,s_1=s)=3s+\eta(w+r_1)>3s+\eta(w+r)= a(r,s).
$$

\noindent{\bf Case 2.2.} Suppose instead that $\epsilon(w+s_1)=1$. Then, $r_1=r$ and
$
3(s_1-s)=6h+3 \geq 3.
$
Thus $s_1 > s$ and a fortiori $a(r_1,s_1)>a(r,s)$.

\noindent{\bf Case 3.} The last possibility is that $\bar{s}=s-2h-1$ with $h \geq 0$ and $\epsilon(w+s)=1$. 
In this case
\begin{equation}\label{deltavcase3}
\bfq(r_1,s_1)-\bfq(r,s) =  \bfv(r,s-2h-1)-\bfv(r,s)=(-1-3h,-1,h+1),
\end{equation}
hence $2 \partial a'= -6h-2$ and $b'(\bfq(r_1,s_1))=b'(\bfq(r,s))-1$, that is, $r_1+\epsilon(w+s_1)=r$.

\noindent{\bf Case 3.1.} Suppose furthermore that $\epsilon(w+s_1)=0$. Then, $r_1=r$ and
$
3(s_1-s)=6h+3 \geq 3.
$
Thus $s_1 > s$ and a fortiori $a(r_1,s_1)>a(r,s)$.

\noindent{\bf Case 3.2.} Suppose finally that $\epsilon(w+s_1)=1$. Then, $r_1=r-1$ and
$
3(s_1-s)=6h+2(\eta(w+r)-\eta(w+r_1))+4 \geq 0.
$
Thus $s_1 \geq s$. If $s_1>s$, a fortiori $a(r_1,s_1)>a(r,s)$. If $s_1=s$, then we must have $\eta(w+r)-\eta(w+r_1)=-2$ and so
$
a(r_1,s_1=s)=3s+\eta(w+r_1)>3s+\eta(w+r)=a(r,s).
$
This completes the  proof of Claim 2.
\\

\noindent {\bf End of the proof of Theorem \ref{unitriangle}.}
Fix a $(r,s) \in D$, 
we have to prove that $\hat{\bfv}(r,s)=\bfv(r,s)$.
By Claim 1 we can find $\bar{s}$ such that $\hat{\bfv} (r,s)=\bfv(r,\bar{s})$.
Consider the $\hat{\varphi}$-multiplier
$$\hat{\bfq}= \hat{\bfq}(r,s)=
\hat{\bfv}(r,s)-\bfu(r,s)
=\bfv(r,\bar{s})-\bfu(r,s).
$$
By assumption $\hat{\bfq}$ is also a $\varphi$-multiplier, hence there is a pair $(r_1,s_1) \in D$ such that
$\hat{\bfq}=\bfq(r_1,s_1)$.
We will prove by descending induction on $a=a(r,s)$ that $\bar{s}=s$, hence $\hat{\bfv}(r,s)=\bfv(r,s)$ as desired.

Suppose first
$a(r,s)$ is the largest value of $a$ on $D$. We must have $s \geq \bar{s}$ because   $a(r,s) \geq a(r,\bar{s})$
by the maximality of $a(r,s)$. 
On the other hand, if $s>\bar{s}$, then Claim 2 implies $a(r_1,s_1)>a(r,s)$, a contradiction.
We conclude that $\bar{s} = s$.

Pick now an arbitrary $(r,s)$.  If $\bar{s}>s$, then $a(r,\bar{s})>a(r,s)$ and by induction
$\hat{\bfv}(r,\bar{s})=\bfv (r,\bar{s})=\hat{\bfv}(r,s)$.
Since $\hat{\bfv}$ is bijective, it follows that $\bar{s}=s$,
contradiction, hence $\bar{s} \leq s$.
Suppose that $\bar{s} < s$.
Let $D_1$ the set of pairs $(r_1,s_1) \in D$
such that $\hat{\bfq}=\bfq(r_1,s_1)$, and let $m_1$ the cardinality of $D_1$: then $\hat{\bfq}$ has
multiplicity $m_1$ as a $\varphi$-multiplier, hence also as a $\hat{\varphi}$-multiplier.
 By Claim 2,
for any pair $(r_1,s_1) \in D_1$ the inequality $a(r_1,s_1)>a(r,s)$ holds - in particular, $(r,s) \notin D_1$.
The descending induction hypothesis implies that
 $\hat{\bfv}(r_1,s_1)=\bfv(r_1,s_1)$, hence $\hat{\bfq}(r_1,s_1)= \bfq(r_1,s_1)$, for all pairs $(r_1,s_1) \in D_1$.
This implies that the $\hat{\varphi}$-multiplier $\hat{\bfq}$ appears already with its multiplicity $m_1$
as the $\hat{\varphi}$-multiplier associated to pairs  $(r_1,s_1) \in D_1$. 
Hence, it cannot be equal to $\hat{\bfq}(r_0,s_0)$
for any pair $(r_0,s_0)$ that does not belong to $D_1$.
This contradicts
$\hat{\bfq}= \hat{\bfq}(r,s)$ because $(r,s) \notin D_1$. Hence, $\bar{s}=s$ and the proof is complete.
\end{proof}

\section{The bijection $\varphi: \rscr'_w \ra \mscr'_w$ and its uniqueness for $w \equiv 0$ modulo $3$} \label{wcongruozero}
\index{aavphiz @ $\varphi$ the bijection $\rscr' \ra \mscr'$ }
The region $\rscr'$  is the disjoint union
$\rscr'= \rtri \cup \rect$ of the triangular and the rectangular regions, which are distinguished by the $\ttr$-invariant, and similarly $\mscr'= \mtri \cup \mrect$.
We define a bijection $\varphi: \rscr' \ra \mscr'$ as $\varphi^\triangle$ on $\rtri$ (cf. Theorem \ref{phitriregion}) and as $\varphi^\hrectangle$ on $\rect$  (cf. Theorem \ref{phirectregion}). 
In particular,
$\varphi(\rtri)=\mtri$ and $\varphi(\rect)=\mrect$. 
By Theorem \ref{phitriregion} and \ref{phirectregion}, the bijection $\varphi$ satisfies divisibility, that is,
for each $\bfu \in \rscr'$ the multiplier $\bfq(\bfu)=\varphi(\bfu)-\bfu$ has non-negative entries.

In this section we prove {\em uniqueness}, in the sense of page \pageref{divis-uniq},
for the restriction of $\varphi$ to $\rscr'_w$ in the case $w \equiv 0 \pmod{3}$.

\bt \label{uniw0}
Assume $w \equiv 0 \pmod{3}$ and $0 \leq w \leq 6m-9$. Let  $\varphi: \rscr'_w
 \ra \mscr'_w$
be the bijection defined as $\varphi^\triangle$ on $\rtri_w$ and as $\varphi^\hrectangle$ on $\rect_w$.
Then, $\varphi$ satisfies {\em uniqueness}:
if $\hat{\varphi}: \rscr'_w \ra \mscr'_w$ is a bijection with the same multiset of multipliers as $\varphi_w$,
then $\hat{\varphi}= \varphi$.
\et

\begin{figure}[ht]
    \centering
    \includegraphics[width=0.8\textwidth]{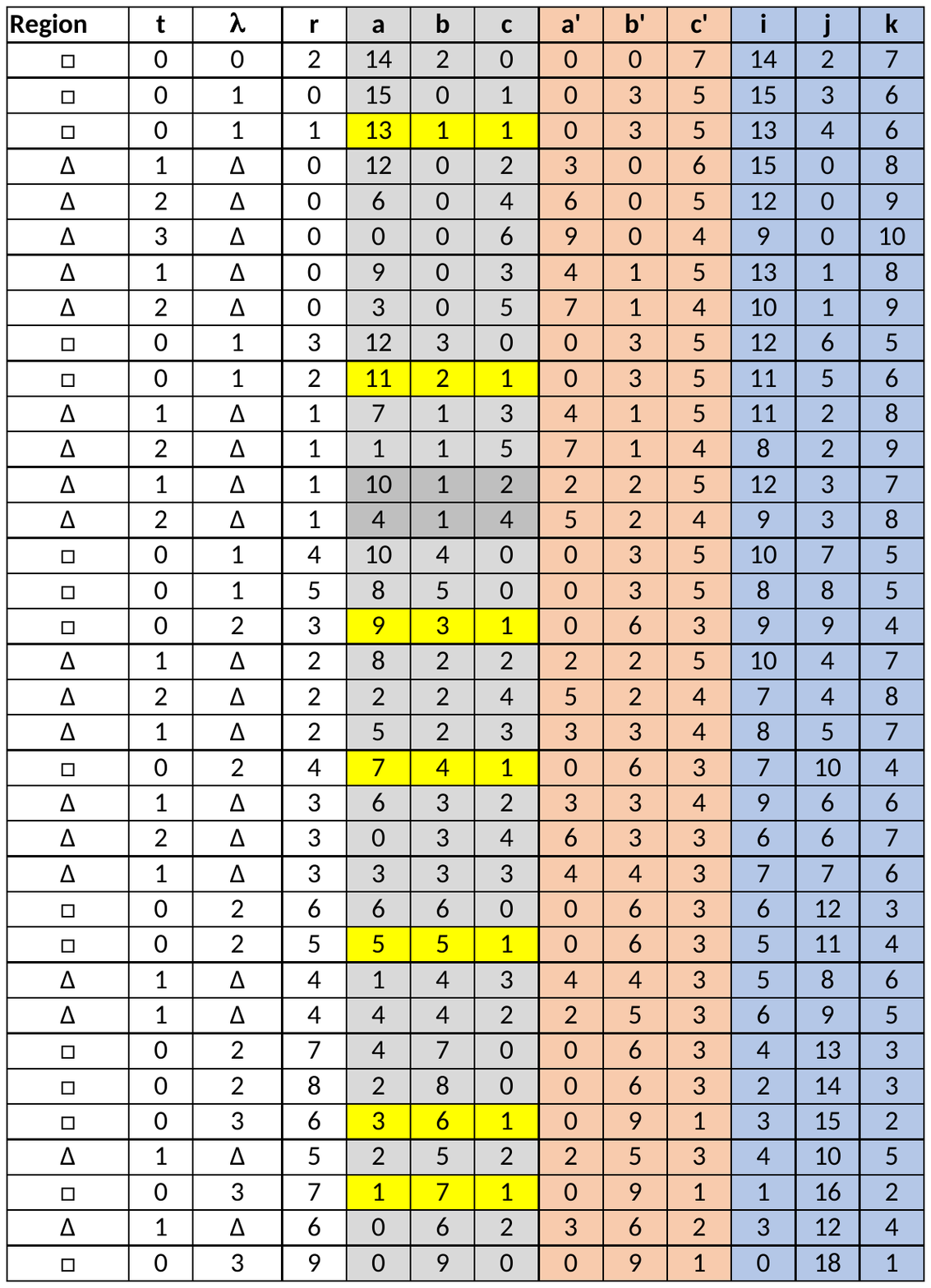}
    \caption{The map $\vphi$ for $m=7$ and $w=18$}
    \label{fig:m7w18}
\end{figure}

\begin{figure}[ht]
    \centering
    \includegraphics[width=0.8\textwidth]{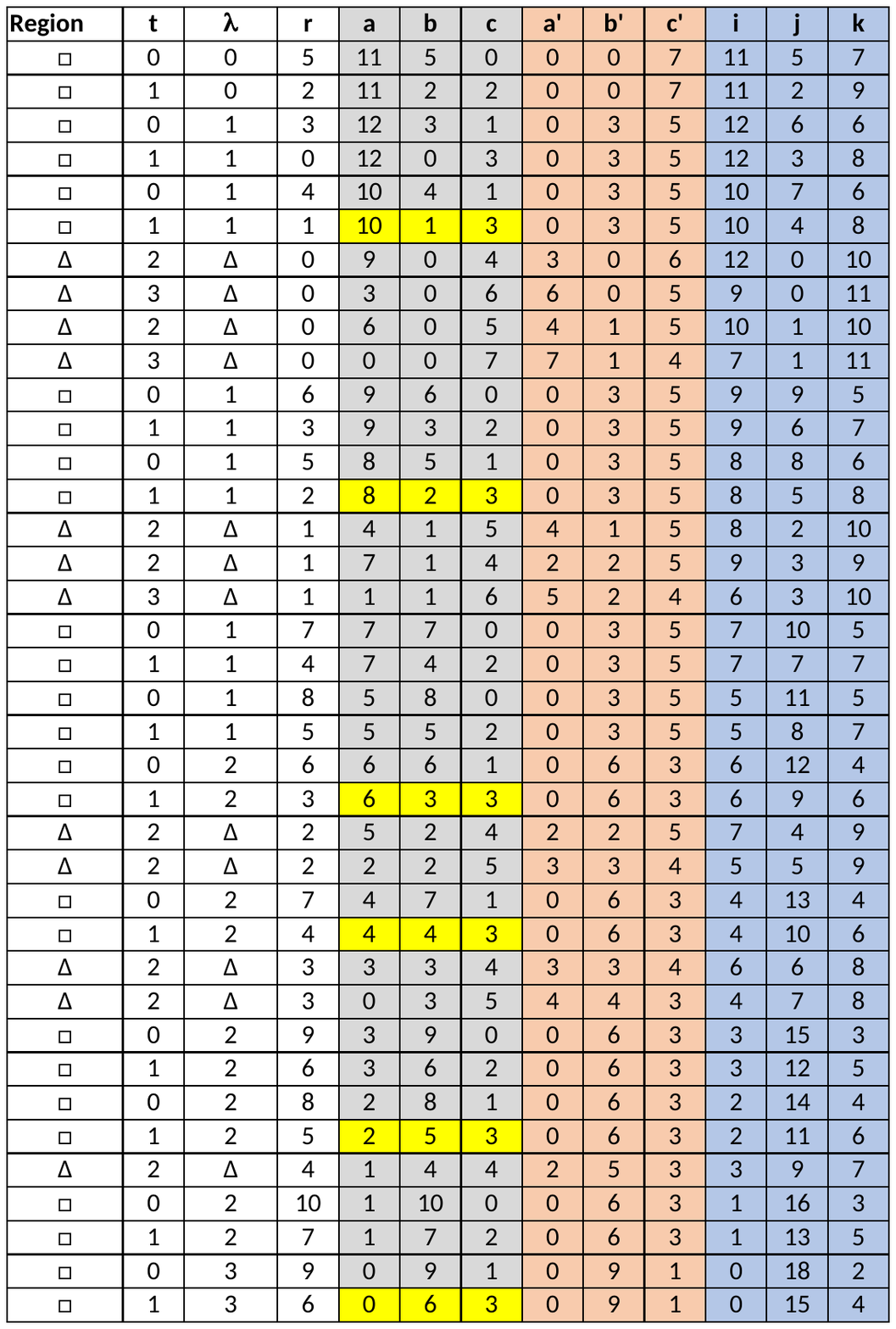}
    \caption{The map $\vphi$ for $m=7$ and $w=21$}
    \label{fig:m7w21}
\end{figure}

\begin{proof}

For $ w \leq 3m-4$ the rectangular regions are empty, hence uniqueness follows from Theorem \ref{unitriangle}. 
Thus, we may assume  $3(m-1) \leq w \leq 6m-9$.

We use the same notation as in Theorem \ref{unitriangle}: we write
$\bfu=(a,b,c)$ for a monomial in $\rscr_w$, $\bfv=(i,j,k)$ for a monomial in $\mscr_w$, $\bfq=(a',b',c')$ for multipliers in $\mathscr{R}_{3m}$.
Given $\bfu \in \rscr'_w$, we let $\bfq(\bfu)=\varphi(\bfu)-\bfu$ and $\hat{\bfq}(\bfu)=\hat{\varphi}(\bfu)-\bfu$
denote the corresponding multipliers. 
The hypothesis on multipliers is that for each $\bfq \in \mathscr{R}_{3m}$ the numbers
of monomials $\bfu \in \rscr'_w$ such that $\bfq(\bfu)=\bfq$ equals the number of  monomials $\bfu \in \rscr'_w$ such that $\hat{\bfq}(\bfu)=\bfq$.
We will say that $\bfq=(a',b',c')$ is a {\em rectangular (resp. triangular) region multiplier}
 if it is of the form $\bfq(\bfu)$ with $\bfu \in \rect_w$ (resp. $\bfu \in \rtri_w$). 
By Theorem \ref{phirectregion}, 
rectangular region multipliers have the form $\bfq=(0, 3\lambda, m -2 \lambda)$ for some integer $\lambda$; 
in particular, they have $a'=0$, while triangular region multipliers have $a' \geq 2$ by Lemma \ref{aprime}, because $w \equiv 0 \pmod{3}$: the sets of rectangular and triangular region multipliers are disjoint, and $a' \neq 1$ for all multipliers  (this is why the case $w \equiv 0 \pmod{3}$ is easier
than the cases $w \equiv 1,2$).

Recall that $\delta_w=6m-w-9$, and set $r_{max}=\floor{\frac{\delta_w}{2}}$. 
We will use the ordering $\prec$ 
of Definition \ref{order-rect}
for the monomials in the rectangular region, and the monomials $\bfu^{(r)}$ introduced in Lemma \ref{urrectangular}.

We  prove that the following  statements hold.

\noindent {\bf Claim} $S_{r,1}$ for $0 \leq r \leq r_{max}$:
$\hat{\varphi}(\bfu)=\varphi(\bfu) \quad \mbox{for every $\bfu \in \rectw$ satisfying $ \bfu \preceq \bfu^{(r)}$}$

\smallskip
\noindent {\bf Claim} $S_{r,2}$  for $0 \leq r \leq r_{max}\!-\!1$:
Denote by $\scrR^{\triangle}_{w,r}$  the set of monomials $\bfu =(a,b,c) \in \rtriw$ for which $b=r$ , and by
$\scrM^{\triangle}_{w,r}$ the set of monomials $\bfv =(i,j,k) \in \mtriw$ for which $j \in \{2r,2r+1\}$.
Then,
$$
\hat{\varphi}\left (\scrR^{\triangle}_{w,r} \right) = \varphi \left (\scrR^{\triangle}_{w,r} \right) = \scrM^{\triangle}_{w,r}
$$
Moreover,  for $\bfu \in \scrR^{\triangle}_{w,r}$,  $\hat{q}(\bfu)=(a',b',c')$ is a triangular region multiplier,
 and $b'$ is equal to $r$ for precisely $\npr$ monomials $\bfu$ in $\scrR^{\triangle}_{w,r}$ and to $r+1$ for the remaining $\nmr$ monomials in $\scrR^{\triangle}_{w,r}$.

\smallskip

\noindent If statements $S_{r,1}$ and $S_{r-1,2}$ hold for $0 \leq r \leq r_{max}$, then the statement of the Theorem follows. To see why, we
observe  that all monomials $\bfu$ in the triangular region $\rtriw$ satisfy $b(\bfu) \leq r_{max}-1$. 
Indeed,
suppose $\bfu=(a,b,c) \in \rtriw$ and let $\bfv=\varphi(\bfu)=(i,j,k)$. 
Then, $j\in\{2b,2b+1\}$ and
inequality (\ref{monewj}) in Proposition \ref{ijtriangle} implies
$$
2b \leq j \leq \delta_w-3.
$$
This implies that $b \leq r_{max}-1$. 
Therefore, all monomials in the triangular region are covered by the statements
$S_{r,2}$.
 In particular, $\hat{\varphi}$ induces a bijection $\rtri \ra \mtri$ with the same multiset of multipliers
as the restriction of $\varphi$ to $\rtri$. 
By Theorem \ref{unitriangle}, it follows that
$\hat{\varphi}(\bfu)=\varphi(\bfu)$ for every $\bfu \in \rtri$. This forces $\hat{\varphi}$ to map $\rectw$ onto $\mrectw$ with the
same multiplier multiset as the restriction $\varphi^\hrectangle$ of $\varphi$ to $\rectw$. 
By Corollary \ref{unirectangle2}, the two bijections coincide also on the rectangular regions.

It remains to prove the statements $S_{r,1}$ and $S_{r-1,2}$ for $0 \leq r \leq r_{max}$, and we do this by induction on $r$.
First, we
establish the initial case $S_{0,1}$ ($S_{-1,2}$ is vacuous). 
Then, we show that, if $S_{k,1}$ and $S_{k-1,2}$
hold for $k \leq r$, then they also hold for $k=r+1$.

Recall that for $\bfu=(a,b,c) \in \rect_w$, we have defined $\lambda(\bfu)= \lambda(\delta_w\!-\!a)$.
We will need the monomials $\bfu^{(r)}=(\delta_w-2r, r,  2\tau_w +1)$ of  Lemma \ref{urrectangular} ,
which satisfy $\lambda (\bfu^{(r)})=\floor{\frac{r+3}{3}}$ and
\begin{equation}\label{blbound}
   b (\bfu) \geq r+1 \quad \mbox{if $\bfu \in \rectw$ and $\bfu \succ \bfu^{(r)}$}.
\end{equation}

\noindent {\bf The initial case}: $r=0$.

We first verify  $S_{0,1}$, that is, $\hat{\varphi}(\bfu) =\varphi(\bfu)$ for all $\bfu =(a,b,c)\in \rectw$
with $ \bfu \preceq \bfu^{(0)}=(\delta_w, 0,  2\tau_w +1)$. We proceed by induction on  $\prec$. Tables \ref{fig:m7w18}
and \ref{fig:m7w21} should help following the argument.

The first batch of monomials $\bfu$ in $\rectw$ with respect to the order $\prec$ are those for which $\lambda=0$, and these all have $a=\delta_w-1$. Since there are no multipliers with $a'=1$ and $i(\bfv) \leq \delta_w$ for every $\bfv \in \mscr'_w$, for any such $\bfu$ we must have $i(\hat{\varphi}(\bfu))=\delta_w-1$. 
By Proposition \ref{ijtriangle}, there is no
$\bfv \in \mtriw$ satisfying $i(\bfv)=\delta_w-1$ (because  such a $\bfv$ would have $j=0$ contradicting $i \equiv j$ modulo $3$).
Hence,  $\hat{\varphi}(\bfu) \in \mrectw$. 
Note that, if $\hbfq=(a',b',c')=\hphi(\bfu)-\bfu$, then $a'=i(\hat{\varphi}(\bfu))-a(\bfu)=0$,
hence $\hbfq=\bfq_n$ for some $n \geq 0=\lambda(\bfu)$. We then  conclude $\hat{\varphi} =\varphi$ on this batch of monomials using induction
on the order $\prec$ of $\rectw$ as in the proof of Corollary  \ref{unirectangle}. 
We call this the {\em standard argument for the rectangular region}.

The  other monomials $\bfu =(a,b,c) \in \rectw$ that satisfy $ \bfu \preceq \bfu^{(0)}$ are those  for which $\lambda=1$ and
$a\geq a(\bfu^{(0)})=\delta_w $. 
But we have $a\leq \delta_w$ for all $\bfu \in \scrR'_w$, hence $a= \delta_w$.
For any such $\bfu$ we must have $i(\hat{\varphi}(\bfu))=\delta_w$. By Proposition \ref{ijtriangle}, if
$\bfv \in \mtriw$ satisfies $i(\bfv)=\delta_w$, then $j(\bfv)=0$ and so $\bfv=(\delta_w,0,\frac{w-\delta_w}{3}+m)$.
So if $\bfv =\hat{\varphi}(\bfu)$ were in $\mtriw$, the corresponding multiplier $\hat{\bfq}=\bfv-\bfu$ would be
$(0,0,m)$, but this multiplier is no longer available because it has already been assigned to the monomials with $\lambda=0$,
which account for all the occurrences of the multiplier $(0,0,m)$.
Thus, $\hat{\varphi}(\bfu) \in \mrectw$ for all monomials $\bfu$ in $\rectw$ for which  $a=\delta_w$.
The standard argument for the rectangular region implies now that  $\hat{\varphi}=\varphi$
 for all monomials $\bfu$ in $\rectw$ for which $a=\delta_w$.

This concludes the proof of $S_{0,1}$, while $S_{-1,2}$ is vacuous.

\noindent {\bf The induction step}: assume  $S_{\ell,1}$ and $S_{\ell-1,2}$ hold for $\ell \leq r$. We will show that
$S_{r+1,1}$ and $S_{r,2}$ also hold. Note that the induction hypothesis implies that  all triangular region multipliers $\bfq =(a',b',c')$ with $b' \leq r-1$
arise as $\hat{\bfq}(\bfu)$ for some $\bfu$ in $\rtri_{w,\ell}$ for some $\ell \leq r-1$.  Furthermore, the rectangular region multipliers have the form
$\bfq_\lambda=(0, 3\lambda, m -2 \lambda)$ for some integer $\lambda$,  and those with $\lambda < \lambda (\bfu^{(r)}) = \floor{\frac{r+3}{3}}$
arise, by the induction hypothesis, as $\hat{\bfq}(\bfu)$ for some $\bfu$  in $\rectw$ satisfying $\bfu \prec \bfu^{(r)}$; the others have
$b' =3 \lambda \geq r+1$
because $\lambda \geq \lambda (\bfu^{(r)}) = \floor{\frac{r+3}{3}}$. 
Thus, the induction hypothesis implies $b'\geq r$ for all multipliers left, that is,
for those that can arise as $\hat{\bfq}(\bfu)$ for some $\bfu$ either in $\rectw$ with $\bfu \succ \bfu^{(r)}$
or in $\rtri_{w,\ell}$ with $\ell \geq r$.

We proceed to proving $S_{r+1,1}$. 
Since $S_{r,1}$ holds, we need to show that, if
$\bfu=(a,b,c) \in \rectw$ and  $\bfu^{(r)} \prec \bfu \preceq \bfu^{(r+1)}$, then
 $\hat{\varphi}(\bfu) =\varphi(\bfu)$. Note that
\begin{equation}\label{lambdarvalues}
\floor{\frac{r+3}{3}} =\lambda (\bfu^{(r)}) \leq \lambda (\bfu) \leq \lambda (\bfu^{(r+1)}) = \floor{\frac{r+4}{3}}
\end{equation}

\noindent {\bf Case 1}: $r \equiv 0 \pmod{3}$, say $r=3h$.

\noindent
In this case, the three monomials $\bfu^{(r)}$, $\bfu$ and $\bfu^{(r+1)}$ all have the same $\lambda=h+1$.
By Lemma \ref{aneworder} (2), the first entry of $\bfu$ is $a=\delta_w-2r-2$.

Now let $\bfv=\hat{\varphi}(\bfu)=(i,j,k)$: then $i \geq a=\delta_w -2r-2$.   
By Lemma \ref{urrectangular},
the monomial $\bfu$ has second entry $b\geq r+1$, while, by the induction hypothesis, the multiplier $(a',b',c')=\hat{\bfq}(\bfu)=\bfv-\bfu$ satisfies $b' \geq r$. Hence
\begin{equation} \label{jbbineq}
j=b+b' \geq (r+1)+r=2r+1.
\end{equation}
We claim that this implies that $\bfv \in \mrectw$. 
Suppose by  contradiction that $\bfv \in \mtriw$. 
Then, by Proposition \ref{ijtriangle}, we have
\begin{equation} \label{diineq}
i \leq \delta_w -j-\eta(j) \leq \delta_w -2r-1.
\end{equation}
Since $a' \neq 1$ for every multiplier  (because of the hypothesis $w \equiv 0$ modulo $3$), we conclude that
$i=\delta_w-2r-2$. 
Then, (\ref{diineq}) implies
\begin{equation} \label{djineq}
\delta_w-2r-2\leq \delta_w -j-\eta(j) \leq \delta_w -2r-1.
\end{equation}
As $j \equiv i \equiv \delta_w-2r-2 \equiv 1 \pmod{3}$, the inequality on the left of (\ref{diineq}) gives
$j \leq 2r+1$. Comparing with (\ref{jbbineq}) we conclude
$j= 2r+1$. This forces $b=r+1$ and $b' =r$. 
Then,
$$
\hat{\bfq}(\bfu)=\bfv-\bfu=(\delta_w\!-\!2r-2,2r+1,k)-(\delta_w,r+1,c)=(0,r,k-c)=(0,3h,m-2h)
$$
where the  second to last equality follows from the fact that  $\w(\hat{\bfq}(\bfu))=3m$.
But this contradicts the induction hypothesis, which implies that all multipliers $(0,3\lambda,m-2\lambda)$ with $\lambda < \lambda (\bfu^{(r)}) = h+1$ have already appeared as multipliers  $\hat{\bfq}(\bfu_1)=\bfq(\bfu_1)$ for monomials $\bfu_1 \in \rectw$ satisfying  $\lambda(\bfu_1)< \lambda (\bfu^{(r)})$ and a fortiori $\bfu_1 \prec \bfu^{(r)}$. We conclude that
$\hat{\varphi}(\bfu)$ must lie in $\mrectw$, and then, using Lemma \ref{lemmaunirectangle} as in the proof of Corollary  \ref{unirectangle}, that
$\hat{\varphi}= \varphi$ on the set of monomials $\preceq \bfu^{(r+1)}$.

\vspace{5pt}
\noindent {\bf Case 2}: $r \equiv 1 \pmod{3}$, say $r=3h+1$.

\noindent
In this case,  $\lambda (\bfu^{(r)}) =\lambda (\bfu) = \lambda (\bfu^{(r+1)})=h+1$,
and, by Lemma \ref{aneworder}(3), the first entry $a$ of $\bfu$ is either $\delta_w-2r-2$ or
$\delta_w -2r-1$. Let $\bfv=\hat{\varphi}(\bfu)=(i,j,k)$ and $\hat{\bfq}=\bfv-\bfu=(a',b',c')$.
Suppose by  contradiction that $\bfv \in \mtriw$. 
Then,
 by Proposition \ref{ijtriangle},  $\delta_w-2r-2 \leq a \leq i \leq \delta_w -j-\eta(j)$, and, as in the previous case, the induction hypothesis
 implies $j=b+b' \geq 2r+1$. 
 Suppose $a=\delta_w-2r-1$. 
The inequalities above imply that $i=\delta_w-2r-1$ and $j=2r+1$. 
Therefore,
 $a'=0$ and $b'=r \equiv 1 \pmod{3}$. But this is a contradiction, because
$a'=0$ implies that $\hat{\bfq}$  is a rectangular region multiplier, and then $b'=3 \lambda$ is congruent to $0$ modulo $3$.

Therefore, if $\bfv \in \mtriw$, then  $a=\delta_w-2r-2$. In this case, Proposition \ref{ijtriangle} implies that
$$\delta_w-2r-2 = a  \leq i \leq \delta_w -j-\eta(j)\leq \delta_w -2r-1-\eta(j).$$
Since we cannot have $i=a+1$, because no multiplier has $a'=1$, we conclude
that $i =\delta_w-2r-2$ and $a'=0$. 
Therefore,
$\hat{\bfq}(\bfu)=(0,3 \lambda, m -2 \lambda)$
for some $\lambda$, and $\lambda \geq \lambda (\bfu^{(r)})=h+1$ by the induction hypothesis. It follows that
$
j= b+ b'  \geq (r+1) + 3h+3=2r+3$
which gives the contradiction
$ \delta_w -2r-3=i \leq \delta_w - j \leq \delta_w -2r-3.$
Therefore, $\hat{\varphi}(\bfu)$ must lie in $\mrectw$, and, as in Case 1, 
we conclude that
$\hat{\varphi}= \varphi$ on the set of monomials $\preceq \bfu^{(r+1)}$.

\noindent {\bf Case 3}: $r \equiv 2 \pmod{3}$, say $r=3h+2$.
\\
\noindent
In this case $\lambda (\bfu^{(r)})= h+1 < \lambda (\bfu^{(r+1)})= h+2$.
By Lemma \ref{aneworder}(4) the first entry $a$ of $\bfu$ assumes one of the values $\delta_w-2r-1$,  $\delta_w-2r-2$,
 $\delta_w-2r-3$.

\noindent {\bf Case 3.1} Suppose $a= \delta_w -2r-1$ or $a= \delta_w -2r-3$, so that $ \lambda (\bfu)= h+1 $.
\\
\noindent
Let $\bfv=\hat{\varphi}(\bfu)=(i,j,k)$ and assume by contradiction that $\bfv \in \mtriw$.
Suppose first that $a=\delta_w -2r-1$.
Using the inequalities $\delta_w-2r-1 = a \leq i \leq \delta_w -j-\eta(j)$ and
 $j=b+b'\geq 2r+1$, we conclude that $i=a$, $j=2r+1$ and $\hat{\bfq}(\bfu)=\bfv-\bfu=(0,r,c')$, which is a contradiction
since $r \not \equiv 0 \pmod{3}$.

Suppose now  that $a=\delta_w -2r-3$. 
Arguing as in the previous cases, we
obtain the inequalities
$$\delta_w-2r-3  \leq i \leq \delta_w -2r-1-\eta(j).$$
Since $a' \neq 1$, either $i =\delta_w-2r-3$ or $i =\delta_w-2r-1$.

If $i=\delta_w-2r-3= a$, the inequality
$\delta_w-2r-3 \leq \delta_w -j-\eta(j)$ gives
$
j+\eta(j) \leq 2r+3,
$
while
$j \equiv i \equiv -2r \equiv 2 \pmod{3}$.
As  $j \geq 2r+1$ by induction, we must have $j=2r+1$, and, as above, this  leads to the contradiction
$\hat{\bfq}(\bfu)=\bfv-\bfu=(0,r,c')$.
Finally, it $i=\delta_w-2r-1= a+2$,  the inequality
$\delta_w-2r-1=i \leq \delta_w -j-\eta(j)$ implies $j=2r+1$. But then
$i \equiv 1$ while $j\equiv 2$ modulo $3$, in contradiction with Proposition \ref{ijtriangle}.

We conclude that $\hat{\varphi}(\bfu)$ must lie in $\mrectw$, and, as in the previous cases, that
$\hat{\varphi}= \varphi$ on the set of monomials $\preceq \bfu^{(r+1)}$.

\noindent {\bf Case 3.2} Suppose $a= \delta_w -2r-1$ , so that $ \lambda (\bfu)= h+2$.
\\
\noindent
Let $\bfv=\varphi(\bfu)=(i,j,k)$, and suppose by  contradiction that $\bfv \in \mtriw$. The inequalities $\delta_w-2r-2 = a  \leq i \leq \delta_w -j-\eta(j)$ and $j\geq 2r+1$
imply either $j=2r+1$ or $j=2r+2$. As $r \equiv 2 \pmod{3}$, the only
possibility is $j=2r+2$ and $i=\delta_w-2r-2$. It follows that $\bfu=\bfu^{(r+1)}$ and
$\hat{\bfq}(\bfu)=(0,3(h+1),m-2(h+1))$.
This cannot hold, because, by induction on the order, we may assume that the rectangular region multiplier $\bfq=(0,3(h+1),m-2(h+1))$ has
already appeared with its multiplicity, since $\hat{\bfq}(\bfu)= \bfq(\bfu)$ for those $\bfu$  in $\rectw$ that satisfy $\lambda(\bfu) =h+1$ and hence
precede  $ \bfu^{(r+1)}$ in the order $\prec$. 
Thus, $\bfv \in \rectw$, and then we conclude as in the previous cases that $\hat{\varphi}(\bfu)= \varphi(\bfu)$ for all $\bfu \preceq \bfu^{(r+1)}$.

\vspace{7pt}
To finish, we need to prove $S_{r,2}$ assuming that $S_{\ell,1}$ holds for $\ell \leq r$ and $S_{\ell,2}$ holds for $\ell \leq r-1$.
Pick $\bfv =(i,j,k) \in \scrM^{\triangle}_{w,r}$, then either $j =2r$  or $j=2r+1$. Let $\bfu=\hat{\varphi}^{-1}(\bfv)=(a,b,c)$.
We have already shown that the induction hypothesis implies that $S_{r+1,1}$ also holds.
If $\bfu \in \rectw$ by contradiction, we would have $\hat{\varphi}(\bfu) \neq \varphi(\bfu)$ hence $\bfu^{(r+1)} \prec \bfu$ by induction,
so that $b \geq r+2$ by Lemma \ref{urrectangular}.
Let $\hat{\bfq}=\bfv-\bfu=(a',b',c')$: by the induction hypothesis, we have $b' \geq r$.
It follows that  $j=b+b' \geq 2r+2$, contradicting $\bfv =(i,j,k) \in \scrM^{\triangle}_{w,r}$. 
We conclude that $\bfu \in \rtriw$. 
By the induction hypothesis, we have
$b \geq r$. As $j=b+b' \leq 2r+1$, it follows $b' \leq r+1$. 
By the induction hypothesis, if $\hat{\bfq}=\bfv-\bfu=(a',b',c')$ were a rectangular region multiplier, then $\hat{\bfq}=(0,3\lambda, m-2\lambda)$ with $\lambda \geq \lambda (\bfu^{(r+1)} )=\floor{\frac{r+4}{3}}$. 
This implies the contradiction
$b'=3\lambda \geq r+2$.
 Thus,  $\hat{\bfq}$ is a triangular region multiplier. 
 Now, the statement $S_{r,2}$ follows as in the proof
of Claim 1 in Theorem \ref{unitriangle}.
\end{proof}

\section{Examples of non-uniqueness} \label{examples}
The bijection $\varphi: \rscr'_{w} \ra  \mscr'_{w}$ constructed in the previous sections does not satisfy uniqueness when
$w\not \equiv 0 \pmod{3}$. 
In this section, we give the  first example of this phenomenon. 
More generally, we show that, in general, there exists {\em no bijection} $\hphi: \rscr'_{w} \ra  \mscr'_{w}$ satisfying both divisibility and uniqueness.

\bex
Let $m=3$ and $w=6m-10=8$. 
We have $w \equiv 2 \pmod{3}$, $\delta_w=1$ and $\tau_w=0$. 
The rectangular region $\rect_w$ consists of
the two monomials $\bfu_1=(0,4,0) \prec \bfu_2=(1,2,1)$,
while the triangular region $\rtri_w$ consists of the single monomial $\bfu_3=(0,1,2)$.
The bijection $\varphi$ sends these monomials respectively to $\bfv_1=(0,4,3)$, $\bfv_2=(1,5,2)$ and $\bfv_3=(1,2,4)$.
If $a$, $b$, $c$ and $d$ denote the coefficients of the monomials $z^3$, $y^3z$, $xyz^2$ and $xy^4$ in the polynomial $g$, then the matrix
that represents $g \cdot$ with respect to the basis $\{\bf{u}_i\}$ of $R'_w$ and the basis $\{\bf{v}_j\}$ of $M'_w$ is
$$
\mathbf{A}_w=
\begin{bmatrix}
  a & 0 & b \\
  c & b & d \\
  0 & a & c
\end{bmatrix}.
$$
The determinant of this matrix is $2abc-a^2d$. This shows that the bijection $\varphi$ does not satisfy uniqueness. In this case, there is
another bijection that satisfies divisibility and uniqueness, having $\{a,a,d\}$ as multiset of multipliers.
\eex 

\bex
Let $m=4$ and $w=6m-10=14$.
We have $w \equiv 2 \pmod{3}$, $\delta_w=1$ and $\tau_w=1$. 
The rectangular region $\rect_w$ consists of
the four monomials
$$\bfu_1=(0,7,0) \prec \bfu_2=(1,5,1)\prec \bfu_3=(0,4,2) \prec \bfu_4=(1,2,3),$$
mapped by  $\vphi$ to
$$\bfv_1=(0,7,4) \prec \bfv_2=(1,8,3)\prec \bfv_3=(0,4,6) \prec \bfv_4=(1,5,5),$$
while the triangular region $\rtri_w$ of the single monomial $\bfu_5=(0,1,4)$, mapped by
$\vphi$ to $(1,2,7)$.

Let $a$, $b$, $c$, $d, e$ denote the coefficients 
of $z^4$, $y^3z^2$, $y^6$, $xyz^3, xy^4z$
in the polynomial $g$. 
The matrix
that represents $g \cdot$ with respect to the bases $\{\bf{u}_i\}$ of $R'_w$ and  $\{\bf{v}_j\}$ of $M'_w$ is
$$
\mathbf{A}_w=
\begin{bmatrix}
  a & 0& b & 0 & c \\
  d & b& e & c & 0 \\
  0 & 0& a & 0 & b \\
  0 & a& d & b & e \\
  0 & 0& 0 & a & d
\end{bmatrix}.
$$
Since $\det(\mathbf{A}_w) = -2a^3be-2a^2cd+3a^2bd$, it follows that no bijection $\hphi: \rscr'_{w} \ra  \mscr'_{w}$ satisfies both divisibility and uniqueness.

\eex

\section{Non Cancellation Lemma} \label{spblocks}

In this section, we fix $m,n,\tau\in \mathbb{N}$ such that
$2n+2 \leq m$.
The  vectors
$$
\beone=(1,1,-1),\; \betwo=(0,3,-2)\index{e1 @ $\beone$}\index{e2 @ $\betwo$}
$$
form a $\Z-$basis of the lattice of integral vectors $(a',b',c')$ of weight $a'+2b'+3c'=0$.
For an integer $h$, define $\bfq_h=(0,0,m)+h \betwo=(0,3n,m-2h)$. The assumption $2n+2 \leq m$
guarantees that the monomials $\bfq_n$ and $\bfq_{n+1}$ have non-negative coordinates.

\bd\label{defabstractspblocks} \index{bb @ $\cab$} \index{cc @ $\cac$}
The special block $\cab$ is the set consisting of the
$2\tau+3$  monomials
\begin{equation*}
  \begin{cases}
    \bfu^0_{\ell} = \ell \betwo
        &
        \mbox{for } \quad  0 \leq \ell \leq \tau+1,
        \\% & \\
    \bfu^1_{\ell} = \beone+\ell \betwo
     & \mbox{for} \quad 0 \leq \ell \leq \tau.
  \end{cases}
\end{equation*}
The special block $\cac$ is the set consisting of the
$2\tau+3$  monomials
\begin{equation*}
  \begin{cases}
 \bfv^0_{\ell}  =\bfq_{n+\ell}=(0,0,m)+(\ell+n) \betwo
        &
        \mbox{for } \quad  1 \leq \ell \leq \tau+1,
        \\% & \\
\bfv^1_{\ell} =  \beone+\bfq_{n+\ell}= (0,0,m)+\beone+(\ell+n) \betwo
     & \mbox{for} \quad 0 \leq \ell \leq \tau+1.
  \end{cases}
\end{equation*}
\ed

Note that the superscript $i$ of $\bfu^i_{\ell}$ and $\bfv^i_{\ell}$ is equal to the first coordinate of the vector.
We let $\cab^i$ and $\cac^i$ denote the set of monomials with superscript equal to $i$ in $\cab$ and in $\cac$, respectively.
The sets $\cab^0$ and $\cac^1$ contain $\tau+2$ monomials, while
$\cab^1$ and $\cac^0$ contain  $\tau+1$ monomials.

For each $\ell_0$ in $\{0,\ldots,\tau+1\}$ we define a bijection $\hphi_{\ell_0}:\cab \ra \cac$ \index{aavphizl @ $\hphi_{\ell_0}$}
sending   $\bfu^0_{\ell_0}$  to  $\bfv^1_{\ell_0}$, and mapping $\cab^0 \setminus \{\bfu^0_{\ell_0}\}$ (resp. $\cab^1$) to
$\cac^0$ (resp. $\cac^1 \setminus \{\bfv^1_{\ell_0}\}$)
 so that order defined by the index $\ell$ is respected:
\begin{equation}\label{blockmonomialsnew}
  \begin{cases}
  \hphi_{\ell_0}(\bfu^0_{\ell}) =  \bfv^0_{\ell+1}
        &
        \mbox{for } \quad  0 \leq \ell < \ell_0,
                \\ %& \\
 \hphi_{\ell_0}(\bfu^0_{\ell_0}) =  \bfv^1_{\ell_0}
        &
                \\ %& \\
  \hphi_{\ell_0}(\bfu^0_{\ell}) =  \bfv^0_{\ell}
        &
        \mbox{for } \quad  \ell_0< \ell \leq \tau+1,
                \\ %& \\
 \hphi_{\ell_0}(\bfu^1_{\ell}) =\bfv^1_{\ell}
 & \mbox{for} \quad 0 \leq \ell < \ell_0,
        \\ %& \\
 \hphi_{\ell_0}(\bfu^1_{\ell}) = \bfv^1_{\ell+1}
 & \mbox{for} \quad \ell_0 \leq \ell \leq \tau.
  \end{cases}
\end{equation}

As in the previous sections, we call $\hphi_{\ell_0}(\bfu)-\bfu$ the
$\hphi_{\ell_0}$-multiplier of $\bfu$.

\bl \label{specialblockbijectionsnew}
The bijections $\hphi_{\ell_0}$  have the same multiset of multipliers, namely
$\bfq_n+\beone$ with multiplicity $1$, $\bfq_n$ and $\bfq_{n+1}$ both with multiplicity $\tau+1$.
In particular, they satisfy divisibility, that is,
the multipliers have non-negative entries.
Any two of these bijections have the same parity, that is,  $(\hphi_2)^{-1} \circ \hphi_1$ is an even permutation.
\el

\begin{proof}
The statement about the multipliers is clear.  To verify the claim about the parity  we note that  $\hphi_{\ell_0+1}^{-1} \circ \hphi_{\ell_0}$ is the $3$-cycle
$(\bfu^0_{\ell_0}, \bfu^1_{\ell_0}, \bfu^0_{\ell_0+1})$.
\end{proof}

The bijections $\hphi_{\ell_0}$ are the unique
bijections $\cab \ra \cac$ having multiset of multipliers as in the previous Lemma; in fact,
the following stronger result holds.

\bl \label{cabninjetive}
Suppose $\hphi: \cab \setminus \{ \bfu^0_{0}\} \ra \cac$ is injective and
\begin{enumerate}
  \item $\hphi$ has the divisibility property;
  \item none of the multipliers $\hat{\bfq} (\bfu)=\hphi(\bfu)-\bfu$ is of the form  $\bfq_h$ or $\bfq_h+\beone$ with $h<n$.
\end{enumerate}
Then, there are at least $(\tau+1)$ monomials
  $\bfu$ in $\cab  \setminus \{ \bfu^0_{0}\}$ such that $\hat{\bfq} (\bfu)=\bfq_n$.

If $\hphi$ can be extended to a bijection $\cab \ra \cac$ that still verifies (1) and (2)
and for which the multiplier
 $\bfq_n$ has multiplicity at most $(\tau+1)$, then the extension is one of the bijections $\hphi_{\ell_0}$ constructed above.
\el

\begin{proof}
Since $\hphi$ is injective, the image of $\hphi$ misses exactly one of the monomials in $\cac$.

Case 1. Suppose first that the monomial not in the image of $\hphi$ is $ \bfv^1_{\ell_1}$, $0 \leq \ell_1\leq \tau+1$. 
In this case, every monomial in $\cac^0$ is in the image of $\hphi$.
By the divisibility property, the inverse image of $\cac^0$ is contained in $\cab^0$,
thus, there is a permutation $\sigma$ of $\{1, \ldots, \tau+1\}$ such that
$$
\hphi (\bfu^0_{\ell}) =\bfv^0_{\sigma(\ell)} \quad \mbox{for all $1 \leq \ell \leq \tau+1$}.
$$
On the other hand, computing the corresponding multipliers
$$
\hphi(\bfu^0_{\ell})-\bfu^0_{\ell}= \bfv^0_{\sigma(\ell)}- \bfu^0_{\ell}=\bfq_{n+\sigma(\ell)-\ell}
$$
we see that $\sigma(\ell) \geq \ell$, using the assumption that no multiplier can be of the form $\bfq_h$ with $h <n$.
We conclude that $\sigma$ is the identity,
and thus the monomials $\bfu^0_{\ell}$, for $1 \leq \ell \leq  \tau+1$, all have multiplier $\bfq_n$, and this proves the first statement
of the lemma in Case 1.

Suppose furthermore that $\hphi$ can be extended to a bijection $\cab \ra \cac$, which we still denote by $\hphi$, that verifies (1) and (2) and for which the multiplier
 $\bfq_n$ has multiplicity at most $(\tau+1)$. 
 Then, $\hphi$ must map $\bfu^0_{0}$ to the monomial $\bfv^1_{\ell_1}$ missing from the image of the original map. Furthermore, there is a bijection $\tilde{\sigma}: \{0, \ldots, \tau\} \ra \{0, \ldots, \tau+1\} \setminus \{\ell_1\}$ such that
 $$\hphi (\bfu^1_{\ell}) =\bfv^1_{\tilde{\sigma}(\ell)} \quad \mbox{for \ $0 \leq \ell \leq \tau$}.$$
We must have $\tilde{\sigma}(\ell) >\ell$ for $0 \leq \ell \leq \tau$, because there cannot be more than $\tau+1$ multipliers equal to $\bfq_n$, and
there are no multipliers $\bfq_h$ with $h <n$.
Then, $0$ cannot be in the image of $\tilde{\sigma}$, so $\ell_1=0$. 
By descending induction on $\ell$, we see that
$\tilde{\sigma}(\ell)=\ell+1$ for $\ell=0,1,\ldots, \tau$.
We conclude  $\hphi=\hphi_{0}$.

Case 2.
Suppose now that the  monomial missing in the image is one of the monomials in $\cac^0$, say $\bfv^0_{p}$.
Then, the image of $\cab^0$ cannot be contained in $\cac^0$, that is, there are integers
$\ell_0$ and $\ell_1$, satisfying $1 \leq \ell_0 \leq \tau+1$ and $0 \leq \ell_1 \leq \tau+1$, such that
$\hphi ( \bfu^0_{\ell_0})= \bfv^1_{\ell_1}$.
We compute
$$
\hbfq (\bfu^0_{\ell_0}) = \bfv^1_{ \ell_1}-\bfu^0_{\ell_0}= \bfq_{n+\ell_1-\ell_0}+ \beone
$$
and we conclude $\ell_1 \geq \ell_0$ by the assumptions on the multipliers.

On the other hand,
by divisibility, $\hphi$ must map $\cab^1_n$ to $\cac^1$, so there is a bijection
$\tilde{\sigma}:\{0, \ldots, \tau\}\ra \{0, \ldots, \tau+1\} \setminus \{\ell_1\}$
such that $\hphi (\bfu^1_{\ell}) =\bfv^1_{\tilde{\sigma}(\ell)}$ for all $\ell \in \{0, \ldots, \tau\}$.
As in Case 1, we must have $\tilde{\sigma}(\ell) \geq \ell$ for all $\ell \in \{0, \ldots, \tau\}$, and this implies,
by ascending induction on $\ell$, that
$$\hphi (\bfu^1_{\ell}) =\bfv^1_{\ell} \quad \mbox{for all $0 \leq \ell < \ell_1$}.$$
Similarly, there is a bijection $\sigma: \{1,\ldots,\tau+1\}  \setminus \{\ell_0\} \ra \{1,\ldots,\tau+1\}  \setminus \{p\}  $
such that  $\hphi (\bfu^0_{\ell}) =\bfv^0_{\sigma(\ell)}$ for all $\ell$ in the domain of $\sigma$. 
Again, $\sigma(\ell) \geq \ell$ by the assumptions on the multipliers,
hence
 $$\hphi (\bfu^0_{\ell}) =\bfv^0_{\ell} \quad \mbox{for all $\ell_0 < \ell \leq \tau+1$}$$ because
we need to find an image $\bfv^0_{\sigma(\ell)}$ with $\sigma(\ell) \geq \ell$ for every $\bfu^0_{\ell}$ satisfying
$\ell_0 < \ell \leq \tau+1$. In particular, $p \leq \ell_0$.
This gives $\hat{\bfq} (\bfu)=\bfq_n$ for at least
$
(\tau+1)-\ell_0 + \ell_1 \geq \tau+1
$
monomials.

Finally,
suppose $\hphi$ can be extended to a bijection $\cab \ra \cac$, which we still denote by $\hphi$, that verifies (1) and (2) and for which the multiplier
 $\bfq_n$ has multiplicity at most $(\tau+1)$. 
 Then, we must have $\ell_1=\ell_0$, that is,  $\hphi ( \bfu^0_{\ell_0})= \bfv^1_{\ell_0}$.
Recall that  $\hphi$  maps $\bfu^1_{\ell}$ to $\bfv^1_{\tilde{\sigma}(\ell)}$
for some injective map
$$\tilde{\sigma}: \{0, \ldots,\tau\} \ra \{0, \ldots, \tau+1\}\setminus \{\ell_0\}$$
that satisfies $\tilde{\sigma}(\ell)=\ell$ if $0 \leq \ell < \ell_0$. The extra information we obtain from the assumption that the multiplier $\bfq_n$ has multiplicity at most $\tau+1$ is that
that $\tilde{\sigma}(\ell) > \ell$ for $\ell \geq \ell_0$.  
We deduce that
$$
\tilde{\sigma}(\tau)=\tau+1,\, \tilde{\sigma}(\tau-1)=\tau,\, \ldots,\, \tilde{\sigma}(\ell_0) = \ell_0+1,
$$
and
therefore $\hphi$ coincide with $\hphi_{\ell_0}$ on $\cab^1$.

In order
to prove that the two bijections coincide also on $\cab^0$,  observe that
 $\hphi$ must  send the  monomial
 $\bfu^0_{0}$ missing in the original domain to the monomial $\bfv^0_{p}$ missing in the image of the original map.
If $\ell \neq 0$ and $\ell \neq \ell_0$, the map $\hphi$  sends $\bfu^1_{\ell}$ to $\bfv^1_{n,\sigma(\ell)}$,
where
$$\sigma: \{1, \ldots,\tau+1 \}\setminus   \{\ell_0\} \ra \{1, \ldots, \tau+1\}\setminus \{p\}$$
is a bijection that satisfies $\sigma(\ell)=\ell$ if $\ell_0 < \ell \leq \tau+1$.
As above, from the assumption on the multiplicity of $\bfq_n$ we deduce that $\sigma(\ell) > \ell$ for $\ell < \ell_0$.
Hence, $p=1$ and $\sigma(\ell)=\ell+1$ for $1  \leq \ell <\ell_0$, that is, $\hphi$ coincides with $\hphi_{\ell_0}$.
\end{proof}

\section{The cases $w \equiv 1$ or $2$ modulo $3$} \label{wcongruounodue}

In this section, we assume the weight $w$ is congruent to either $1$ or  $2 \pmod{3}$, and satisfies $3m-2 \leq w \leq 6m-10$. 
In this range,
the rectangular regions $\rectw$ ad $\mrectw$ are nonempty.  
While 
the bijection $\vphi$ does not satisfy uniqueness,
we will  prove that $\vphi$ satisfies {\bf Non cancellation}, by showing that uniqueness only fails due to the presence of the special blocks analyzed in Section \ref{spblocks}. Since we assume that $\Bbbk$ has characteristic 0, we conclude that $g_w \cdot$ is an isomorphism.

We begin by introducing  filtrations on
$\rscr_w'$ and  $\mscr'_w$ that are needed to set up the inductive procedure of the proof.

\bd \label{defdeltan} \index{nmax @ $n_{max}$ }
Recall that $\delta_{w,n}=6m-9-w-6n$ (see page \pageref{orderdeltawn}).
We let $n_{max}$ denote the largest integer $n$ such that
$\delta_{w,n} \geq 0$, that is,
 $$n_{max}=\floor{\frac{6m-9-w}{6}}= m-1+\floor{\frac{-3-w}{6}}.$$
\ed

\br \label{nmaxw2}
We observe
\begin{enumerate}
 \item if $w \equiv 1$  (mod. 6), then $n_{max}=m\!-\!1\!-\!\frac{w+5}{6}$  and $\delta_{w,n_{max}}=2$;
 \item if $w \equiv 2$  (mod. 6), then $n_{max}=m\!-\!1\!-\!\frac{w+4}{6}$  and $\delta_{w,n_{max}}=1$;
 \item if $w \equiv 4$  (mod. 6), then $n_{max}=m\!-\!1\!-\!\frac{w+8}{6}$  and $\delta_{w,n_{max}}=5$;
 \item if $w \equiv 5$  (mod. 6), then $n_{max}=m\!-\!1\!-\!\frac{w+7}{6}$  and $\delta_{w,n_{max}}=4$.
\end{enumerate}
In any case, $2n_{max}+2 \leq m-1$, and the assumption $2n+2 \leq m$ of Section \ref{spblocks} is verified
by all $n \leq n_{max}$.
\er

\bd \label{hnkndefinition}
For  $0 \leq n \leq n_{max}$, define
\begin{equation}\label{hnformula}
\mathscr{H}_n=
\left\{
\bfu=(a,b,c) \in \rscr'_w:
\mbox{ $a \geq \delta_{w,n} -5$ or $b \leq 3n+2$}
\right\},
\end{equation}
\begin{equation}\label{knformula}
\mathscr{K}_n=
\left\{
\bfv=(i,j,k) \in \mscr'_w:
\mbox{ $i \geq \delta_{w,n} -5$ or $j \leq 6n+5$}
\right\}.
\end{equation}
\ed
It is clear that the sets $\scrhn$ (resp. $\scrkn$) increase with $n$, and
$\scrH_{n_{max}}=\rscr'_w$ (resp. $\scrK_{n_{max}}=\mscr'_w$).

Set $c_w=2(\tau_w+1)$.\index{cw @ $c_w$ the smallest value of $c(\bfu)$ for $\bfu \in \rtriw$} 
A monomial $\bfu=(a,b,c) \in \rscr'_w$ lies in the rectangular region
if and only if $c<c_w$. 
Monomials $\bfu$ with $c(\bfu)=c_w$ lie in the triangular region but are close to the rectangular region.

\bp \label{hkproperties} \
 Let $\vphi:  \rscr'_w \ra \mscr'_w$
be the bijection defined as $\varphi^\triangle$ on $\rtri_w$ and as $\varphi^\hrectangle$ on $\rect_w$.
Then, $\vphi (\scrhn)=\scrkn$ for every $n$. Furthermore,
\begin{enumerate}
  \item $\scrhn \cap \rectw  = \left\{\bfu=(a,b,c) \in \rectw: a\geq \delta_{w,n}-5 \right\},$
  \item $\scrkn \cap \mrectw  = \left\{\bfv=(i,j,k) \in \mrectw: i\geq \delta_{w,n}-5 \right\},$
  \item $\scrhn \cap \rtriw  = \left\{\bfu=(a,b,c) \in \rtriw: b \leq 3n+2 \right\} $ if $w \equiv 0$ or $w \equiv 1 \pmod{3}$, while
 $\scrhn \cap \rtriw  = \left\{\bfu=(a,b,c) \in \rtriw: b \leq 3n+2 \right\}  \cup \left\{ \bfu^{\triangle}_{2,n} \right\} $ if $w \equiv 2$  $\pmod{3}$ where  $\bfu^{\triangle}_{2,n}=(\delta_{w,n}-5,3n+3,c_w)$,
  \item $\scrkn \cap \mtriw  = \left\{\bfv=(i,j,k) \in \mtriw: j \leq 6n+5 \right\} $ if $w \equiv 0$ or $w \equiv 1 \pmod{3}$, while
 $\scrkn \cap \mtriw  = \left\{\bfv=(i,j,k) \in \mtriw: j \leq 6n+5 \right\} \cup \left\{\vphi \left( \bfu^{\triangle}_{2,n}\right)\right\}$ if $w \equiv 2$  $\pmod{3}$.
 \end{enumerate}
\ep

\begin{proof}
Suppose $\bfu=(a,b,c) \in \rectw$ and $a \leq \delta_{w,n}-6=\delta_w-6(n+1)$. Then,
$\bfu \succ \bfu^{(3n+2)}$ by Lemma \ref{aneworder}, hence $b \geq 3n+3$ by Lemma \ref{urrectangular}.
Therefore, on $\rectw$ the condition $b \leq 3n+2$ implies $a \geq \delta_{w,n}-5$, and
 $\scrhn \cap \rectw$  consists of those  $\bfu=(a,b,c) \in \rectw$ that satisfy $a\geq \delta_{w,n}-5$.

Suppose that $\bfv=(i,j,k) \in \mrectw$ and $i \leq \delta_{w,n}-6=\delta_w-6(n+1)$. The condition that $\bfv$ belongs to
$\mrectw$ is that $\ttm(\bfv) \leq \tau_w$, that is,
\begin{equation}\label{ttmcondition}
6m -6-w \leq i+j - \eta(-i-j)-\eta(w).
\end{equation}
If we add the condition $i \leq \delta_{w,n}-6$, we obtain $j \geq 6n+9$.
Thus, on $\mrectw$ the condition $j \leq 6n+5$ implies $i \geq \delta_{w,n}-5$, and
 $\scrkn \cap \mrectw$  consists of the  $\bfv=(i,j,k) \in \mrectw$ such that $i\geq \delta_{w,n}-5$.

 Suppose that $\bfu=(a,b,c) \in \rtriw$ and $b \geq 3n+3$. Since $\bfu \in \rtriw$,  formula (\ref{ttrconditionexplicit}) implies
 \begin{equation}\label{prot}
  a+w=2a+2b+3c \leq 6m-12-3\epsilon(c)+2 \eta(a+2b) -2b.
\end{equation}
 If we add the condition $b \geq 3n+3$, then we see that
 \begin{equation}\label{prot}
  a \leq 6m-12-w-3\epsilon(c)+2 \eta(a+2b) -6n-6=\delta_{w,n}-9-3\epsilon(c)+2\eta(a+2b).
\end{equation}
This implies $a \leq \delta_{w,n}-5$, and $a \leq \delta_{w,n}-6$ unless $b=3n+3$, $\epsilon(c)=0$ and $\eta(a)=2$.
Therefore, any $\bfv$ in $\scrhn$ which belongs to $\rtriw$ satisfies either $b \leq 3n+2$ or
$a=\delta_{w,n}-5$, $b=3n+3$,  $\epsilon(c)=0$ and $\eta(a)=2$. As $\delta_{w,n} \equiv -w \pmod{3}$, the latter conditions
imply $2 \equiv -w-5 \pmod{3}$, hence $w \equiv 2 \pmod{3}$ and
$\bfu=(\delta_{w,n}-5,3n+3, c_w)$.

Suppose that $\bfv=(i,j,k) \in \mtriw$ and $j \geq 6n+6$.
We have $i \leq \delta_w +\eta(w) -j - \eta(j-w)$
 by Lemma \ref{ijtriangle}.
If $j=6n+7$ or if $j=6n+6$ and $w \equiv 0$ or $1 \pmod{3}$, we deduce $i \leq \delta_{w,n}-6$.
If $j=6n+6$ and $w \equiv 2 \pmod{3}$, we deduce $i \leq \delta_{w,n}-5$.
Thus, $\scrkn \cap \mtriw  = \left\{\bfv=(i,j,k) \in \mtriw: j \leq 6n+5 \right\} $  when  $w \equiv 0$ or $1 \pmod{3}$,
while for  $w \equiv 2 \pmod{3}$ the intersection $\scrkn \cap \mtriw$ also contains  the unique monomial
of the triangular region for which $ i \geq \delta_{w,n}-5$ and $j \geq 6n+6$, which is
\begin{equation}
(\delta_{w,n}-5,6n+6,c_w+m-2n-2) = \vphi \left( \bfu^{\triangle}_{2,n}\right),
\end{equation}
see \ref{table trouble w2 v2}. 
We can easily check that $\vphi$ maps $\scrhn$ onto $\scrkn$:  by construction,  $\vphi$ maps a monomial
$\bfu=(a,b,c)$ in the rectangular region $\rectw$ to  a monomial
$\bfv=(i,j,k)$ in the rectangular region $\mrectw$ having $i=a$; and $\vphi$ maps a monomial
$\bfu=(a,b,c)$ in the triangular region $\rtriw$ to  a monomial
$\bfv=(i,j,k)$ in the triangular region $\mrectw$ with the same $r$-invariant, that is, with $\ds b=\floor{\frac{j}{2}}$.
 \end{proof}

\subsection{Special blocks}

\bd \label{defcabn}
Assume $w \equiv 1,2 \pmod{3}$ and let  $0 \leq n \leq n_{\max}$. 
The $n^{th}$  {\em special block}
\index{bn @ $\cab_n$}
$\cab_n$ of $\rscr'_w$ is the set of  monomials $\bfu=(a,b,c)$ in $\rscr'_w$ for which
$a$ is either $\delta_{w,n}$ or $\delta_{w,n}-1$, and
$b \geq 3n-1+\eta(w)$.
Similarly,   the $n^{th}$  {\em special block} \index{cn @ $\cac_n$}
$\cac_n$ of $\mscr'_w$ is the set of  monomials $\bfv=(i,j,k)$ in $\mscr'_w$ for which
$i$ is either $\delta_{w,n}$ or $\delta_{w,n}-1$ and
$j \geq 6n+2$.
\ed

\bl \label{specialblock}
Assume $w \equiv 1,2 \pmod{3}$ and let  $0 \leq n \leq n_{\max}$. 
 Let
\begin{equation}\label{vertexspecialblock} \index{uncorner @ $\bfu^0_{n,0}$ the corner monomial of the special block $\cab_n$}
  \bfu^0_{n,0}= (\delta_{w,n}\!-\!1, 3n-1+\eta(w),c_w ).
\end{equation}
Then, in the notation of Section \ref{spblocks} with $\tau=\tau_w$, we have
$\cab_n=\bfu^0_{n,0}+ \cab$ and $ \cac_n= \bfu^0_{n,0}+\cac$.
The restriction of $\vphi$ to $\cab_n$ is given by the formula
$\vphi(\bfu^0_{n,0}+\mathbf{z})=
\bfu^0_{n,0}+\hphi_{0} (\mathbf{z}) \quad (\mathbf{z} \in \cab).$
In particular, $\cac_n= \vphi (\cab_n)$.
\el

\begin{proof}
Explicitly, the statement is that the special block $\cab_n$ consists of the monomials
\begin{equation*}
  \begin{cases}
    \bfu^0_{n,\ell} = \left(\delta_{w,n}\!-\!1, 3(n+\ell)-1+\eta(w),c_w \!-\! 2 \ell\right)=
    \bfu^0_{n,0}+\ell \betwo
        &
        \mbox{for } \quad  0 \leq \ell \leq \tau_w+1,
        \\% & \\
    \bfu^1_{n,\ell} = \left(\delta_{w,n}, 3(n+\ell)+\eta(w),c_w \!-\! 2 \ell\!-\!1\right)=\bfu^0_{n,0}+\ell \betwo+\beone
     & \mbox{for} \quad 0 \leq \ell \leq \tau_w;,
  \end{cases}
\end{equation*}
that the special bloc $\cac_n$ consists of the monomials
\begin{equation*}
  \begin{cases}
 \bfv^0_{n,\ell}  = \left(\delta_{w,n}\!-\!1, 6n+3\ell-1+\eta(w),c_w +m\!-\! 2(n+ \ell)\right)=\bfu^0_{n,0}+\bfq_{n+\ell}
        &
        \mbox{for } \quad  1 \leq \ell \leq \tau_w+1,
        \\% & \\
 \bfv^1_{n,\ell} = \left(\delta_{w,n}, 6n+3\ell+\eta(w),c_w +m\!-\! 2(n+ \ell)-1\right)=
 \bfu^0_{n,0}+\bfq_{n+\ell}+\beone
     & \mbox{for} \quad 0 \leq \ell \leq \tau_w+1,
  \end{cases}
\end{equation*}
and finally that the bijection $\vphi$ on $\cab_n$ is given by
\begin{equation}\label{blockmonomialsphi}
  \begin{cases}
 \varphi(\bfu^0_{n,0}) =  \bfv^1_{n,0},
        &
                \\% & \\
 \varphi(\bfu^0_{n,\ell}) =  \bfv^0_{n,\ell}
        &
        \mbox{for } \quad  1 \leq \ell \leq \tau_w+1,
        \\% & \\
 \varphi(\bfu^1_{n,\ell}) = \bfv^1_{n,\ell+1}
 & \mbox{for} \quad 0 \leq \ell \leq \tau_w.
  \end{cases}
\end{equation}
The $\vphi$-multipliers $\bfq(\bfu)=\vphi(\bfu)-\bfu$ are therefore
\begin{equation}\label{blockmonomialsmult}
  \begin{cases}
    \bfq(\bfu^0_{n,0}) = \bfq_n+\beone=(1,3n+1,m-2n-1),
        &
                \\
    \bfq(\bfu^0_{n,\ell}) = \bfq_n=(0,3n,m-2n)
        &
        \mbox{for } \quad  1 \leq \ell \leq \tau_w+1,
        \\% & \\
    \bfq(\bfu^1_{n,\ell}) = \bfq_{n+1}=(0,3n+3,m-2n-2)
     & \mbox{for} \quad 0 \leq \ell \leq \tau_w.
  \end{cases}
\end{equation}
To prove the statement, observe that monomials in $\rscr'_w$ with given first coordinate $a$ differ from one another by a multiple of $\betwo=(0,3,-2)$.
One monomial of weight $w$ and first coordinate $\delta_{w,n}\!-\!1$ is
$\bfu^0_{n,0}$, and one monomial of weight $w$ and first coordinate $\delta_{w,n}$
is   $\bfu^1_{n,0}= \bfu^0_{n,0}+ \beone$. 
Adding $\ell \betwo$, with $\ell$ negative, gives a monomial with
$b \leq 3n-4+\eta(w)$ in the former case and $b \leq 3n-3+\eta(w)$ in the latter. The upper bound on $\ell$ is what is needed to keep the third coordinate
of $\bfu$ non-negative.

Similarly,  in $\mscr'_w$ one monomial of $\mscr$-weight  $w$ and first coordinate $\delta_{w,n}\!-\!1$  is
$\bfv^0_{n,0}=(\delta_{w,n}\!-\!1, 6n+\eta(w)-1,c_w+m-2n)$, and one monomial of weight $w$ and first coordinate $\delta_{w,n}$
is   $\bfv^1_{n,0}= \bfv^0_{n,0}+ \beone$. 
To be in $\cac_n$, a monomial $\bfv=(i,j,k)$
must have $j \geq 6n+\eta(w)$, thus, for $\bfv^0_{n,0}+\ell \betwo$
(resp.   $\bfv^1_{n,0}+\ell \betwo$) to be in $\cac_n$ we need $\ell \geq 1$ (resp. $\ell \geq 0$).
The upper bound on $\ell$ is what is needed for $i+j$ to be less or equal than $3(m-1)$.

As to the computation of the multipliers, note that all monomials different from
$\bfu^0_{n,0}$ in the special block $\cab_n$ are in the rectangular region because their third coordinate is $< c_w$.  For the monomials $\bfu$ in the rectangular region the multiplier is $\bfq_h$ where $h=\lambda(\bfu)=\lambda(\delta_w-a(\bfu))$,
and our claim follows from $\lambda(6n+1)=n$ and $\lambda(6n)=n+1$. The fact that the $\vphi$-multiplier of $\bfu^0_{n,0}$ is $\bfq_n+\beone$
will be shown below when we compute the $\vphi$-multipliers of monomials in the triangular region with $c=c_w$.
\end{proof}

\bd \label{cornermonomial}
Fix an integer $n$ satisfying $0 \leq n \leq n_{\max}$.
We say that the monomial
$\bfu^0_{n,0}= (\delta_{w,n}\!-\!1, 3n-1+\eta(w),c_w )$ of formula (\ref{vertexspecialblock}) is
the {\em corner} monomial of the special block $\cab_n$.
It is the unique monomial $\bfu$ in $\cab_n$ that lies in $\rtriw$.

For $0 \leq \ell \leq \tau_w+1$ we define a bijection $\hphi_{n,\ell}:\cab_n \ra \cac_n$
by the formula \index{aavphizln @ $\hphi_{n,\ell}$}
$$\hphi_{n,\ell}(\bfu^0_{n,0}+\mathbf{z})=
\bfu^0_{n,0}+\hphi_{\ell_0} (\mathbf{z}) \quad (\mathbf{z} \in \cab),
$$
where $\hphi_{\ell_0}: \cab \ra \cac$ is as in Section \ref{spblocks} with $\tau=\tau_w$. 
\ed 

\br
The restriction of $\vphi$ to $\cab_n$ is $\hphi_{n,0}$, by Lemma \ref{specialblock}.
By Lemmas \ref{specialblockbijectionsnew} and \ref{cabninjetive},
among the bijections $\cab_n \ra \cac_n$, the $\hphi_{n,\ell}$ are precisely those
having the same multiset of multipliers as $\hphi_{n,0}$, and they have the same parity as $\hphi_{n,0}$.
In other words, $\hphi_{n,0}: \cab_n \ra \cac_n$ satisfies {\bf Non cancellation}.
\er

We can now state the result that concludes the proof of Theorem \ref{mainthm}.
\bt \label{uniquew2}
Suppose $w \equiv 1, 2$ modulo $3$ and  $3m-1 \leq w \leq 6m-10$. Let $\hphi:\rscr'_w \ra \mscr'_w$ be a bijection having the same multipliers multiset as $\vphi$. Then
\begin{enumerate}
  \item $\hphi$ and $\vphi$ coincide on monomials that do no belong to any special block $\cab_n$;
  \item   for each $0 \leq n \leq n_{\max}$, the restriction of $\hphi$ to the special block $\cab_n$ is one of the $(\tau_w+2)$ bijections $\hphi_{n,\ell_0}: \cab_n \ra \cac_n$ of Definition \ref{cornermonomial}.
\end{enumerate}
In particular, $\hphi$ and $\vphi$ have the same parity, and  the bijection $\vphi$ satisfies {\bf Non cancellation}.
\et

\subsection{Proof of Theorem \ref{uniquew2} in case $w \equiv 2 \pmod{3}$} \label{wcongruodue}
We need an  auxiliary result, namely, we
need to list the monomials of the triangular region whose $\vphi$-multiplier $\bfq=(a',b',c')$ has
$a'=0,1$.

\bp \label{table trouble w2 v2}
Suppose $w \equiv 2 \pmod{3}$.
We list in Table \ref{tavolabase} the monomials $\bfu$ of the triangular region $\rtriw$ whose $\vphi$-multiplier has $a'$ equal to either $0$ or $1$,
together with their $\vphi$-multipliers $\bfq(\bfu)$ and their images $\vphi(\bfu)$ in $\mtriw$. It follows
there is no $\bfu$ in the triangular region $\rtriw$ whose $\vphi$-multiplier is $\bfq_0=(0,0,m)$.

\noindent
\hoffset -0.6truecm
\begin{table}[h]
\caption{$\delta_{w,n}=\delta_w-6n$, $c_w=2(\tau_w+1)$, $\bfq_n=(0,3n,m-2n)$, $\beone=(1,1,-1)$  }
\label{tavolabase}
\footnotesize{
\begin{tabular}{|l|l|l|l|}
\hline
&&&
\\
  $\bfu$ & $\bfq(\bfu)$ & $\varphi(\bfu)$ & range of $n$
\\
&&&
\\
\hline
  \multirow{2}{*}{$\bfu^{\triangle}_{1,n}=(\delta_{w,n}-3, 3n+2,c_w) $}
 & \multirow{2}{*}{$\bfq_{n+1}$}
 &  \multirow{2}{*}{$(\delta_{w,n}-3, 6n+5,c_w+m-2n-2)$}
 &  $0 \leq n \leq n_{max}-1$ if $w \equiv 2 \pmod{6}$
  \\
  &&&
 $0 \leq n \leq n_{max}$ if $w \equiv 5 \pmod{6}$
 \\
\hline
$\bfu^{\triangle}_{2,n}=(\delta_{w,n}-5, 3n+3,c_w) $
& $\bfq_{n+1}$
&$(\delta_{w,n}-5, 6n+6,c_w+m-2n-2)$
&  $0 \leq n \leq n_{max}-1$
\\
\hline
 $\bfu^{\triangle}_{3,n}=(\delta_{w,n}-1, 3n+1,c_w) $
 & $\bfq_{n}+\bfe_1$
 &$(\delta_{w,n}, 6n+2,c_w+m-2n-1)$
 &
 $0 \leq n \leq n_{max}$
 \\
\hline
  \multirow{2}{*}{$\bfu^{\triangle}_{4,n}=(\delta_{w,n}-2, 3n,c_w+1) $}
 & \multirow{2}{*}{$\bfq_{n}+\bfe_1$}
 &  \multirow{2}{*}{$(\delta_{w,n}-1, 6n+1,c_w+m-2n)$}
 &  $0 \leq n \leq n_{max}-1$ if $w \equiv 2 \pmod{6}$
 \\
  &&&
 $0 \leq n \leq n_{max}$ if $w \equiv 5 \pmod{6}$
\\
\hline
\end{tabular}
}
\end{table}
\ep

\begin{proof}
We first explain the notation used in the table. The multipliers are $\bfq_n=(0,3n,m-2n)$ and
$\bfq_n+\beone=(1,3n+1,m-2n-1)$.
Recall that $c_w=2(\tau_w+1)$ is the smallest possible value of the third coordinate $c(\bfu)$ for monomials $\bfu$ in the triangular region $\rtriw$.

The statement follows easily from Lemma \ref{aprime}. Write $\bfu=(a,b,c)$ and $\bfq(\bfu)=(a',b',c')$. If $a'=0$, then by Lemma \ref{aprime} we have  that
$t_\rscr(\bfu)=\tau_w+1$, $a$ is congruent to either $1$ or $2$ modulo $3$, and $c \equiv a+w \equiv 0$ modulo $2$. 
Thus,
$c=2(\tau_w+1)=c_w$ and $\delta_w-a$ is congruent to either $3$ or $5$ modulo $6$, therefore either
$\bfu=\bfu^{\triangle}_{1,n}=(\delta_{w,n}-3, 3n+2,c_w) $ or $\bfu=\bfu^{\triangle}_{2,n}=(\delta_{w,n}-5, 3n+3,c_w)$ for some $n$.

If $a'=1$, then by Lemma \ref{aprime} we have $t_\rscr(\bfu)=\tau_w+1$ and either $a$ is congruent to $0$ modulo $3$ and $c \equiv a+w \equiv 0$ modulo $2$,
or $a$ is congruent to $2$ modulo $3$ and $c \equiv a+w \equiv 1$ modulo $2$. It follows that
either $c=2(\tau_w+1)=c_w$ and $\delta_w-a$ is congruent to $1$ modulo $6$, or
or $c=2(\tau_w+1)+1=c_w+1$ and $\delta_w-a$ is congruent to $2$ modulo $6$.
Therefore, either
$\bfu=\bfu^{\triangle}_{3,n}=(\delta_{w,n}-1, 3n+1,c_w) $ or $\bfu^{\triangle}_{4,n}=(\delta_{w,n}-2, 3n,c_w+1)$ for some $n$.
\end{proof}

\br \label{2cornermonomial}
Observe that the monomial
$  \bfu^\triangle_{3,n}= (\delta_{w,n}\!-\!1, 3n+1,c_w )$
is the corner monomial $\bfu^{0}_{n,0}$ of the special block $\cab_n$. It is also one of the (at most two) monomials whose $\vphi$-multiplier is
$\bfq_n+\beone$. The monomial $ \bfu^\triangle_{2,n}$ has already appeared in Proposition \ref{hkproperties}.
\er

\begin{proof}[Proof of Theorem \ref{uniquew2} in case $w\equiv 2$]
The proof proceeds by induction on $n$ examining the behavior of $\hphi$ on the {\em boundary}
$\partial \hscr_n= \hscr_n \setminus \hscr_{n-1}$.
 By Proposition \ref{hkproperties} , for $0 \leq n \leq n_{max}$
the boundary $\partial \hscr_n$ is the disjoint union of
\begin{enumerate}
\item the special block $\cab_n$,
\item the set of monomials $\bfu=(a,b,c) \in \rectw$ satisfying $\delta_{w,n}-2 \geq a \geq  \delta_{w,n}-5$,
\item the set of monomials $\bfu=(a,b,c) \in \rtriw$ satisfying $3n \leq b \leq 3n+2$, with the exception
of $\bfu^{\triangle}_{2,n-1}$, that has $b=3n$, but is in $\hscr_{n-1}$, and of the corner monomial
$\bfu^{\triangle}_{3,n}$ that has $b=3n+1$, but is in the special block $\cab_n$, and
\item the monomial $\bfu^{\triangle}_{2,n}=(\delta_{w,n}-5,3n+3,c_w)$ (for $n < n_{max}$).
\end{enumerate}

Before we can write the induction statement, we need to describe in detail the multipliers multiset common to $\vphi$ and $\hphi$.
Recall  that multipliers are monomials $\bfq=(a',b',c')$ satisfying $a'+2b'+3c'=3m$. If
$a'=0$, then  $\bfq=\bfq_n=(0,3n,m-2n)$ for some integer $n$, because the weight is $3m$ and therefore $b'$ is a multiple of $3$.
The monomials $\bfq_n$ occur as $\vphi$-multipliers, that is, as $\bfq (\bfu)=\vphi(\bfu)-\bfu$ for some monomial $\bfu \in \rscr'_w$,
with the following multiplicities:
\begin{enumerate}
  \item[a)] $\bfq_0=(0,0,m)$ is the $\vphi$-multiplier $\bfq (\bfu)=\vphi(\bfu)-\bfu$ of $(t_w+1)$ monomials $\bfu$, namely of
  those monomials $\bfu$  that belong to the rectangular region and have $a=\delta_w-1$ (these monomial belong to the special block $\cab_0$);
 $\bfq_0=(0,0,m)$ is not the $\vphi$-multiplier of any monomial of the triangular region by Proposition \ref{table trouble w2 v2}.
  \item[b)]
  $\bfq_n=(0,3n,m-2n)$ for $1 \leq n \leq n_{max}$ is the $\vphi$-multiplier of  $6(t_w+1)$ monomials of the rectangular region
 (those $\bfu =(a,b,c) \in \rectw$ for which $a$ assumes one of the values $\delta_{w,n\!-\!1}, \delta_{w,n\!-\!1}-2,\delta_{w,n\!-\!1}-3, \delta_{w,n\!-\!1}-4, \delta_{w,n\!-\!1}-5, \delta_{w,n}-1$ so that $\lambda(\bfu)=n)$);
  it also appears twice as the multiplier of $\bfu^\triangle_{1,n-1}$ and $\bfu^\triangle_{2,n-1}$. Note that exactly $(t_w+1)$ of these monomials,
  namely those in the rectangular region having $a=\delta_{w,n}-1$, belong to the boundary $\partial \hscr_n$, the rest belong to $ \hscr_{n-1}$
  (see table \ref{tavolabase}).
  \item[c)] (not needed for the proof)
  When $w \equiv 2 \pmod{6}$, $\bfq_{n_{max}+1}$ is the
 $\vphi$-multiplier of   exactly $(\tau_w+1)$ monomials, those in the rectangular region having $a=\delta_{n_{max}}=1$.
 When $w \equiv 5 \pmod{6}$, it is the
 $\vphi$-multiplier of   exactly   $4(\tau_w+1)$ monomials, those
  in the rectangular region having $a\in \{4,2,1,0\}$; it is also the $\vphi$-multiplier of one monomial of the triangular region, namely of
  $\bfu^{\triangle}_{1,n_{max}}=( 1,3 n_{max}+2, c_w)$.
\end{enumerate}

On the other end, the multipliers $\bfq=(a',b',c')$ with $a' >0$ arise only for the triangular region, and, we count them collectively according to the value of $b'$. As in Theorem \ref{unitriangle} let
$$
\scrR^{\triangle}_{w,r}= \big\{\bfu \in \rtriw: b(\mf{u})=r\big\}, \quad
\scrM^{\triangle}_{w,r}= \big\{\bfv \in \mtriw:  \floor{j(\mf{v})/2}=r\big\},
$$
$$
\npr=\#\big\{\bfv \in \scrM^{\triangle}_{w,r}: j(\mf{v})=2r\big\}, \quad
\nmr=\#\big\{\bfv \in \scrM^{\triangle}_{w,r}: j(\mf{v})=2r+1\big\},
$$
so that $\# \scrR^{\triangle}_{w,r}= \# \scrM^{\triangle}_{w,r}= \npr+\nmr$.
Furthermore, given an integer $h$, the number of $\bfu \in \rtriw$ whose
$\vphi$-multiplier  $\bfq=(a',b',c')=\vphi(\bfu)-\bfu$ has $b'=h$ is equal to $\nu^{+}_{h}+\nu^{-}_{h\!-\!1}$;
note that when $h=3n$ this number counts also the occurrences of $\bfq_n$ as the $\vphi$-multiplier of
 the monomials $\bfu^\triangle_{1,n-1}$ and $\bfu^\triangle_{2,n-1}$.

By induction we will show that for each $0 \leq n \leq n_{max}$ the following statements hold:
\begin{enumerate}[label=(\roman*)]
\item
the restriction of $\hphi$ to the special block $\cab_n$ is, up to translation as in the statement, one of the bijections $\hphi_{\ell_0}$
 of Section \ref{spblocks};

\item
$\hphi(\bfu)=\vphi(\bfu)$ if $\bfu=(a,b,c) \in \rectw$ satisfies $\delta_{w,n}-2 \geq a \geq \delta_{w,n}-5$;

\item
$\hphi$ sends monomials $\bfu=(a,b,c) \in \rtriw$ satisfying $3n \leq b \leq 3n+2$, with the possible exception of the corner monomial
$\bfu^{\triangle}_{3,n}$, to monomials $\hat{\bfv}=(i,j,k) \in \mtriw$ having the same $r$-invariant (that is, $\floor{\frac{j}{2}}=b$);

\item
$\hphi\left(\bfu^{\triangle}_{2,n}\right) =\vphi\left(\bfu^{\triangle}_{2,n}\right)$ (for $n < n_{max}$); note that
from statements (i)-(iv) it follows that $\hphi$ maps the boundary $\partial \hscr_n$ of $\hscr_n$ bijectively onto the
boundary $\partial \kscr_n=\kscr_n \setminus \kscr_{n-1}$ of $\kscr_n$.

\item
if $\bfu \in \partial \hscr_n$ and $\hat{\bfq}(\bfu)=\hphi(\bfu)-\bfu=(a',b',c')$ is the corresponding
$\hphi$-multiplier, then $3n \leq b' \leq 3(n+1)$ and
\begin{enumerate}[label=(\Alph*)]
\item
$\bfq_n=\hat{\bfq}(\bfu)$ for exactly $(\tau_w+1)$ monomials $\bfu$ in $\cab_n$;
\item if $n>0$ (resp. if $n=0$)
there are  $\nu^{+}_{3n}-1$ (resp. $\nu^{+}_{0}$) monomials $\bfu  \in \partial \hscr_n$ whose $\hphi$ multiplier has $b'=3n$ and $a'>0$;
\item
the $\hphi$-multipliers having $b'=3n+1$ all arise from monomials $\bfu  \in \partial \hscr_n$;
\item
the $\hphi$-multipliers having $b'=3n+2$ all arise from monomials $\bfu  \in \partial \hscr_n$; 
\item
 if $n < n_{max}$, the $\hphi$-multipliers having $b'=3(n+1)$ arise from monomials $\bfu  \in \partial \hscr_n$, except
 for  $(\tau_w+1) + \nu^{+}_{3(n+1)}-1$ many, of which at most  $(\tau_w+1)$ are equal to $\bfq_{n+1}$.
\end{enumerate}
\end{enumerate}
The Theorem follows from these statements:  this is clear, except for the equality $\hphi(\bfu)=\vphi(\bfu)$ for all monomials $\bfu$ in the triangular region that do not belong to any special block (that is, are different from the corner monomials). To prove this
last point, note that the restrictions of $\hphi$ and $\vphi$ to the special block $\cab_n$ have the same multiset of multipliers by induction and
Lemma \ref{specialblockbijectionsnew}, while their restrictions to the rectangular region minus the special blocks  coincide by induction.
It follows
that the  restrictions of $\hphi$ and $\vphi$ to the triangular region minus the union of the special blocks also have the same multiset of multipliers,
hence they coincide by Theorem \ref{unitriangle}.

We break the inductive proof of statements $(i)-(v)$ into seven steps. These steps include the initial case $n=0$. For $n>0$ we proceed by induction,  in particular, we assume that $\hphi$ maps $\hscr_{n-1}$ to $\kscr_{n-1}$,
and that all $\hphi$-multipliers $\hbfq=(a',b',c')$ of monomials not in $\hscr_{n-1}$ have $b' \geq 3n$.
The whole argument is summarized in table \ref{outlinepfthm2}.

\noindent{\bf Step 1.}
The bijection $\hphi$ maps the special block $\cab_n$ to the special block $\cac_n$, and  the restriction of $\hphi$ to
$\cab_n$ is, up to translation, one of the $(\tau_w+2)$ bijections of Lemma \ref{specialblockbijectionsnew}, hence it has the same multiset of multipliers
as the restriction of $\vphi$ to $\cab_n$. 

\noindent{\bf Proof of Step 1.}
Suppose $\bfu=(a,b,c) \in \cab_n$ and let $\hat{\bfv}=(i,j,k)=\hphi(\bfu)$. We claim that $\hat{\bfv}\in\cac_n$.
As $\bfu \notin \hscr_{n-1}$, by induction
$\hat{\bfv} \notin \kscr_{n-1}$, and in particular $i \leq \delta_{w,n}$. As $\bfu \in \cab_n$, we
have $\delta_{w,n} \geq a \geq \delta_{w,n}-1$ and $b \geq 3n+1$. By divisibility, $i \geq a$, so that $\delta_{w,n} \geq i \geq \delta_{w,n}-1$,
and to check that $\hat{\bfv}$ belongs to $\cac_n$ it remains to show $j \geq 6n+2$.
By the induction hypothesis on the multipliers,
if $\hat{q}=(a',b',c')=\hat{\bfv}-\bfu$, then $b' \geq 3n$.
Therefore $j =b+b' \geq 6n+1$. Thus we have to exclude the case $j=6n+1$.

Suppose first $\hat{\bfv}$ is in the rectangular region. The inverse image $\bfu_{0}=(a_0,b_0,c_0)=\vphi^{-1}(\hat{\bfv})$
under $\vphi$ of $\hat{\bfv}$
is then in the rectangular region, so that  $\hat{\bfv}= \vphi(\bfu_0)=\bfu_0+\bfq_{\lambda(\bfu_0)}$, that is,
$i=a_0$ and $j=b_0+3 \lambda (\bfu_0)$. As $a_0$ is equal to either $\delta_{w,n}$ or $\delta_{w,n}-1$ and
 $\bfu_{0}$ belongs to the rectangular region,  $b_0 \geq 3n+2$ by Lemmas \ref{urrectangular} and \ref{aneworder}.
The value of $\lambda(\bfu_0)$ is
$n$ if $a_0=\delta_{w,n}-1$, and is $n+1$ if $a_0=\delta_{w,n}$. Thus,
$j =b_0+3 \lambda \geq 3n+2+3n=6n+2$, which means that $\hat{\bfv}$ belongs to $\cac_n$ as claimed.

Suppose next that $\hat{\bfv}$ is in the triangular region and $j=6n+1$. 
Recall that $i \in \{\delta_{w,n},\delta_{w,n}-1\}$.
Since $i\equiv w-2j \equiv 0$ and $\delta_{w,n} \equiv -w \equiv 1 \pmod{3}$, we must have $i=\delta_{w,n}-1$. Then,
\begin{equation}\label{vfn}
\hat{\bfv}=(\delta_{w,n}-1,6n+1,c_w+m-2n)= \bfu^{0}_{n,0}+\bfq_n.
\end{equation}
(Note that, under the original bijection $\vphi$,  this monomial is the image of $\bfu^\triangle_{4,n}$ with multiplier
$\bfq_n+\beone$ by Proposition \ref{table trouble w2 v2}, and $\bfu^\triangle_{4,n}$ {\em does not} belong to $\cab_n$).

Summarizing, we have shown that the image under $\hphi$ of the special block $\cab_n$
is contained in $\cac_n \cup \{ \bfu^{0}_{n,0}+\bfq_n\}$. 
Suppose by contradiction that
$\hphi(\cab_n) \neq  \cac_n$. 
Then, there is $\bfu \in \cab_n$ such that $\hphi(\bfu)=  \bfu^{0}_{n,0}+\bfq_n$.
Since $j( \bfu^{0}_{n,0}+\bfq_n)= 6n+1$, the monomial $\bfu$ must have $b \leq 3n+1$ because
of the induction hypothesis that all $\hphi$-multipliers of monomials not in $\hscr_{n-1}$  have $b' \geq 3n$.
The only monomial in $\cab_n$ that satisfies $b \leq 3n+1$ is the corner monomial $\bfu^{0}_{n,0}$. We conclude that the $\hphi$-multiplier of the corner monomial of $\cab_n$ is $\bfq_n$. On the other hand, the bijection $\hphi$ must send the other monomials of $\cab_n$
injectively in $\cac_n$, and by Lemma \ref{cabninjetive} the multiplier $\bfq_n$ occurs
as $\hbfq(\bfu)$ for at least  $(\tau_w+1)$ monomials in $\cab_n$, the corner monomial excluded. Adding the corner monomial, we count at least  $(\tau_w+2)$ occurrences of $\bfq_n$ as a multiplier, contradicting the induction hypothesis (vE) that
at most $(\tau_w+1)$ monomials not in $ \hscr_{n-1}$ have  $\hphi$-multiplier equal to  $\bfq_n$ (for $n=0$, the multiplicity of the multiplier
$\bfq_0$ is exactly $(\tau_w+1)$, as  noted at the beginning of the proof).

We conclude that $\hphi (\cab_n)=\cac_n$ and, using Lemma \ref{cabninjetive} again, that the restriction of
$\hphi$ to $\cab_n$  is, up to translation, one of the bijections $\hphi_{\ell_0}$ of Lemma \ref{specialblockbijectionsnew}.
In particular, $\bfq_n$ is the $\hphi$-multiplier of exactly $(\tau_w+1)$ monomials $\bfu$ in $\cab_n$, and by
the induction hypothesis this accounts for all the remaining
occurrences of  $\bfq_n$ as a $\hphi$-multiplier: if $\hat{\bfq} (\bfu)=\bfq_n$, then $\bfu \in \hscr_{n-1} \cup \cab_n$.

\noindent {\bf Step 2.}
The bijection $\hphi$  maps $\scrR^{\triangle}_{w,3n}$ onto $\scrM^{\triangle}_{w,3n}$.
The set $\scrR^{\triangle}_{w,3n}$  contains one monomial
from  $\hscr_{n-1}$, namely $\bfu^\triangle_{2,n-1}=(\delta_{w,n}+1,3n,c_w)$,
and this monomial, by induction,  has $\hphi$-multiplier $\bfq_{n}$.
If $n>0$ (resp. if $n=0$), there are $\nu^+_{3n}-1$ (resp. $\nu^+_{3n}$)  monomials in $\scrR^{\triangle}_{w,3n} \setminus \left\{\bfu^\triangle_{2,n-1}\right\}$ whose $\hphi$-multiplier
has $b'=3n$, while $\nu^-_{3n}$ monomials in $\scrR^{\triangle}_{w,3n}$ have $\hphi$-multiplier with $b'=3n+1$. This accounts for all
multipliers left with $b'=3n$, that is, all monomials not in $\hscr_{n-1} \cup \cab_n \cup \scrR^{\triangle}_{w,3n}$ have $\hphi$-multiplier $\hat{\bfq}=(a',b',c')$ with $b' \geq 3n+1$.

\noindent{\bf Proof of Step 2.}
The argument we use to prove Step 2 is similar to that in the proof of Statement $S_r$ of Theorem \ref{unitriangle}. 
By definition,
every monomial $\bfu$ that does not belong to  $\hscr_{n-1} $ has $b \geq 3n$ and by induction the corresponding $\hphi$-multiplier
$\hat{\bfq}=(a',b',c')$ has  $b' \geq  3n$; furthermore, by the first step of the proof  the  $\hphi$-multiplier of
a monomial $\bfu$ not in  $\hscr_{n-1} \cup \cab_n$ is different from $\bfq_n$.

Pick a monomial $\hat{\bfv}=(i,j,k)$ in $\scrM^{\triangle}_{w,3n}= \big\{\bfv\in \mtriw:  \floor{j(\mf{v})/2}=3n\big\}$,
and assume $\hat{\bfv} \neq \vphi\left(\bfu^\triangle_{2,n-1} \right)$ if $n>0$.
Consider $\bfu=\hphi^{-1}(\hat{\bfv})=(a,b,c)$.  As
$\hat{\bfv}$ does not belong to $\kscr_{n-1} \cup \cac_n$, by induction and Step $1$ the monomial
$\bfu$ does not belong to $\hscr_{n-1} \cup \cab_n$. If $\bfu=(a,b,c)$ were in $\rectw$, then  it would have $a \leq \delta_{w,n}-2$, hence $b \geq 3n+3$ by Lemmas \ref{urrectangular} and \ref{aneworder}. But then $j =b+b' \geq 6n+3$, a contradiction. Hence $\bfu$ belongs to the triangular region.

If $j= 6n$, then the equality $j=b+b'$ together with the induction hypothesis
implies $b=3n$ and $b'=3n$.  It follows that $\bfu \in \scrR^{\triangle}_{w,3n}$, and we conclude that there are precisely $\nu^+_{3n}-1$ monomials if $n>0$
(resp. $\nu^+_{0}$  if $n=0$) in $\scrR^{\triangle}_{w,3n} \setminus \left\{\bfu^\triangle_{2,n-1}\right\}$ whose $\hphi$-multiplier
has $b'=3n$. 
If $n>0$, by the induction hypothesis (vE)$_{n-1}$,
 this accounts for all the remaining $\hphi$-multipliers
having $b'=3n$; for $n=0$ the same conclusion follows from the computation of the multiplicities of the multipliers at the beginning of the proof.

If $j = 6n+1$, the $\hphi$-multiplier $\hbfq=\hat{\bfv}-\bfu$ must have $b' \geq 3n+1$, and so
for the equality $j=b+b'$ to hold we must have $b=3n$ and $b'=3n+1$:
this means that $\bfu\in\scrR^{\triangle}_{w,3n}$, thus
the number of monomials
in $\scrR^{\triangle}_{w,3n}$ whose $\hphi$-multiplier has
 $b'=3n+1$ is $\nu^-_{3n}$.

\noindent {\bf Step 3.} The bijection $\hphi$ maps
$\rtri_{w,3n+1} \setminus \left\{\bfu^{\triangle}_{3,n}\right\}$ onto $\mtri_{w,3n+1} \setminus \left\{\vphi\left(\bfu^{\triangle}_{3,n}\right)\right\}$.
Note that $\bfu^{\triangle}_{3,n}$ is the corner monomial of the special block $\cab_n$ and has to be excluded
because all we can say about $\hphi\left(\bfu^{\triangle}_{3,n}\right)$ is that it belongs to the special block $\cac_n$,
 as does $\vphi\left(\bfu^{\triangle}_{3,n}\right)$, but the two need not coincide. 
 The
 $\vphi$-multiplier of $\bfu^{\triangle}_{3,n}$ is $\bfq_n+\bfe_1$ which, by Step 1, arises also as the $\hphi$-multiplier of exactly one
monomial in $\cab_n$  and has $b'=3n+1$. 
The number of monomials
in $\scrR^{\triangle}_{w,3n+1} \setminus \left\{\bfu^{\triangle}_{3,n}\right\}$ whose $\hphi$-multiplier
has $b'=3n+1$ (resp. $b'=3n+2$) is $\nu^+_{3n+1}-1$  (resp. $\nu^-_{3n+1}$).
After this step, all available multipliers have $b' \geq 3n+2$: if a monomial $\bfu$ has
$\hphi$-multiplier with $b' \leq 3n+1$, then $\bfu$ belongs to one of the regions already examined, that is,
$\hscr_{n-1}$, $\cab_n$, $\rtri_{w,3n}$, $\rtri_{w,3n+1}$.

\noindent{\bf Proof of Step 3.}
Pick a monomial $\hat{\bfv}=(i,j,k)$  in $\scrM^{\triangle}_{w,3n+1}$,
different from $\vphi\left(\bfu^{\triangle}_{3,n}\right)=(\delta_{w,n},6n+2,c_w+m-2n-1)$,
and  let $\bfu=\hphi^{-1}(\hat{\bfv})=(a,b,c)$. 
Note that by induction and the previous steps,
if $\bfu$ belongs to the rectangular region, then $a \leq \delta_{w,n} -2$  hence $b \geq 3n+3$,
while $b \geq 3n+1$ if $\bfu$ belongs to the triangular region; furthermore the corresponding
$\hphi$-multiplier has second entry $b' \geq 3n+1$.   As $j=b+b'$ is either $6n+2$ or $6n+3$, it follows that
$\bfu$ belongs to the triangular region.

If $j = 6n+2$, then we must have $b= 3n+1$ and $b'= 3n+1$. 
Thus, $\bfu \in \scrR^{\triangle}_{w,3n+1}$.
Note that $\bfu$ cannot be $\bfu^{\triangle}_{3,n}=(\delta_{w,n}-1,3n+1,c_w)$ because
by Step 1 the monomial $\hphi\left(\bfu^{\triangle}_{3,n}\right)$ belongs to $\cac_n$, and the only
monomial in $\cac_n$ having $j=6n+2$ is $\vphi\left(\bfu^{\triangle}_{3,n}\right)$ which we have assumed to be different from $\hat{\bfv}$.
The number of monomials $\hat{\bfv}$  in the triangular region having $j=6n+2$ is by definition  $\nu^+_{3n+1}$, thus there are $\nu^{+}_{3n+1}-1$ monomials
in $\scrR^{\triangle}_{w,3n+1} \setminus \left\{\bfu^{\triangle}_{3,n}\right\}$ whose $\hphi$-multiplier has $b'=3n+1$.
Now, the total number of monomials whose $\hphi$-multiplier has $b'=3n+1$ is $\nu^{+}_{3n+1}+\nu^-_{3n}$, and in Step 1 we have shown that in $\cab_n$ there is exactly one monomial whose $\hphi$-multiplier has $b'=3n+1$, and in Step 2 that in $\scrR^{\triangle}_{w,3n+1}$ there are $\nu^-_{3n}$ such monomials, 
thus, the remaining $\hphi$-multipliers with $b'=3n+1$ are the $\nu^{+}_{3n+1}-1$ ones we have just identified in correspondence to
monomials in $\scrM^{\triangle}_{w,3n+1}$ having $j=6n+2$ and different from $\vphi\left(\bfu^{\triangle}_{3,n}\right)$.

If $j=6n+3$, the corresponding $\hphi$-multiplier has $b'\geq 3n+2$. 
Since $j=b+b'=6n+3$, $b \geq 3n+1$ and $b' \geq 3n+2$,
the only possibility is $b=3n+1$ and $b'=3n+2$.  
Thus, the number of monomials in $\scrR^{\triangle}_{w,3n+1}$ whose $\hphi$-multiplier has $b'=3n+2$ is $\nu^-_{3n+1}$. This concludes the proof of Step 3.

\noindent {\bf Step 4.}
If $\bfu =(a,b,c)\in \rectw$ and $\delta_{w,n} -2\geq a \geq \delta_{w,n}-4$,
then  $\hphi(\bfu)=\vphi(\bfu)$ and
the corresponding $\hphi$-multiplier is
$\bfq_{n+1}$.

\noindent {\bf Proof of Step 4.}
Suppose first $a=\delta_{w,n}-2$.
Let $\hat{\bfv}=(i,j,k)=\hphi(\bfu)$ and assume by contradiction that $\hat{\bfv}$ belongs to the triangular region $\mtriw$.
Then, $i \geq \delta_{w,n}-2$  by divisibility, and by Lemma \ref{ijtriangle}
\begin{equation}\label{sn31}
j + \eta(1+j) \leq \delta+2-i \leq  6n+4
\end{equation}
On the other hand, monomials in $\scrM^{\triangle}_{w,r}$ with $r \leq 3n+1$ have already been proven to be the image under $\hphi$
of a monomial either in $\hscr_{n-1}$ or in the triangular region with $r \leq 3n+1$,
so we must have  $j \geq 2(3n+2)=6n+4$. 
This contradicts (\ref{sn31}), 
because $\eta(1+j)=2 $ if $j=6n+4$.
We conclude that $\hat{\bfv} \in \mrectw$. 
Monomials $\bfv$ in the rectangular region $\mrectw$ with $i(\bfv) \geq \delta_{w,n}-1$
belong to either the special block $\cac_n$ or to $\kscr_{n-1}$, and, by either Step 1 or induction, they are the image of
some monomial in either $\cab_n$ or in $\hscr_{n-1}$. 
It follows that $i(\hat{\bfv}) \leq \delta_{w,n} -2$. 
Since $a=\delta_{w,n}-2$,
we have $i=\delta_{w,n}-2$ and the
corresponding $\hphi$-multiplier has $a'=0$, hence is of the form $\bfq_\lambda$ for some integer $\lambda$.  
Now,
the standard argument based on Corollary \ref{unirectangle} shows that $\hphi(\bfu)=\vphi(\bfu)$ for all $\bfu =(a,b,c)\in \rectw$ with $a=\delta_{w,n}-2$.

Next, suppose  that $\bfu=(a,b,c) \in \rectw$ and $a \in \{\delta_{w,n}-3,\delta_{w,n}-4\}$.
Let us show that $\hat{\bfv}=(i,j,k)=\hphi(\bfu)$ belongs to the rectangular region $\mrectw$. Assume by contradiction that $\hat{\bfv}$ belongs to the triangular region $\mtriw$.
As $ i \geq a \geq \delta_{w,n}-4$, by Lemma \ref{ijtriangle} we have
\begin{equation}\label{sn32}
j + \eta(1+j)\leq \delta+2-i \leq  6n+6
\end{equation}
while, as above, $j \geq 6n+4$. It follows that either $j=6n+4$ or $j=6n+5$.

If $j=6n+4$,  then $i \leq \delta_{w,n}-4$ by (\ref{sn32}), hence equality  $i =\delta_{w,n}-4$ must hold, and this implies that $a=\delta_{w,n}-4$ too.
Now observe that $b=b(\bfu) \geq 3n+4$ by Lemmas \ref{urrectangular} and \ref{aneworder}, and $b'(\hat{\bfv}-\bfu) \geq 3n$ by the induction hypothesis.
We conclude that the multiplier $\hbfq=\hat{\bfv}-\bfu$ is $(0,3n,m-2n)=\bfq_n$. This is a contradiction since $\bfq_n$ is no longer available as a multiplier after Step 1, so $\hat{\bfv} \in \mrectw$.

If $j=6n+5$,  then $i \leq \delta_{w,n}-3$ by (\ref{sn32}), and this inequality together with $i \geq a$ and $i \equiv j+2 \pmod{3}$ forces
$i=\delta_{w,n}-3$. Then, $(a,b)$ is either $(\delta_{w,n}-3,3n+5)$ or $(\delta_{w,n}-4,3n+4)$, 
and the corresponding multiplier
is either $\bfq_n$ or $\bfq_n+\beone$. This contradicts the fact that $\bfq_n$ and $\bfq_n+\beone$ are no longer available as  multipliers after Step 1 and Step 2 respectively. We conclude that $\hat{\bfv} \in \mrectw$.

To summarize, we have shown that $\hat{\bfv}=\hphi(\bfu)$ belongs to $\mrectw$ for every $\bfu=(a,b,c) \in \rectw$ with $a \in \{\delta_{w,n}-3,\delta_{w,n}-4\}$.
On the other hand, monomials in $\mrectw$ with $i \geq \delta_{w,n}-2$ are no longer available. 
Then, we can use the standard argument based on Corollary \ref{unirectangle} to conclude that $\hphi(\bfu)=\vphi(\bfu)$ for every $\bfu=(a,b,c) \in \rectw$ with $a \in \{\delta_{w,n}-3,\delta_{w,n}-4\}$.

\noindent {\bf Step 5.}
The bijection $\hphi$ satisfies  $\hat{\varphi} \left(\scrR^{\triangle}_{w,3n+2}\right)=\scrM^{\triangle}_{w,3n+2}$, and
the number of monomials in $\scrR^{\triangle}_{w,3n+2} $ whose $\hphi$-multiplier
has $b'=3n+2$ (resp. $b'=3n+3$) is $\nu^+_{3n+2}$  (resp. $\nu^-_{3n+2}$). 
Furthermore,
the monomial $\bfu^\triangle_{1,n}$ belongs to $\scrR^{\triangle}_{w,3n+2}$, its $\hphi$-multiplier
is $\bfq_{n+1}$, and
$\hphi\left(\bfu^\triangle_{1,n}\right)=\vphi \left(\bfu^\triangle_{1,n}\right)$.
After this step, all available multipliers have $b' \geq 3n+3$.

\noindent {\bf Proof of Step 5.}
Pick a monomial $\bfv=(i,j,k)$ in $\scrM^{\triangle}_{w,3n+2}$, so that $j=6n+4$ or $j=6n+5$, and let $\bfu=(a,b,c)=\hphi^{-1}(\bfv)$.

Suppose first that $j= 6n+4$. 
By induction and the previous steps, if $\bfu$ belongs to the rectangular region, then $ a \leq \delta_{w,n} - 5$,
hence $b \geq 3n+5$ by Lemmas \ref{urrectangular} and \ref{aneworder}, while if $\bfu$ belongs to the triangular region, then $b \geq 3n+2$. 
On the other hand, all multipliers left
after Step 3  have $b' \geq 3n+2$. 
We conclude that $b= 3n+2$, $\bfu$ is in the triangular region and the corresponding $\hphi$-multiplier has $b'=3n+2$, hence it arises from the triangular region. 
Note that this accounts for all the $\nu^{+}_{3n+2}$ multipliers left with $b'=3n+2$.

Next, suppose that $j= 6n+5$. Then, by the same argument we conclude that $b= 3n+2$, $\bfu$ is in the triangular region and the corresponding multiplier has $b'=3n+3$.

The monomial $\bfu^{\triangle}_{1,n}=(\delta_{w,n}-3, 3n+2,c_w) $   has $a=\delta_{w,n}-3$.
By Lemma \ref{ijtriangle}, all  $\bfv=(i,j,k)\in \scrM^{\triangle}_{w,3n+2}$ have $i \leq \delta_{w,n}-3$,
and, if equality $i=\delta_{w,n}-3$ holds, then $j=6n-5$.
It follows that the $\hphi$-multiplier of $\bfu^{\triangle}_{1,n}$ has $a'=0$ and $b'=3n+3$, that is, it coincides
with the $\vphi$-multiplier $\bfq_{n+1}$ of  $\bfu^{\triangle}_{1,n}$ .

\noindent {\bf Step 6.}
If $\bfu=(a,b,c) \in \rectw$
and $a=\delta_{w,n}-5$, then $\hphi(\bfu)=\vphi(\bfu)$ and
its $\hphi$-multiplier is
$\bfq_{n+1}$.

\noindent {\bf Proof of Step 6.}
Let $\bfu=(a,b,c) \in \rectw$
and $a=\delta_{w,n}-5$. 
Then, $b \geq 3n+5$ by Lemmas \ref{urrectangular} and \ref{aneworder}. As $b \equiv a-2 \equiv 0$
$\pmod{3}$ because $a+2b+3c=w$, we actually have $b \geq 3n+6$.

Suppose by contradiction that $\hat{\bfv}=(i,j,k)=\hphi(\bfu)$
belongs to the triangular region. Then,  
$j \geq 6n+6$ by induction and the previous steps,
 $i \geq \delta_{w,n}-5$ by divisibility, and by Lemma \ref{ijtriangle} we have
$j + \eta(1+j) \leq \delta+2-i \leq  6n+7$.
We conclude that $j=6n+6$ and $i=\delta_{w,n}-5$. 
But then  $b=3n+6$ by divisibility, and the multiplier would be
$(0,3n,m-2n)=\bfq_n$, contradicting Step 1.
Thus, $\hat{\bfv} \in \mrectw$, and then the standard argument based on Corollary \ref{unirectangle} shows that $\hphi(\bfu)=\vphi(\bfu)$
and $\hbfq(\bfu)=\bfq_{n+1}$ for all $\bfu=(a,b,c) \in \rectw$ having
 $a=\delta_{w,n}-5$.

\noindent {\bf  Step 7.}
If $\bfu=(a,b,c) \in \rtriw$ satisfies  $a\geq \delta_{w,n} -5$ and $b\geq 3n+3$, 
then
 $\bfu= \bfu^{\triangle}_{2,n}=(\delta_{w,n}-5, 3n+3,c_w)$,  $\hphi(\bfu)=\vphi(\bfu)$ and $\hbfq(\bfu)=\bfq_{n+1}$.

\noindent {\bf Proof of Step 7.}
Suppose $\bfu=(a,b,c) \in \rtriw$ satisfies  $a\geq \delta_{w,n} -5$ and $b\geq 3n+3$.
Since all $\hat{\bfv}\in\mrectw$ that satisfy $i \geq \delta_{w,n}-5$
have already been treated and shown to be images under $\hphi$ of monomials of either the rectangular region
or one of the special blocks,
we see that $\hat{\bfv}=(i,j,k)=\hphi(\bfu)$ is in the triangular region.
Furthermore, all remaining multipliers  have $b' \geq 3n+3$, hence $j \geq 6n+6$. But then, as in the previous step, we conclude
$j=6n+6$ and $i=\delta_{w,n}-5$. It follows that $b=3n+3$,
$
\bfu= \bfu^{\triangle}_{2,n}=(\delta_{w,n}-5, 3n+3,c_w)$,
and $\hphi(\bfu)=\vphi(\bfu)$, that is, $\hbfq(\bfu)=\bfq_{n+1}$.

\noindent{\bf Conclusion of the proof.}
Let us check that the statements (i)-(v) follow from Steps 1-7. Statement (i) follows from Step 1, (ii) from Steps 4 and 6, (iii) from Steps 2, 3 and
5, (iv) from Step 7, (vA) from Step 1, (vB) from Step 2, (vC) from Steps 2 and 3, (vD) from Steps 3 and 5. 
Let us discuss (vE).
Suppose $n<n_{max}$. The multiplier $\bfq_{n+1}$ has total multiplicity $6(\tau_w+1) + 2$ as we pointed out in b) at the beginning of the proof. It occurs as the $\hphi$-multiplier of {\em at least} $5(\tau_w+1)+2$ monomials in $\partial \hscr_n$ by Steps 1 and 4-7, hence  it can occur {\em at most} $(\tau_w+1)$ times as  the $\hphi$-multiplier of monomials not in $\hscr_n$. Finally, the total multiplicity of
multipliers having $b'=3(n+1)$ is $6(\tau_w+1)+\nu^{+}_{3(n+1)}+\nu^{-}_{3n}$; by Steps 1 and 4-7   exactly
$5(\tau_w+1)+\nu^{-}_{3(n+2)}+1$ monomials in $\partial \hscr_n$ have $\hphi$-multiplier with $b'=3(n+1)$. This proves (vE) and concludes the proof.
\end{proof}

\noindent
\hoffset -0.6truecm
\begin{table}[h]
\caption{Outline of the proof of Theorem \ref{uniquew2} for $w \equiv 2 \pmod{3}$}
\label{outlinepfthm2}
\footnotesize{
\begin{tabular}{|c|l|l|}
\hline
&&\\
  Step & Action of $\hphi$ & Multipliers \\
&&\\
\hline
\multirow{2}{*}{Step 1}&$\hphi\left(\cab_n\right)=\cac_n$ & $\bfq_n+\mathbf{e}_1$ once \\
 & $\hphi_{\mkern 1mu \vrule height 2ex\mkern2mu \cab_n}=\hphi_{n,\ell}$ & $\bfq_n$, $\bfq_{n+1}$ with mult. $(\tau_w+1)$
\\
\hline
\multirow{4}{*}{Step 2}&
\multirow{4}{*}
{
$\hphi\left(\rtri_{w,3n} \setminus \left\{\bfu^\triangle_{2,n-1}\right\}\right)=\mtri_{w,3n} \setminus \left\{\vphi\left(\bfu^\triangle_{2,n-1}\right)\right\} $
}
& $\hat{\bfq} \neq \bfq_{\lambda}$  \\
&& $\#\{\hat{\bfq}: b' = 3n \} = \nu^+_{3n}-1$  \\
&& $\#\{\hat{\bfq}: b' = 3n+1 \} = \nu^-_{3n}$  \\
&& \\
\hline
\multirow{3}{*}{Step 3}&
\multirow{3}{*}
{
$\hphi\left(\rtri_{w,3n+1} \setminus \left\{\bfu^{\triangle}_{3,n}\right\}\right)=\mtri_{w,3n+1} \setminus \left\{\vphi\left(\bfu^{\triangle}_{3,n}\right)\right\} $
}
& $\hat{\bfq} \neq \bfq_{\lambda}$\\
&& $\#\{\hat{\bfq}: b' = 3n+1 \} = \nu^+_{3n+1}-1$  \\
&& $\#\{\hat{\bfq}: b' = 3n+2 \} = \nu^-_{3n+1}$  \\
\hline
\multirow{2}{*}{Step 4}& $\forall \bfu \in \rectw$ with $\delta_{w,n}-2 \geq a(\bfu) \geq \delta_{w,n}-4$,
& $\hat{\bfq}=\bfq_{n+1}$ \\
& $\hphi\left(\bfu\right)= \vphi\left(\bfu\right)$  & with mult. $3(\tau_w+1)$\\
\hline
\multirow{4}{*}{Step 5}&
\multirow{4}{*}{$\hphi\left(\rtri_{w,3n+2} \right) = \mtri_{w,3n+2}$}& $\hat{\bfq}\left(\bfu^{\triangle}_{1,n} \right) =\bfq_{n+1}$ \\
& & $\#\{\hat{\bfq}: b' = 3n+2 \} = \nu^+_{3n+2}$  \\
& & $\#\{\hat{\bfq}: b' = 3n+3 \} = \nu^-_{3n+2}$ \\
\hline
\multirow{2}{*}{Step 6}& $\forall \bfu \in \rectw$ with $ a(\bfu) = \delta_{w,n}-5$,
& $\hat{\bfq}=\bfq_{n+1}$ \\
& $\hphi\left(\bfu\right)= \vphi\left(\bfu\right)$  & with mult. $(\tau_w+1)$\\
\hline
\multirow{2}{*}{Step 7}&
$\rtri_{w,3n+3} \cap \left\{ \bfu: \, a(\bfu) \geq \delta_{w,n}-5\right\}= \left\{\bfu^{\triangle}_{2,n}\right\}$&
\multirow{2}{*}{$\hat{\bfq}\left(\bfu^{\triangle}_{2,n} \right) =\bfq_{n+1}$} \\
& $\hphi\left(\bfu^{\triangle}_{2,n} \right)=\vphi\left(\bfu^{\triangle}_{2,n} \right)$& \\
\hline
\end{tabular}
}
\end{table}

\subsection{Proof of Theorem \ref{uniquew2} in case $w \equiv 1 \pmod{3}$}
\label{wcongruouno}
We now assume $w \equiv 1 \pmod{3}$ and $3m-2 \leq w \leq 6m-11$.

\bp \label{table trouble w1}
Suppose $w \equiv 1 \pmod{3}$. There is no $\bfu$ in the triangular region $\rtriw$ whose $\vphi$-multiplier has $a'=0$.
We list in Table \ref{1tavolabase} the monomials $\bfu$ of the triangular region $\rtriw$ whose $\vphi$-multiplier has $a'$ equal to either $1$ or $2$,
together with their $\vphi$-multipliers $\bfq(\bfu)$ and their images $\vphi(\bfu)$ in $\mtriw$.

\noindent
\hoffset -0.6truecm
\begin{table}[h]
\caption{Here $\delta_{w,n}=\delta_w-6n$ and $c_w=2(\tau_w+1)$ }
\label{1tavolabase}
\footnotesize{
\begin{tabular}{|l|l|l|l|}
\hline
&&&
\\
  $\bfu$ & $\bfq(\bfu)$ & $\varphi(\bfu)$ & range of $n$
\\
&&&
\\
\hline
 $\bfu^{\triangle}_{5,n}=(\delta_{w,n}-1, 3n,c_w) $
 & $\bfq_{n}+\bfe_1$
 &$(\delta_{w,n}, 6n+1,c_w+m-2n-1)$
 &
 $0 \leq n \leq n_{max}$
 \\
\hline
  \multirow{2}{*}{$\bfu^{\triangle}_{6,n}=(\delta_{w,n}-3, 3n+1,c_w) $}
 & \multirow{2}{*}{$\bfq_{n}+\bfe_1$}
 &  \multirow{2}{*}{$(\delta_{w,n}-2, 6n+2,c_w+m-2n-1)$}
 &  $0 \leq n \leq n_{max}-1$ if $w \equiv 1$ (mod. $6$)
 \\
  &&&
 $0 \leq n \leq n_{max}$ if $w \equiv 4$ (mod. $6$)
\\
\hline
  \multirow{2}{*}{$\bfu^{\triangle}_{7,n}=(\delta_{w,n}-5, 3n+2,c_w) $}
 & \multirow{2}{*}{$\bfq_{n}+2\bfe_1$}
 &  \multirow{2}{*}{$(\delta_{w,n}-3, 6n+4,c_w+m-2n-2)$}
 &  $0 \leq n \leq n_{max}-1$ if $w \equiv 1$ (mod. $6$)
 \\
  &&&
 $0 \leq n \leq n_{max}$ if $w \equiv 4$ (mod. $6$)
 \\
\hline
$\bfu^{\triangle}_{8,n}=(\delta^w_{n+1}, 3n+1,c_w+1) $
 & $\bfq_{n}+2\bfe_1$
 & $(\delta_{w,n}-4, 6n+3,c_w+m-2n-1)$
 &  $0 \leq n \leq n_{max}-1$
\\
\hline
\end{tabular}
}
\end{table}
\ep

\begin{proof}
The proof is similar to that of Proposition \ref{table trouble w2 v2}.
\end{proof}

\br \label{1cornermonomial}
Observe that the monomial
$  \bfu^\triangle_{5,n}= (\delta_{w,n}\!-\!1, 3n,c_w )$
is the corner monomial $\bfu^{0}_{n,0}$ of the special block $\cab_n$. It is also one of the (at most two) monomials whose $\vphi$-multiplier is
$\bfq_n+\beone$.
\er

\begin{proof}[Proof of Theorem \ref{uniquew2} in case $w\equiv 1$]
The proof proceeds by induction on $n$ examining the behavior of $\hphi$ on
$\partial \hscr_n= \hscr_n \setminus \hscr_{n-1}$ as in the case $w \equiv 2$. By Proposition \ref{hkproperties} , for $0 \leq n \leq n_{max}$
the boundary $\partial \hscr_n$, is the disjoint union of
\begin{enumerate}
\item the special block $\cab_n$,
\item the set of monomials $\bfu=(a,b,c) \in \rectw$ satisfying $\delta_{w,n}-2 \geq a \geq  \delta_{w,n}-5$, and
\item the set of monomials $\bfu=(a,b,c) \in \rtriw$ satisfying $3n \leq b \leq 3n+2$ with the exception of the corner monomial
$\bfu^{\triangle}_{5,n}=(\delta_{w,n}-1,3n,c_w)$ that has $b=3n$, but is in the special block $\cab_n$.
\end{enumerate}

Before we can write the induction statement, we need to describe in detail the multiset of multipliers  common to $\vphi$ and $\hphi$.
If a multiplier $\bfq=(a',b',c')$ has
$a'=0$, then  $\bfq=\bfq_n=(0,3n,m-2n)$ for some integer $n$.
By Proposition \ref{table trouble w1}, the monomials $\bfq_n$ occur as $\vphi$-multipliers $\bfq(\bfu)$ only for monomials $\bfu$ of the rectangular region. Their multiplicities are as follows:
\begin{enumerate}
  \item[a)] $\bfq_0=(0,0,m)$ is the $\vphi$-multiplier $\bfq (\bfu)=\vphi(\bfu)-\bfu$ of $(t_w+1)$ monomials $\bfu$, namely of
  those monomials $\bfu$  that belong to the rectangular region and have $a=\delta_w-1$ (these monomial belong to the special block $\cab_0$);
  \item[b)]
  $\bfq_n=(0,3n,m-2n)$ for $1 \leq n \leq n_{max}$ is the $\vphi$-multiplier of  $6(t_w+1)$ monomials, all from the rectangular region:
 those $\bfu =(a,b,c) \in \rectw$ for which $a$ assumes one of the values $\delta_{w,n\!-\!1}, \delta_{w,n\!-\!1}-2,\delta_{w,n\!-\!1}-3, \delta_{w,n\!-\!1}-4, \delta_{w,n\!-\!1}-5, \delta_{w,n}-1$ so that $\lambda(\bfu)=n)$. Exactly $(t_w+1)$ of these,
namely those having $a=\delta_{w,n}-1$, belong to the boundary $\partial \hscr_n$, the rest  to $ \hscr_{n-1}$.
  \item[c)] We do not need to specify the multiplicity of the multiplier   $\bfq_{n_{max}+1}$.
\end{enumerate}

On the other end, the $\vphi$-multipliers $\bfq=(a',b',c')$ arising from the triangular region are those for which $a'>0$,
and the number of $\bfu \in \rtriw$ whose
$\vphi$-multiplier  $\bfq=(a',b',c')$ has $b'=h$ is equal to $\nu^{+}_{h}+\nu^{-}_{h\!-\!1}$.
By Corollary \ref{table trouble w1}, the monomial $\bfq_n+\bfe_1$ is the $\vphi$-multiplier of exactly two monomials,
namely the corner monomial $\bfu^{\triangle}_{5,n}=(\delta_{w,n}-1,3n,c_w)$ and $\bfu^{\triangle}_{6,n}=(\delta_{w,n}-3, 3n+1,c_w)$.

We will show by induction that, for each $0 \leq n \leq n_{max}$, the following statements hold:
\begin{enumerate}[label=(\roman*)]
\item
the restriction of $\hphi$ to the special block $\cab_n$ is, up to translation as in the statement, one of the bijections $\hphi_{\ell_0}$
 of Section \ref{spblocks};

\item
$\hphi(\bfu)=\vphi(\bfu)$ if $\bfu=(a,b,c) \in \rectw$ satisfies $\delta_{w,n}-2 \geq a \geq \delta_{w,n}-5$;

\item
$\hphi$ sends monomials $\bfu=(a,b,c) \in \rtriw$ satisfying $3n \leq b \leq 3n+2$, with the possible exception of the corner monomial
$\bfu^{\triangle}_{5,n}$, to monomials $\hat{\bfv}=(i,j,k) \in \mtriw$ having the same $r$-invariant (that is, $\floor{\frac{j}{2}}=b$);

\item
if $\bfu \in \partial \hscr_n$ and $\hat{\bfq}(\bfu)=\hphi(\bfu)-\bfu=(a',b',c')$ is the corresponding
$\hphi$-multiplier, then $3n \leq b' \leq 3(n+1)$ and
\begin{enumerate}[label=(\Alph*)]
\item
$\bfq_n=\hat{\bfq}(\bfu)$ for exactly $(\tau_w+1)$ monomials $\bfu$ in $\cab_n$;
\item
there are  $\nu^{+}_{3n}$  monomials $\bfu  \in \partial \hscr_n$ whose $\hphi$ multiplier has $b'=3n$ and $a'>0$;
\item
the $\hphi$-multipliers having $b'=3n+1$ all arise from monomials $\bfu  \in \partial \hscr_n$;
\item
the $\hphi$-multipliers having $b'=3n+2$ all arise from monomials $\bfu  \in \partial \hscr_n$;
\item
 if $n < n_{max}$, the $\hphi$-multipliers having $b'=3(n+1)$ arise from monomials $\bfu  \in \partial \hscr_n$, except
 for  $(\tau_w+1) + \nu^{+}_{3(n+1)}$ many, of which at most  $(\tau_w+1)$ are equal to $\bfq_{n+1}$.
 \end{enumerate}
\end{enumerate}
Once these statements  are proven, the Theorem follows as in the proof of Theorem \ref{uniquew2}.

We break the inductive proof of statements $(i)-(iv)$ into several steps. These steps include the initial case $n=0$. For $n>0$ we proceed by induction,  in particular, we assume that $\hphi$ maps $\hscr_{n-1}$ to $\kscr_{n-1}$,
and that all $\hphi$-multipliers $\hbfq=(a',b',c')$ of monomials not in $\hscr_{n-1}$ have $b' \geq 3n$.

\noindent{\bf Step 1.}
The bijection $\hphi$ maps the special block $\cab_n$ to the special block $\cac_n$, and  the restriction of $\hphi$ to
$\cab_n$ is, up to translation,  one of the $(\tau_w+2)$ bijections of Lemma \ref{specialblockbijectionsnew}, hence it has the same multiset of multipliers
as the restriction of $\vphi$ to $\cab_n$.

\noindent{\bf Proof of Step 1.}
Suppose $\bfu=(a,b,c) \in \cab_n$ and let $\hat{\bfv}=(i,j,k)=\hphi(\bfu)$. We claim that $\hat{\bfv}\in \cac_n$.
As $\bfu \notin \hscr_{n-1}$, by induction
$\hat{\bfv} \notin \kscr_{n-1}$, and in particular $i \leq \delta_{w,n}$. As $\bfu \in \cab_n$, we
have $\delta_{w,n} \geq a \geq \delta_{w,n}-1$ and $b \geq 3n$. By divisibility, $i \geq a$, so that $\delta_{w,n} \geq i \geq \delta_{w,n}-1$,
and to check that $\hat{\bfv}$ belongs to $\cac_n$ it remains to show $j \geq 6n+1$.
By the induction hypothesis on the multipliers,
if $\hat{q}=(a',b',c')=\hat{\bfv}-\bfu$, then $b' \geq 3n$.
Therefore $j =b+b' \geq 6n$. Thus we have to exclude the case $j=6n$.

Suppose first that $\hat{\bfv}$ is in the rectangular region. The inverse image $\bfu_{0}=(a_0,b_0,c_0)=\vphi^{-1}(\hat{\bfv})$
under $\vphi$ of $\hat{\bfv}$
is then in the rectangular region, so that  $\hat{\bfv}= \vphi(\bfu_0)=\bfu_0+\bfq_{\lambda(\bfu_0)}$, that is,
$i=a_0$ and $j=b_0+3 \lambda (\bfu_0)$. As $a_0$ is equal to either $\delta_{w,n}$ or $\delta_{w,n}-1$ and
 $\bfu_{0}$ belongs to the rectangular region,  $b_0 \geq 3n+1$ by Lemmas \ref{urrectangular} and \ref{aneworder} (with $r=3n-2$).
The value of $\lambda(\bfu_0)$ is
$n$ if $a_0=\delta_{w,n}-1$, and is $n+1$ if $a_0=\delta_{w,n}$. Thus,
$j =b_0+3 \lambda \geq 6n+1$, which means that $\hat{\bfv}$ belongs to $\cac_n$ as claimed.

Suppose next $\hat{\bfv}$ is in the triangular region and $j=6n$. Recall $i \in \{\delta_{w,n},\delta_{w,n}-1\}$.
Since $i\equiv w-2j \equiv 1$ and $\delta_{w,n} \equiv -w \equiv 2 \pmod{3}$, we must have $i=\delta_{w,n}-1$. Then,
\begin{equation}\label{1vfn}
\hat{\bfv}=(\delta_{w,n}-1,6n,c_w+m-2n)= \bfu^{0}_{n,0}+\bfq_n% = \bfu^\triangle_{4,n}+\bfq_n+\beone.
\end{equation}

Summarizing, we have shown that the image under $\hphi$ of the special block $\cab_n$
is contained in $\cac_n \cup \{ \bfu^{0}_{n,0}+\bfq_n\}$. Suppose by contradiction that
$\hphi(\cab_n) \neq  \cac_n$. Then, there is $\bfu \in \cab_n$ such that $\hphi(\bfu)=  \bfu^{0}_{n,0}+\bfq_n$.
Since $j( \bfu^{0}_{n,0}+\bfq_n)= 6n$, the monomial $\bfu$ must have $b \leq 3n$ because
of the induction hypothesis that all $\hphi$-multipliers of monomials not in $\hscr_{n-1}$  have $b' \geq 3n$.
The only monomial in $\cab_n$ that satisfies $b \leq 3n$ is the corner monomial $\bfu^{0}_{n,0}$. We conclude that the $\hphi$-multiplier of the corner monomial of $\cab_n$ is $\bfq_n$. On the other hand, the bijection $\hphi$ must send the other monomials of $\cab_n$
injectively in $\cac_n$, and by Lemma \ref{cabninjetive} the multiplier $\bfq_n$ occurs
as $\hbfq(\bfu)$ for at least  $(\tau_w+1)$ monomials in $\cab_n$, the corner monomial excluded. Adding the corner monomial, we count at least  $(\tau_w+2)$ occurrences of $\bfq_n$ as a multiplier, contradicting the induction hypothesis (ivE) that
at most $(\tau_w+1)$ monomials not in $ \hscr_{n-1}$ have  $\hphi$-multiplier equal to  $\bfq_n$ (for $n=0$, the multiplicity of the multiplier
$\bfq_0$ is exactly $(\tau_w+1)$, as  noted at the beginning of the proof).

We conclude that $\hphi (\cab_n)=\cac_n$ and, using Lemma \ref{cabninjetive} again, that the restriction of
$\hphi$ to $\cab_n$  is, up to translation,  one of the bijections $\hphi_{\ell_0}$ of Lemma \ref{specialblockbijectionsnew}.
In particular, $\bfq_n$ is the $\hphi$-multiplier of exactly $(\tau_w+1)$ monomials $\bfu$ in $\cab_n$, and by
the induction hypothesis this accounts for all the remaining
occurrences of  $\bfq_n$ as a $\hphi$-multiplier: if $\hat{\bfq} (\bfu)=\bfq_n$, then $\bfu \in \hscr_{n-1} \cup \cab_n$.

\noindent {\bf Step 2.}
The bijection $\hphi$  maps $\scrR^{\triangle}_{w,3n}\setminus \left\{\bfu^\triangle_{5,n}\right\}$ into
$\scrM^{\triangle}_{w,3n}\setminus \left\{\vphi\left(\bfu^\triangle_{5,n} \right)\right\}$.
In $\scrR^{\triangle}_{w,3n}\setminus \left\{\bfu^\triangle_{5,n}\right\}$
there are $\nu^+_{3n}$ monomials whose $\hphi$-multiplier
has $b'=3n$, and $\nu^-_{3n}-1$ whose $\hphi$-multiplier has $b'=3n+1$. 
This accounts for all remaining multipliers  with $b'=3n$, hence all monomials not in $\hscr_{n-1} \cup \cab_n \cup \scrR^{\triangle}_{w,3n}$ have $\hphi$-multiplier $\hat{\bfq}=(a',b',c')$ with $b' \geq 3n+1$.

\noindent{\bf Proof of Step 2.}
By definition,
every monomial $\bfu$ that does not belong to  $\hscr_{n-1} $ has $b \geq 3n$, and by induction the corresponding $\hphi$-multiplier
$\hat{\bfq}=(a',b',c')$ has  $b' \geq  3n$. 
Furthermore, by the first step of the proof,  the  $\hphi$-multiplier of
a monomial $\bfu$ not in  $\hscr_{n-1} \cup \cab_n$ is different from $\bfq_n$.

Let $\hat{\bfv}=(i,j,k) \in \scrM^{\triangle}_{w,3n}= \big\{\bfv\in \mtriw:  \floor{j(\mf{v})/2}=3n\big\}$,
and assume that $\hat{\bfv} \neq \vphi\left(\bfu^\triangle_{5,n} \right)=(\delta_{w,n},6n+1,c_w+m-2n-1)$.
Consider $\bfu=\hphi^{-1}(\hat{\bfv})=(a,b,c)$.  
Since
$\hat{\bfv}\notin \kscr_{n-1} \cup \cac_n$, by induction and Step $1$ the monomial
$\bfu$ does not belong to $\hscr_{n-1} \cup \cab_n$. In particular,
 $\bfu$ cannot be the corner monomial $\bfu^{\triangle}_{5,n}=(\delta_{w,n}-1,3n,c_w)$.

 If $\bfu=(a,b,c)	\in\rectw$, then  it would have $a \leq \delta_{w,n}-2$, hence $b \geq 3n+2$ by Lemmas \ref{urrectangular} and \ref{aneworder} with $r=3n+1$. But then $j =b+b' \geq 6n+2$,  contradicting the assumption
$\hat{\bfv}=(i,j,k) \in \scrM^{\triangle}_{w,3n}$. Hence, $\bfu$ belongs to the triangular region.

If $j= 6n$, then the equality $j=b+b'$ together with the induction hypothesis
implies $b=3n$ and $b'=3n$.  
It follows that $\bfu \in \scrR^{\triangle}_{w,3n}$, and we conclude there are precisely $\nu^+_{3n}$ monomials
in $\scrR^{\triangle}_{w,3n}\setminus \left\{\bfu^\triangle_{5,n}\right\}$  whose $\hphi$-multiplier
has $b'=3n$. 
If $n>0$, by the induction hypothesis (ivE)$_{n-1}$, 
this accounts for all the remaining $\hphi$-multipliers
having $b'=3n$; for $n=0$, the same conclusion follows from the computation of the multiplicities of the multipliers at the beginning of the proof.

If $j = 6n+1$, the $\hphi$-multiplier $\hat{\bfq}=\hat{\bfv}-\bfu$ must have $b' \geq 3n+1$, and so
for the equality $j=b+b'$ to hold we must have $b=3n$ and $b'=3n+1$:
this means that $\bfu\in\scrR^{\triangle}_{w,3n}$, thus
the number of monomials
in $\scrR^{\triangle}_{w,3n}\setminus \left\{\bfu^\triangle_{5,n}\right\}$ whose $\hphi$-multiplier has
 $b'=3n+1$ is $\nu^-_{3n}-1$.

\noindent {\bf Step 3.} 
The bijection $\hphi$ maps
$\rtri_{w,3n+1} $ onto $\mtri_{w,3n+1}$.  The number of monomials
in $\scrR^{\triangle}_{w,3n+1}
$ whose $\hphi$-multiplier
has $b'=3n+1$ (resp. $b'=3n+2$) is $\nu^+_{3n+1}$  (resp. $\nu^-_{3n+1}$).
After this step, all available multipliers  have $b' \geq 3n+2$: if a monomial $\bfu$ has
$\hphi$-multiplier with $b' \leq 3n+1$, then $\bfu$ belongs to one of the regions already examined, that is,
$\hscr_{n-1}$, $\cab_n$, $\rtri_{w,3n}$, $\rtri_{w,3n+1}$.

\noindent{\bf Proof of Step 3.}
Let $\hat{\bfv}=(i,j,k)\in\scrM^{\triangle}_{w,3n+1}$
and  $\bfu=\hphi^{-1}(\hat{\bfv})=(a,b,c)$.
Note that, by induction and the previous steps,
if $\bfu$ belongs to the rectangular region, then $a \leq \delta_{w,n} -2$  and $b \geq 3n+2$,
while $b \geq 3n+1$ if $\bfu$ belongs to the triangular region; furthermore, the corresponding
$\hphi$-multiplier has  $b' \geq 3n+1$.

Suppose $\hat{\bfv}$ is one of the $\nu^+_{3n+1}$
monomials $\hat{\bfv}$  in the triangular region having $j=6n+2$.
Then, we must have $b= 3n+1$ and $b'= 3n+1$, thus, $\bfu \in \scrR^{\triangle}_{w,3n+1}$.
Now the total number of monomials whose $\hphi$-multiplier has $b'=3n+1$ is $\nu^{+}_{3n+1}+\nu^-_{3n}$; in Step 1 we have shown that in $\cab_n$ there is exactly one monomial whose $\hphi$-multiplier has $b'=3n+1$, and in Step 2 that in $\scrR^{\triangle}_{w,3n}$ there are $\nu^-_{3n}-1$ such monomials, thus the remaining $\hphi$-multipliers with $b'=3n+1$ are the $\nu^{+}_{3n+1}$ ones associated to
monomials in $\scrM^{\triangle}_{w,3n+1}$ having $j=6n+2$.

As a consequence, if we now suppose $\hat{\bfv}$ has $j=6n+3$, the corresponding $\hphi$-multiplier must have $b'\geq 3n+2$. 
As $j=b+b'=6n+3$ with $b \geq 3n+1$ and $b' \geq 3n+2$,
the only possibility is $b=3n+1$ and $b'=3n+2$. 
This implies that  $\bfu$ is in the triangular region, hence belongs to $\scrR^{\triangle}_{w,3n+1}$.
Thus, in $\scrR^{\triangle}_{w,3n+1}$ there are $\nu^-_{3n+1}$ monomials whose $\hphi$-multiplier has $b'=3n+2$.
Since the cardinality of $\scrR^{\triangle}_{w,3n+1}$ is  $\nu^{+}_{3n+1}+\nu^-_{3n+1}$,
this concludes the proof of Step 3.

\noindent {\bf Step 4.}
If $\bfu =(a,b,c)\in \rectw$ and $\delta_{w,n} -2\geq a \geq \delta_{w,n}-4$,
then $\hphi(\bfu)=\vphi(\bfu)$ and
the corresponding $\hphi$-multiplier is
$\bfq_{n+1}$.

\noindent {\bf Proof of Step 4.}
Suppose first $a=\delta_{w,n}-2$.
Let $\hat{\bfv}=(i,j,k)=\hphi(\bfu)$ and assume by contradiction that $\hat{\bfv}$ belongs to the triangular region $\mtriw$.
Then, by divisibility $i \geq \delta_{w,n}-2$ and by Lemma \ref{ijtriangle}
\begin{equation}\label{1sn31}
j + \eta(j-1) \leq \delta+1-i \leq  6n+3.
\end{equation}
On the other hand, monomials in $\scrM^{\triangle}_{w,r}$ with $r \leq 3n+1$ have already been proven to be the image under $\hphi$
of a monomial either in $\hscr_{n-1}$ or in the triangular region with $r \leq 3n+1$,
so we must have  $j \geq 2(3n+2)=6n+4$. This contradicts (\ref{1sn31}).
We conclude that $\hat{\bfv} \in \mrectw$. Monomials $\bfv$ in the rectangular region $\mrectw$ with $i(\bfv) \geq \delta_{w,n}-1$
belong to either the special block $\cac_n$ or to $\kscr_{n-1}$, and by either Step 1 or induction are the image of
some monomial in either $\cab_n$ or in $\hscr_{n-1}$. 
Hence, $i(\hat{\bfv}) \leq \delta_{w,n} -2$. 
Since $a=\delta_{w,n}-2$,
it follows that $i=\delta_{w,n}-2$ and the
corresponding $\hphi$-multiplier has $a'=0$, hence is of the form $\bfq_\lambda$ for some integer $\lambda$.  
Now,
the standard argument based on Corollary \ref{unirectangle} shows that $\hphi(\bfu)=\vphi(\bfu)$ for all $\bfu =(a,b,c)\in \rectw$ with $a=\delta_{w,n}-2$.

Suppose next that $\bfu=(a,b,c) \in \rectw$ and $a \in \{\delta_{w,n}-3,\delta_{w,n}-4\}$.
Then, $b=b(\bfu) \geq 3n+3$ by Lemmas \ref{urrectangular} and \ref{aneworder}, and $b' \geq 3n+2$ by the induction hypothesis
and Step 3. Let  $\hat{\bfv}=(i,j,k)=\hphi(\bfu)$. Then, $j=b+b' \geq 6n+5$.
Assume by contradiction that $\hat{\bfv}$ belongs to the triangular region $\mtriw$.
As $ i \geq a \geq \delta_{w,n}-4$, by Lemma \ref{ijtriangle}  we get
$
j + \eta(j-1)\leq \delta^w+1-i \leq  6n+5$,
which implies $j \leq 6n+4$, a contradiction. We conclude that $\hat{\bfv} \in \mrectw$.

To summarize, we have shown that $\hat{\bfv}=\hphi(\bfu)$ belongs to $\mrectw$ for every $\bfu=(a,b,c) \in \rectw$ with $a \in \{\delta_{w,n}-3,\delta_{w,n}-4\}$.
On the other hand, monomials in $\mrectw$ with $i \geq \delta_{w,n}-2$ are no longer available. Then, we can use the standard argument based on Corollary \ref{unirectangle} to conclude that $\hphi(\bfu)=\vphi(\bfu)$ for every $\bfu=(a,b,c) \in \rectw$ with $a \in \{\delta_{w,n}-3,\delta_{w,n}-4\}$.

\noindent {\bf Step 5.}
The bijection $\hphi$ satisfies  $\hat{\varphi} \left(\scrR^{\triangle}_{w,3n+2}\right)=\scrM^{\triangle}_{w,3n+2}$, and
the number of monomials in $\scrR^{\triangle}_{w,3n+2} $ whose $\hphi$-multiplier
has $b'=3n+2$ (resp. $b'=3n+3$) is $\nu^+_{3n+2}$  (resp. $\nu^-_{3n+2}$).
After this step all available multipliers have $b' \geq 3n+3$.

\noindent {\bf Proof of Step 5.}
Let $\bfv=(i,j,k)\in\scrM^{\triangle}_{w,3n+2}$, so 
 $j=6n+4$ or $j=6n+5$, and  $\bfu=(a,b,c)=\hphi^{-1}(\bfv)$.

Suppose first that $j= 6n+4$. By induction and the previous steps, if $\bfu$ belongs to the rectangular region, then $ a \leq \delta_{w,n} - 5$,
hence $b \geq 3n+4$ by Lemmas \ref{urrectangular} and \ref{aneworder}, while if $\bfu$ belongs to the triangular region, then $b \geq 3n+2$. On the other hand, all multipliers left
after Step 3  have $b' \geq 3n+2$. We conclude that $b= 3n+2$, $\bfu$ is in the triangular region and the corresponding $\hphi$-multiplier has $b'=3n+2$, hence it arises from the triangular region. Note that this accounts for all the $\nu^{+}_{3n+2}$ multipliers left with $b'=3n+2$.

Next, suppose that $j= 6n+5$. 
By the same argument, we conclude that $b= 3n+2$, $\bfu$ is in the triangular region and the corresponding multiplier has $b'=3n+3$. Since the cardinality of $\scrR^{\triangle}_{w,3n+2}$ is  $\nu^{+}_{3n+2}+\nu^-_{3n+2}$,
this concludes the proof of Step 5.

\noindent {\bf Step 6.}
If $\bfu=(a,b,c) \in \rectw$,
and $a=\delta_{w,n}-5$, 
then, $\hphi(\bfu)=\vphi(\bfu)$ and
the corresponding $\hphi$-multiplier is
$\bfq_{n+1}$.

\noindent {\bf Proof of Step 6.}
Suppose $\bfu=(a,b,c) \in \rectw$
and $a=\delta_{w,n}-5$. 
By Lemmas \ref{urrectangular} and \ref{aneworder}, we have $b \geq 3n+4$.
Since $b \equiv a+2 \equiv 2$
$\pmod{3}$, because $a+2b+3c=w$, actually $b \geq 3n+5$ holds.
Suppose by contradiction that $\hat{\bfv}=(i,j,k)=\hphi(\bfu)$
belongs to the triangular region. 
Then, $j \geq 6n+6$ by induction and the previous steps. On the other hand,
$i \geq \delta_{w,n}-5$ by divisibility, and $j + \eta(j-1) \leq \delta^w+1-i \leq  6n+6
$ by Lemma \ref{ijtriangle}.
This implies $j \leq 6n+5$, a contradiction.
We conclude that $\hat{\bfv} \in \mrectw$, and then the standard argument based on Corollary \ref{unirectangle} shows that $\hphi(\bfu)=\vphi(\bfu)$
and $\hbfq(\bfu)=\bfq_{n+1}$ for all $\bfu=(a,b,c) \in \rectw$ having
 $a=\delta_{w,n}-5$.

\noindent{\bf Conclusion of the proof.}
Let us check that the statements (i)-(iv) follow from Steps 1-6 (see also table \ref{tableunique1}). 
Statement (i) follows from Step 1, (ii) from Steps 4 and 6, (iii) from Steps 2, 3 and 5,  (ivA) from Step 1, (ivB) from Step 2, (ivC) from Steps 2 and 3, (ivD) from Steps 3 and 5. 
We  discuss (ivE).
Suppose $n<n_{max}$. The multiplier $\bfq_{n+1}$ has total multiplicity $6(\tau_w+1) $ as we pointed out at the beginning of the proof. It occurs as the $\hphi$-multiplier of {\em at least} $5(\tau_w+1)$ monomials in $\partial \hscr_n$ by Steps 1, 4 and 6, hence  it can occur {\em at most} $(\tau_w+1)$ times as  the $\hphi$-multiplier of monomials not in $\hscr_n$. Finally, the total multiplicity of
multipliers having $b'=3(n+1)$ is $6(\tau_w+1)+\nu^{+}_{3(n+1)}+\nu^{-}_{3n+2}$; by Steps 1-6   exactly
$5(\tau_w+1)+\nu^{-}_{3n+2}$ monomials in $\partial \hscr_n$ have $\hphi$-multiplier with $b'=3(n+1)$. This proves (ivE) and concludes the proof.
\end{proof}

\noindent
\hoffset -0.6truecm
\begin{table}[h]
\caption{Outline of the proof of Theorem \ref{uniquew2} for $w \equiv 1 \pmod{3}$}
\label{tableunique1}
\footnotesize{
\begin{tabular}{|c|l|l|}
\hline
&&\\
  Step & Action of $\hphi$ & Multipliers \\
&&\\
\hline
\multirow{2}{*}{Step 1}&$\hphi\left(\cab_n\right)=\cac_n$ & $\bfq_n+\mathbf{e}_1$ once \\
 & $\hphi_{\mkern 1mu \vrule height 2ex\mkern2mu \cab_n}=\hphi_{n,\ell}$ & $\bfq_n$, $\bfq_{n+1}$ with mult. $(\tau_w+1)$
\\
\hline
\multirow{3}{*}{Step 2}&
\multirow{3}{*}
{
$\hphi\left(\rtri_{w,3n} \setminus \left\{\bfu^\triangle_{5,n}\right\}\right)=\mtri_{w,3n} \setminus \left\{\vphi\left(\bfu^\triangle_{5,n}\right)\right\} $
} & \\
&&  $\#\{\hat{\bfq}: b' = 3n \} = \nu^+_{3n}$  \\
&& $\#\{\hat{\bfq}: b' = 3n+1 \} = \nu^-_{3n}-1$  \\
\hline
\multirow{3}{*}{Step 3}&
\multirow{3}{*}
{
$\hphi\left(\rtri_{w,3n+1} \right)=\mtri_{w,3n+1} $
}
& \\
&& $\#\{\hat{\bfq}: b' = 3n+1 \} = \nu^+_{3n+1}$  \\
&& $\#\{\hat{\bfq}: b' = 3n+2 \} = \nu^-_{3n+1}$  \\
\hline
\multirow{2}{*}{Step 4}& $\forall \bfu \in \rectw$ with $\delta_{w,n}-2 \geq a(\bfu) \geq \delta_{w,n}-4$,
& $\hat{\bfq}=\bfq_{n+1}$ \\
& $\hphi\left(\bfu\right)= \vphi\left(\bfu\right)$  & with mult. $3(\tau_w+1)$\\
\hline
\multirow{4}{*}{Step 5}&
\multirow{4}{*}{$\hphi\left(\rtri_{w,3n+2} \right) = \mtri_{w,3n+2}$}&  \\
& & $\#\{\hat{\bfq}: b' = 3n+2 \} = \nu^+_{3n+2}$  \\
& & $\#\{\hat{\bfq}: b' = 3n+3 \} = \nu^-_{3n+2}$ \\
\hline
\multirow{2}{*}{Step 6}& $\forall \bfu \in \rectw$ with $ a(\bfu) = \delta_{w,n}-5$,
& $\hat{\bfq}=\bfq_{n+1}$ \\
& $\hphi\left(\bfu\right)= \vphi\left(\bfu\right)$  & with mult. $(\tau_w+1)$\\
\hline
\end{tabular}
}
\end{table}

\appendix

\section{Tables of values for the bijection of Theorem \ref{mainthm}}\label{AppendixTable}
In this appendix, we illustrate our results 
 when $m=7$.
 For each value  $18 \leq w \leq 23$, 
so as to include one example for each congruence class of the weight $w$ modulo $6$,
we collect all the relevant data of the bijection $\varphi$ of Theorem \ref{mainthm} in a table.
Specifically, the integral vector $(a,b,c)$, corresponding to the monomial $\bfu=x^ay^bz^c \in \rscr_w$,
and the vector $(i,j,k)$, corresponding to $\vphi(\bfu)=x^iy^jz^k$, are displayed on the same row. 
The order in which the monomials appear is a crucial technical point of our argument, and is illustrated by the tables.
 The ordering of monomials in the rectangular regions is introduced in Section \ref{orderrectangular}. 
 The ordering of monomials $\bfv=x^iy^jz^k$ in the triangular region is by increasing $j$ and then decreasing $i$.
When $w\equiv 1,2 \pmod{3}$, 
the monomials in the rectangular and triangular regions are %previously 
reordered according to a filtration defined in \ref{hnkndefinition} that allows an inductive proof of the main theorem.
The {\em special blocks} of  \cref{wcongruounodue} are marked in green. The monomials $\bfu^{(r)}$ of Lemma \ref{urrectangular} are marked in yellow.

\begin{figure}[ht]
    \centering
    \includegraphics[width=0.8\textwidth]{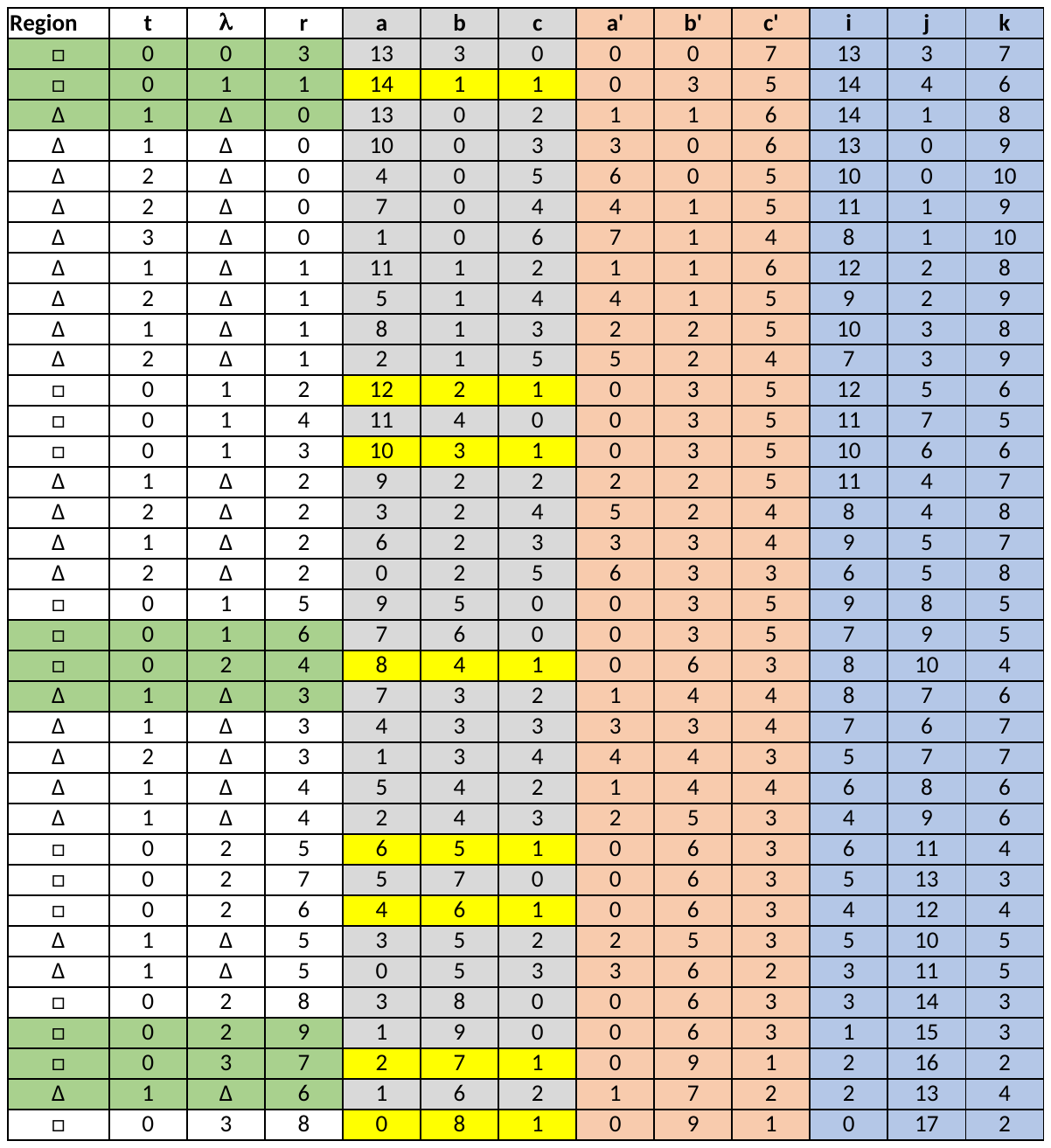}
    \caption{The map $\vphi$ for  $m=7$ and $w=19$}
    \label{fig:m7w19}
\end{figure}

\begin{figure}[ht]
    \centering
    \includegraphics[width=0.8\textwidth]{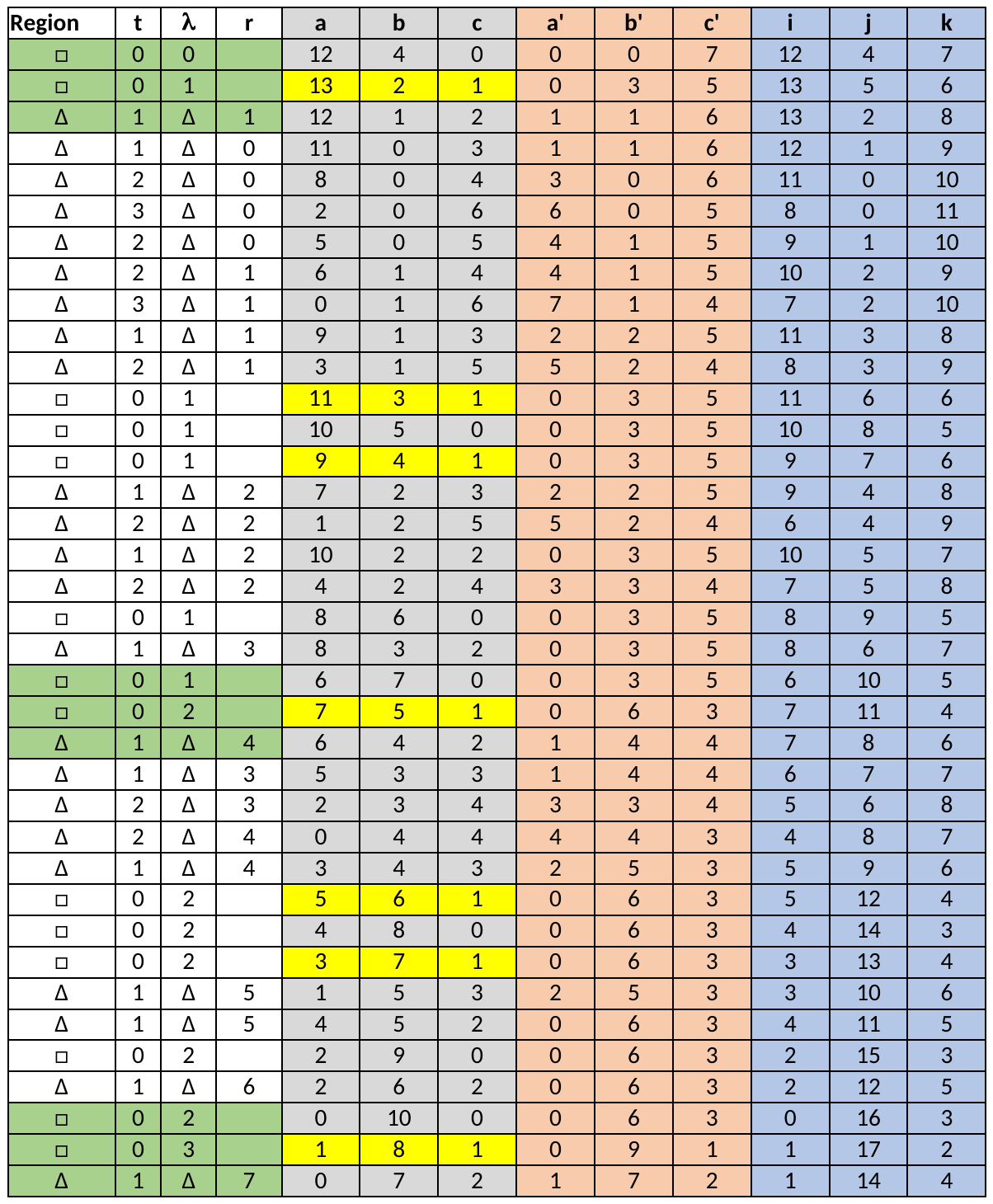}
    \caption{The map $\vphi$ for $m=7$ and $w=20$}
    \label{fig:m7w20}
\end{figure}

\begin{figure}[ht]
    \centering
    \includegraphics[width=0.8\textwidth]{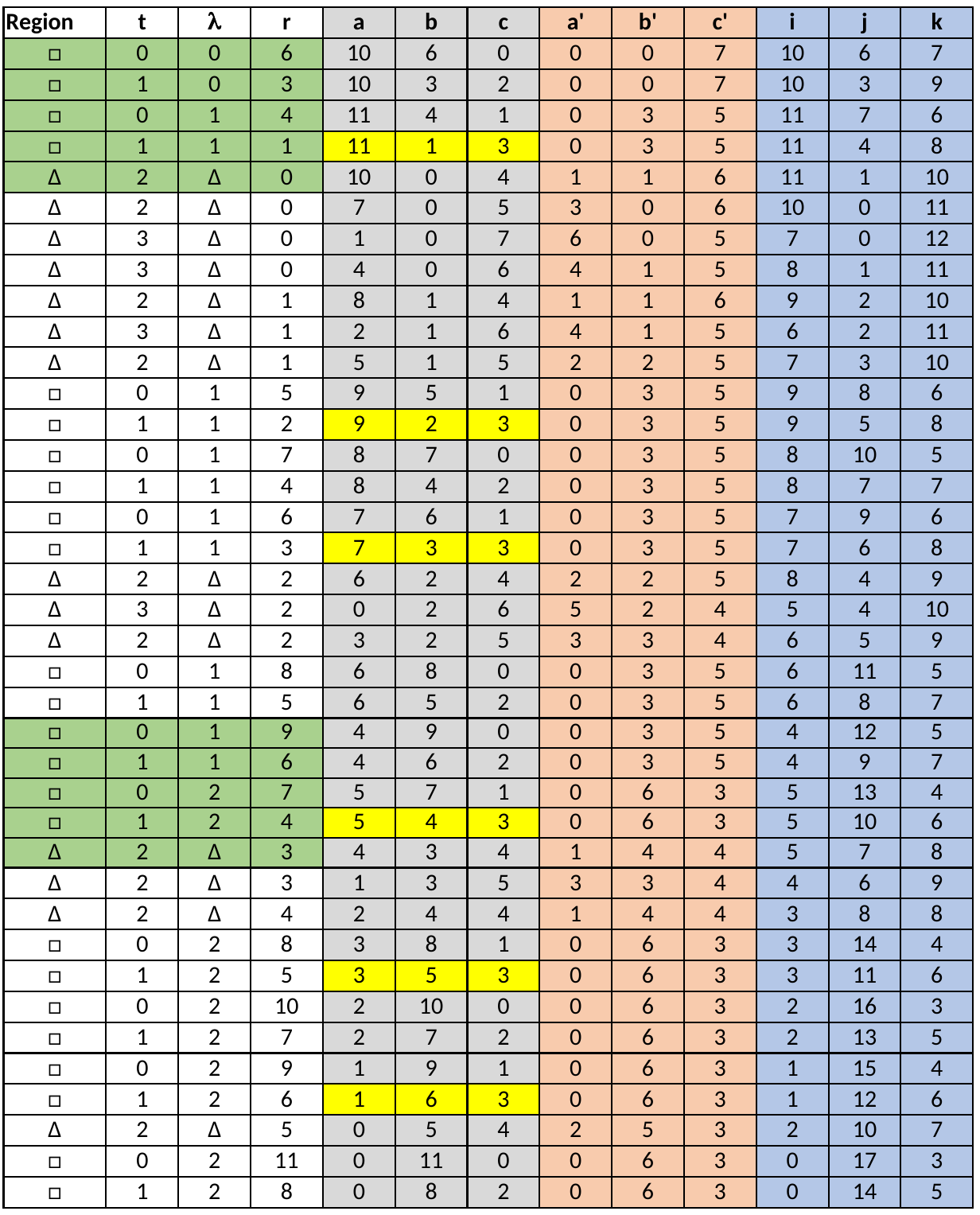}
    \caption{The map $\vphi$ for  $m=7$ and $w=22$}
    \label{fig:m7w22}
\end{figure}

\begin{figure}[ht]
    \centering
    \includegraphics[width=0.8\textwidth]{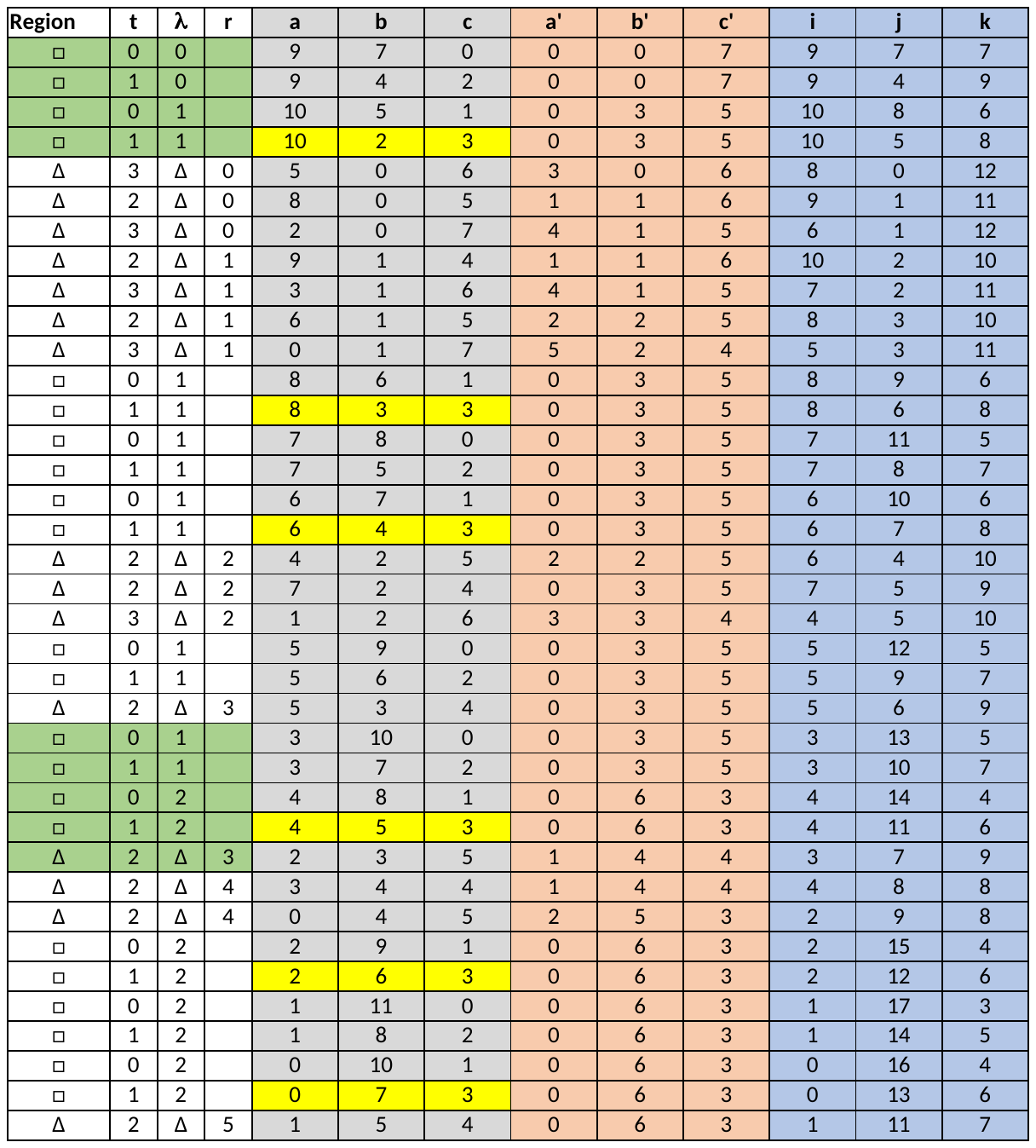}
    \caption{The map $\vphi$ for $m=7$ and $w=23$}
    \label{fig:m7w23}
\end{figure}

\printindex

\newpage

\section*{Funding}
The authors were supported by the grant PRIN 2020355B8Y
{\em Square-free Gr\"obner degenerations, special varieties and related topics},
and are members of  INdAM – GNSAGA.

\end{document}